\documentclass[reqno,12pt]{amsart}
\usepackage{amssymb}
\usepackage{amsmath}
\usepackage{amsfonts}

\newtheorem{theorem}{Theorem}
\theoremstyle{plain}

\newtheorem{corollary}{Corollary}
\newtheorem{definition}{Definition}
\newtheorem{example}{Example}
\newtheorem{lemma}{Lemma}
\newtheorem{proposition}{Proposition}
\newtheorem{remark}{Remark}
\numberwithin{equation}{section}

\input{tcilatex}

\begin{document}
\def\ss{\hskip1.11mm}

\title[Totally ordered compact sets]{Geometric properties of some totally
ordered compact sets}
\author{Mohammad DAHER and Khalil SAADI}
\email{daher.mohammad@ymail.com\\
kh\_saadi@yahoo.fr}

\begin{abstract}
In this paper, we show that there are a totally ordered compact $K$
separable, a Hausdorff topology $\tau ^{\prime }$ on $C(K)$ and two closed
subspaces $Y_{1},Y_{2}$ of ($C(K),\tau _{p})$ such that $(C(K),\tau ^{\prime
})$ is not universally measurable, $(C(K),\tau _{p})=(Y_{1},\tau _{p})\oplus
(Y_{2},\tau _{p}),$\textbf{\ } $(Y_{1},\tau _{p})$ is isomorphic to $%
(Y_{2},\tau _{p})$, $(Y_{j},\tau _{p})=(Y_{j},\tau ^{\prime })$, $j\in
\{1,2\}$ and $Bor(C(K),\tau ^{\prime })\otimes Bor(C(K),\tau ^{\prime })\neq
Bor(C(K))\times C(K)),\tau ^{\prime }\otimes \tau ^{\prime }),$ this is the
main result of this work.

We start this work to construct totally ordered non metrisable compact sets $%
K_{\mu }(E)$ from a reference set $E$ which is totally ordered, and from a
positive Borel measure on $E$ satisfying some reasonable assumptions.
\end{abstract}

\maketitle

At first, we introduce a totally ordered set~$E$, which is equipped with the
order topology~$\tau _{0}$, and a non-negative Borel measure on $(E,\tau
_{0})$. To these two elements we associate a compact set~$K$=$K_{\mu }(E)$,
obtained as spectrum of a certain $C^{\ast }$-subalgebra of $L^{\infty
}(E,\mu )$. We impose some natural properties on $E$ and $\mu $ that make
the compact $K$ not metrisable. \vskip2pt plus 1pt minus 4pt

Afterward, we introduce a certain function $h_{\mu ,E}^{{}}=h$ defined on a
subset of $K$ with values in $E$ or in $\overline{E}=E\cup \{-\infty
,+\infty \}.$ We study the surjectivity and the continuity of the function $h
$; this function plays an important role for understanding the geometric
properties of $K$. \vskip2pt plus 1pt minus 4pt

We will see that $E$ is homeomorphic to $K$, When $E$ is a separable
connected strongly Lindel\"{o}f space and if $\mu (I)>0$ for every nonempty
interval of~$E$. \vskip2pt plus 1pt minus 4pt

We study in detail a particular case when $E$ is a locally compact totally
ordered group satisfying the following conditions: 
\begin{equation}
\forall a,b\in E\text{ such that }a<b,\,\,\,\exists c\in E\text{ verifying }%
a<c<b.  \tag{*}
\end{equation}%
\begin{equation}
a\leq b\Longleftrightarrow a-b\leq 0,\quad a,b\in E.  \tag{**}
\end{equation}%
\bigskip

Our main result (theorem 6) will be shown, for this case and we will see
that the compact $K$ ($K$ is a Rosenthal compact set) has the following
property:

\bigskip $Bor(C(K),\tau _{p})\otimes Bor(C(K),\tau _{p})=Bor(C(K))\times
C(K)),\tau _{p}\otimes \tau _{p})$ (with the continuous hypothesis) but $%
Bor(K)\otimes Bor(K))\neq Bor(K\times K)$.\vskip2pt plus 1pt minus 4pt

Recall that if $L$ is a metrisable compact set, then $Bor(L)\otimes
Bor(L))=Bor(L\times L)$.\vskip2pt plus 1pt minus 4pt

We will see in part 6 (lemma 14) that if $L$ is a Rosenthal compact set ,
then $(C(L),\tau _{p})$ is universally measurable. In theorem 6 the
Hausdorff topology $\tau ^{\prime }$ on $C(K)$ satsifying $(Y_{j},\tau
_{p})=(Y_{j},\tau ^{\prime })$, $j\in \{1,2\}$ and $(C(K),\tau
_{p})=(Y_{1},\tau _{p})\oplus (Y_{2},\tau _{p}),$\textbf{\ } $(Y_{1},\tau
_{p})$ but $(C(K),\tau ^{\prime })$ is not universally measurable (here $%
Y_{1}$ and $Y_{2}$ are two closed subspaces of $(C(K),\tau _{p})).$\vskip2pt
plus 1pt minus 4pt

Let $(E, \leq)$ be a totally ordered set. The order topology $\tau_0$ on $E$
is the topology generated by the intervals of the form $\{x \in E ; x < a\}$
or of the form $\{x \in E ; x > b\}$, $a, b \in E$. This topology is
Hausdorff (if $x < a < y$ in~$E$, $V_x = \{ t \in E; t < a\}$ and $V_y = \{
t \in E; t > a\} $ are disjoint neighborhoods of $x$ and $y$; if there is no 
$a$ between $x$ and $y$, $V_x = \{ t \in E; t < y\}$ and $V_y = \{ t\in E; t
> x\}$ do the job). \smallskip

We denote by $\tau_1$ (resp. $\tau_2)$ the topology on $E$ generated by the
intervals of the form $\{x \in E ; x > a\}$ or $\{x \in E ; x \leq b\}$
(resp. of the form $\{x \in E ; x < a\}$ or $\{x \in E ; x \geq b\}$), $a, b
\in E$. \smallskip

Let $(X,\tau )$ be a topological space and $X_{1}$ a subset of~$X$. $%
(X_{1},\tau )$ (since unless is said otherwise) denotes the topological
space formed by the relative open subsets (a relative open subset of $X_{1}$
is of the form $X_{1}\cap V$, where $V$ is an open subset of $(X,\tau )$).

\section{Construction of totally ordered compact sets}

In this part we construct totally ordered compact sets $K=K_{\mu }(E)$ non
metrisables, in passing, we introduce the function $h$, this function help
us to know the topological structure of $K$.

Let $(E,\leq )$ be a totally ordered set. If $E$ admits a maximum (resp. a
minimum) we denote

$+\infty = \max \{x \in E; x \in E\}$ (resp. $-\infty = \min \{x \in E; x
\in E\})$.

Let $\mu $ be a positive Borel measure on $(E, \tau_0)$. One supposes that
the measure $\mu $ verifies the following condition: 
\begin{equation}
\mu (\{u \in E ; u > t \}) > 0, \quad \forall t \in E \setminus \{+\infty\}.
\label{y}
\end{equation}

Let $t\in E$ and let $f^{t}$ be the characteristic function of $\{u\in
E;u>t\}$. Denote by $A$ the $C^{\ast }$-subalgebra of $L^{\infty }(\mu )$
(with unit) generated by the family $\{f^{t},t\in E\}$ (the scalars are
complex scalars). Let $K(E)$ be the set of characters on~$A$. By \cite{Godem}%
, \cite{Sto}, $K_{\mu }(E)$ is $\sigma (A^{\ast },A)$ compact and $%
A=C(K_{\mu }(E))$.

As $f^s f^t = f^{s \vee t}$ when $s, t \in E$, the family $(f^t)_{t \in E}$
is stable under product, and $A$ is simply the closure in $L^\infty(\mu)$ of
the vector space generated by the constant function $\mathbf{1}$ (the unit~$%
e $ of~$A$) and the family $(f^t)_{t \in E}$. As a consequence, a character $%
\theta \in K$ is fully determined by its values on the $f^t$s: if $%
\theta(f^t) = \eta(f^t)$ for all $t \in E$, then $\theta = \eta$.

If $\theta \in K_{\mu }(E)$ is a character and $\chi _{F}$ the
characteristic function of a subset $F\subset E$, then $\theta (\chi _{F})=0$
or~$1$ because $(\chi _{F})^{2}=\chi _{F}$. In particular, $\theta
(f^{t})\in \{0,1\}$ and furthermore, when $s<t$ in $E$, we have $\theta
(f^{s})\geq \theta (f^{t})$ since $f^{s}=f^{t}+\chi _{(s,t]}$.

For every $t\in E$ let 
\begin{equation*}
\lbrack t]=\Bigl\{t^{\prime }\in E;\bigl\|f^{t}-f^{t^{\prime }}\bigr\|%
_{L^{\infty }(E,\mu )}=0\Bigr\}=\{t^{\prime }\in E;\mu (]t\wedge t^{\prime
},t\vee t^{\prime }])=0\}.
\end{equation*}%
In the following, we suppose that 
\begin{equation}
\{[t];t\in E\}\text{ is uncountable.}  \label{x}
\end{equation}%
Since $\Vert f^{t}-f^{t^{\prime }}\Vert _{L^{\infty }(E,\mu )}\geq 1$ for $%
[t]\neq \lbrack t^{\prime }]$, $C(K_{\mu }(E))=A$ is not a separable space,
hence $K_{\mu }(E)$ is not metrisable. \vskip2pt plus 1pt minus 1pt

We define the function $\psi _{\mu ,E}:K_{\mu }(E)\times E\rightarrow
\{0,1\} $ by 
\begin{equation*}
\psi _{_{\mu ,}E}(\theta ,t)=\theta (f^{t}),\quad \theta \in K,\,t\in E.
\end{equation*}%
When we have that $-\infty ,+\infty \in E$, the function $\psi _{\mu ,E}$ is
defined on $K_{\mu }(E)\times \lbrack E\setminus \{-\infty ,+\infty \}]$,
with values in $\{0,1\}$ as was said above.

Let $B_{1}^{\mu ,E}$ be the subset of all characters $\theta \in K$ such
that the map $(E,\tau _{0})\ni t\mapsto \psi _{\mu ,E}(\theta ,t)$ is
continuous; this means that the two sets $\{t\in E;\psi _{\mu ,E}(\theta
,t)=\varepsilon \}$, for $\varepsilon =0,1$, are both closed (and both
open). Let now $\theta \in K_{\mu }(E)\setminus B_{1}^{\mu ,E}$; there are
two types of discontinuities.

\emph{Case 1}, $\{t\in E;\psi _{\mu ,E}(\theta ,t)=0\}$ is not closed. Let $%
h_{\mu ,E}(\theta )$ be a point in the closure of $\{t\in E;\psi _{\mu
,E}(\theta ,t)=0\}$ in $E$ such that $\psi _{\mu ,E}(\theta ,h_{\mu
,E}^{{}}(\theta ))=1$.

\begin{lemma}
\label{hl} The point $h_{\mu ,E}^{{}}(\theta )$ is unique.
\end{lemma}

\noindent \textit{Proof}. Suppose that there exists $h^{\prime }(\theta )\in
E$ in the closure of the set $E_{0}(\theta )=\{t\in E;\psi _{\mu ,E}(\theta
,t)=0\}$, such that $\psi _{\mu ,E}(\theta ,h^{\prime }(\theta ))=1$ and $%
h_{\mu ,E}^{{}}(\theta )<h^{\prime }(\theta )$. The set $\{u\in
E;u<h^{\prime }(\theta )\}$ is then an open neighborhood of $h_{\mu
,E}^{{}}(\theta )$, but $h_{\mu ,E}^{{}}(\theta )$ is in the closure of $%
E_{0}(\theta )$, hence there is $u_{1}\in E_{0}(\theta )$ such that $%
u_{1}<h^{\prime }(\theta )$. This is impossible, because the map $t\mapsto
\psi _{\mu ,E}(\theta ,t)$ is a non-increasing function. Thus $h_{\mu
,E}^{{}}(\theta )$ is unique. $\blacksquare $ \bigskip

\emph{Case 2}, $\{t\in E;\psi _{\mu ,E}(\theta ,t)=1\}$ is not closed. We
let now $h_{\mu ,E}^{{}}(\theta )$ be the unique point in the closure of $%
\{t\in E;\psi _{\mu ,E}(\theta ,t)=1\}$ in $E$ such that $\psi _{\mu
,E}(\theta ,h_{\mu ,E}^{{}}(\theta ))=0$. \bigskip

Let $B_{2}^{\mu ,E}=\{\theta \in K_{\mu }(E)\setminus B_{1}^{\mu ,E}$ $%
;\{t\in E;\psi _{\mu ,E}(\theta ,t)=0\}\text{ is not closed}\}$ and let $%
B_{3}^{\mu ,E}=\{\theta \in K_{\mu }(E)\setminus B_{1}^{\mu ,E};$ $\{t\in
E;\psi _{\mu ,E}(\theta ,t)=1\}\text{ is not closed}\}$.

When the space $E$ is fixed, we write $K=K_{\mu }(E),$ $h_{\mu ,E}=h,$ $\psi
=\psi _{\mu ,E}$ and $B_{j}^{\mu ,E}=B_{j},$ $j\in \left\{ 1,2,3\right\} .$

\begin{lemma}
\label{kt} For every $\theta \in B_{2}^{{}}\cup B_{3}^{{}}$, we have \ \ 
\begin{equation}
h(\theta )=\sup \{t\in E;\psi (\theta ,t)=1\}.  \label{1}
\end{equation}
\end{lemma}

\noindent \textit{Proof}. Remember that the map $t\mapsto \psi (\theta ,t)$
is a non-increasing function; let $E_{j}(\theta )=\{t\in E;\psi (\theta
,t)=j\}$, for $j=0,1$.

If $\theta \in B_{2}^{{}}$, we have $\psi (\theta ,h_{{}}(\theta ))=1$ while 
$h(\theta )$ is in the closure of~$E_{0}(\theta )$ ; when $v>h(\theta )$,
the set $\{u\in E;u<v\}$ is a neighborhood of $h(\theta )$, hence there is $%
u_{1}\in E$ such that $u_{1}<v$ and $\psi (\theta ,u_{1})=0$, thus $\psi
(\theta ,v)=0$: when $\theta \in B_{2}^{{}}$, the point $h(\theta )$ is
actually the largest element of $E_{1}(\theta )$. \smallskip

If $\theta \in B_{3}^{{}}$, then $\psi (\theta ,h(\theta ))=0$; if $\psi
(\theta ,v)=1$, we have $v<h(\theta )$, showing already that $h(\theta )$ is
an upper bound of $E_{1}(\theta )$. If now we have $a<h(\theta )$, then $%
V=\{t\in E;t>a\}$ is a neighborhood of $h(\theta )$, implying that there is $%
v_{1}\in E$ such that $v_{1}>a$ and $\psi (\theta ,v_{1})=1$. This shows
that $a<h(\theta )$ cannot be an upper bound of~$E_{1}(\theta )$. $%
\blacksquare $ \vskip2pt plus 1pt minus 1pt

By an argument similar to that of Lemma~\ref{kt}, we see that:

\begin{lemma}
\label{oo} For every $\theta \in B_{2}^{{}}\cup B_{3}^{{}}$ : 
\begin{equation*}
h(\theta )=\inf \{t\in E;\psi (\theta ,t)=0\}.
\end{equation*}
\end{lemma}

Note that $B_{2}^{{}}\cap B_{3}^{{}}=\emptyset $, $B_{2}^{{}}=\{\theta \in
K\setminus B_{1}^{{}};\psi (\theta ,h(\theta ))=1\}$ and $B_{3}=\{\theta \in
K\setminus B_{1}^{{}};\psi (\theta ,h(\theta ))=0\}$.

\begin{remark}
\label{gd} Assume that $(E, \tau_0)$ is separable and that $(E, \leq)$
satisfies the condition~$(*)$. Then $(E, \tau_0)$ has a countable basis.
\end{remark}

\bigskip

\noindent\textit{Proof}. Indeed, let $(a_n)_{n \in \mathbb{N}}$ be a dense
sequence in $(E, \tau_0)$. For every $m, n \in \mathbb{N}$, consider $%
V_{m,n} = \{ x \in E; a_n > x > a_m\}$ and let $M = \{(m, n) \in \mathbb{N}%
^2 ; V_{m, n} \neq \emptyset\}$. It is easy to see that $(V_{m, n})_{(m, n)
\in M}$ forms a basis of $(E, \tau_0)$. $\blacksquare $

\begin{example}
\label{yu} Let $E$ be a separable uncountable abelian locally compact
totally ordered group \cite{Hew}. Assume that $E$ satisfies the condition~$%
(\ast )$ and that the order topology coincides with the topology of~$E$.
Choose $\mu =m$ the Haar measure on~$E$. We observe that every non empty
open set in $E$ is of strictly positive measure, and since $\mu (\{t\})=0$
for every $t\in E$, the conditions~(\ref{y}) and~(\ref{x}) are satisfied. We
shall study this example in detail in the last part.
\end{example}

\begin{example}
\label{UT}Let $E=\left[ 0,1\right] $ and $\mu =m$ the Lebesgue measure, in
corollary \ref{kv} we will see that $K=K_{m}(E)$ is homeomorphic to the
two-arrows space $\left[ 0,1\right] \times \left\{ 0,1\right\} .$
\end{example}

\begin{proposition}
\label{ah} The compact set $K$ is totally ordered (i.e., the topology on $K$
is defined by a total order relation).
\end{proposition}

\noindent \textit{Proof}. We define on $K$ the following relation, for all $%
\theta _{1},\theta _{2}\in K$: 
\begin{equation*}
\theta _{1}\leq \theta _{2}\Leftrightarrow \psi (\theta _{1},t)=\theta
_{1}(f^{t})\leq \theta _{2}(f^{t})=\psi (\theta _{2},t),\quad \forall t\in E.
\end{equation*}%
Let us show that $(K,\leq )$ is totally ordered. Pick $\theta _{1},\theta
_{2}\in K$, $\theta _{1}\neq \theta _{2}$. We may assume that there exists $%
t_{0}\in E$ such that $\psi (\theta _{1},t_{0})=0$ and $\psi (\theta
_{2},t_{0})=1$. Let $t\in E$. If $t<t_{0}$, $\psi (\theta _{2},t)=1$, hence $%
\psi (\theta _{1},t)\leq \psi (\theta _{2},t)$. If $t\geq t_{0}$, $\psi
(\theta _{1},t)=0$, hence $\psi (\theta _{1},t)\leq \psi (\theta _{2},t)$,
this implies that $\theta _{1}<\theta _{2}$. It remains to show that the
topology of~$K$ and the order topology are identical.

Since $K$ is compact and ($K,\tau _{0})$ is Hausdorff (see the
introduction), it is enough to show that the order topology is weaker than
the topology of~$K$. Given $\theta \in K$, let 
\begin{equation*}
Z=\{\alpha \in K;\alpha <\theta \}.
\end{equation*}%
Suppose that $Z$ is not empty and $\alpha \in Z$. Since $\alpha ,\theta $
are $\{0,1\}$-valued on the family $(f^{t})_{t\in E}$ and $\alpha <\theta $,
it is clear that there is $t\in E$ such that $0=\alpha (f^{t})<\theta
(f^{t})=1$, hence 
\begin{equation*}
Z=\cup _{t\in E}\{\alpha \in K;\psi (\alpha ,t)=0\text{ and }\psi (\theta
,t)=1\}.
\end{equation*}%
On the other hand, $K\ni \alpha \mapsto \alpha (f^{t})$ is continuous for
the topology of~$K$, namely, the $\sigma (A^{\ast },A)$ topology; for every $%
t\in E$ the set 
\begin{equation*}
\{\alpha \in K;\psi (\alpha ,t)=0\text{ and }\psi (\theta ,t)=1\}
\end{equation*}%
is an open subset of $K$, hence $Z$ is an open subset of $(K,\sigma (A^{\ast
},A))$. By the same proof, $Z^{\prime }=\{\alpha \in K;\alpha >\theta \}$ is
an open subset of~$K$. $\blacksquare $

\begin{definition}
\label{vo} Let $K$ be a Hausdorff compact space. We say that $C(K)$ admits
an equivalent $\tau _{p}$-Kadec norm, if there exists an equivalent norm $%
\rho $ on $C(K))$ such that the strong topology and the $\tau _{p}$ topology
coincide on $\{f\in C(K));\rho (f)=1\}$.
\end{definition}

\begin{corollary}
\label{mn} The space $C(K)$ admits an equivalent $\tau _{p}$-Kadec norm.
\end{corollary}

\noindent\textit{Proof}. This is a consequence of Proposition~\ref{ah}, and
of Theorem~A in \cite{Hay-Jay-Nam-Rog}. $\blacksquare $ \vskip2pt plus 1pt
minus 1pt

When $(X, \tau)$ is a topological space, we denote by $\mathop{\rm Bor}%
\nolimits(X, \tau)$ the Borel $\sigma$-algebra of subsets of~$X$ that is
generated by the class of open subsets of~$(X, \tau)$.

\begin{corollary}
\label{mr} We have $\mathop{\rm Bor}\nolimits(C(K),\Vert .\Vert )=%
\mathop{\rm
Bor}\nolimits(C(K)),\tau _{p})$.
\end{corollary}

\noindent \textit{Proof}. Corollary~\ref{mn} shows us that there exists on $%
C(K)$ an equivalent $\tau _{p}$-Kadec norm. By \cite{Jay-Nam-Rog}, it
follows that $\mathop{\rm Bor}\nolimits(C(K),\Vert .\Vert )=\mathop{\rm
Bor}\nolimits(C(K)),\tau _{p})$. $\blacksquare $

\section{Surjectivity of $h$}

In this part, we study the surjectivity of $h,$ We find sufficient
conditions for $h$ to be surjective.

On the set $\{ [t] ; t \in E\}$ we define an order relation by 
\begin{equation*}
[t] < [u] \text{ if } t < u \text{ and } \mu (]t, u]) > 0
\end{equation*}
(we know that $[t] = [u]$ when $\mu (]t \wedge u, t \vee u]) = 0)$. It is
obvious that the set $\{ [t] ; t \in E\}$ is totally ordered by this
relation.

\begin{definition}
\label{nj} A totally ordered set is said to be complete if every nonempty
subset that has an upper bound, has a least upper bound.
\end{definition}

\vskip2pt plus 1pt minus 1pt

One sees then that every nonempty subset $F$ that has a lower bound has a
greatest lower bound, by applying the preceding definition to the set $%
F^{\prime}$ of lower bounds of~$F$. \medskip

We denote by $\theta _{E}^{\prime }=\theta ^{\prime }$ the maximal element
of $K$ and by $\theta _{E}^{\prime \prime }=\theta ^{\prime \prime }$ the
minimal element of $K$. Note that $\psi (\theta ^{\prime },t)=1$ and $\psi
(\theta ^{\prime \prime },t)=0$, for every $t\in E\setminus \{-\infty
,+\infty \}$. It is clear that $\theta ^{\prime },\theta ^{\prime \prime
}\in B_{1}^{{}}$, the subset of $K$ consisting of those $\theta $ with $%
t\mapsto \theta (f^{t})$ continuous on $(E,\tau _{0})$.

\begin{remark}
\label{cnn} Denote by $h_{j}^{{}}$ the restriction of $h$ to $B_{j}^{{}}$, $%
j\in \{2,3\}$; then $h_{j}^{{}}:B_{j}^{{}}\rightarrow E$ is an injective
function, for $j\in \{2,3\}$.
\end{remark}

\noindent \textit{Proof}. Let us show for example that $h_{2}^{{}}$ is an
injective function. Let $\theta _{1},\theta _{2}\in B_{2}^{{}}$ such that $%
\theta _{1}<\theta _{2}$. There exists $t\in E$ such that $\psi (\theta
_{1},t)=0$ and $\psi (\theta _{2},t)=1$, it follows that $h_{2}^{{}}(\theta
_{1})<t\leq h_{2}^{{}}(\theta _{2})$, (by Lemma~\ref{kt}, we know that $%
t\leq h_{2}^{{}}(\theta _{2}))$ hence $h_{2}^{{}}(\theta _{1})\neq
h_{2}^{{}}(\theta _{2})$. $\blacksquare $ \bigskip \goodbreak

\begin{remark}
\label{lmv} Suppose that $(E,\leq )$ is complete. Let $\theta \in
B_{1}^{{}}\setminus \{\theta ^{\prime },\theta ^{\prime \prime }\}$. It is
clear that $\{t\in E;\psi (\theta ,t)=1\}$ and $\{u\in E;\psi (\theta
,u)=0\} $ are both nonempty, for otherwise $\theta $ would coincide on every 
$f^{t}$, $t\in E$, with either $\theta ^{\prime }$ or $\theta ^{\prime
\prime }$, implying that $\theta =\theta ^{\prime }$ or $\theta =\theta
^{\prime \prime }$.
\end{remark}

Suppose that $(E,\leq )$ is complete and let $\theta \in B_{1}^{{}}\setminus
\{\theta ^{\prime },\theta ^{\prime \prime }\}$; there exists $u^{\prime
}\in E$ with $\psi (\theta ,u^{\prime })=0$. Note that we have $u^{\prime
}\geq u$ for every $u\in E$ such that $\psi (\theta ,u)=1$, hence $\sup
\{t\in E;\psi (\theta ,t)=1\}$ exists by completeness. Let us extend to $%
B_{1}\setminus \{\theta ^{\prime },\theta ^{\prime \prime }\}$ the
definition of $h$ by setting $h(\theta )=\sup \{t\in E;\psi (\theta ,t)=1\}$
(now, this formula holds for every $\theta \in K\setminus \{\theta ^{\prime
},\theta ^{\prime \prime }\}$, see Lemma~\ref{kt}).

\begin{lemma}
\label{nf} Suppose that the order on $E$ is complete. Then for every $\theta
\in B_{1}^{{}}\setminus \{\theta ^{\prime },\theta ^{\prime \prime }\}$, one
has $\psi (\theta ,h(\theta ))=1$.
\end{lemma}

\noindent \textit{Proof}. Let $\theta \in B_{1}^{{}}\setminus \{\theta
^{\prime },\theta ^{\prime \prime }\}$ and set $E_{1}(\theta )=\{t\in E;\psi
(\theta ,t)=1\}$. It is easy to see that $h(\theta )=\sup E_{1}(\theta )$
belongs to the closure of $E_{1}(\theta )$ for the order topology~$\tau _{0}$%
. But $E_{1}(\theta )=\{t\in E;\theta (f^{t})=1\}$ is closed since $\theta
\in B_{1}^{{}}$ and $t\mapsto \theta (f^{t})$ is continuous. It follows that 
$h(\theta )\in E_{1}(\theta )$. $\blacksquare $ \medskip

Denote $\overline{E}=E\cup \{-\infty ,+\infty \}$ ($\overline{E}=E$, if $%
-\infty ,+\infty \in E$). It is obvious that $\overline{E}$ is totally
ordered. Denote by $e$ the unit element of~$C(K)$. When $F$ is a subset of $%
E $, denote by $\chi _{F}$ the characteristic function of~$F$. \medskip

\begin{lemma}
\label{charac}Let $t\in E\diagdown \left\{ -\infty \right\} $, then there
exists a character $\theta _{t,E}\in K$ (we denote $\theta _{t,E}=\theta
_{t},$ When $E$ is fixed) such that for every $u\in E$, we have $\theta
_{t}(f^{u})=1$ and only if $[u]\leq \lbrack t]$ (the value is~$0$ otherwise,
of course).
\end{lemma}

\vskip 2pt plus 1pt minus 4pt \vskip 2pt plus 1pt minus 4pt \goodbreak

\noindent \textit{Proof}. We shall define a function $\theta _{t}$ on~$A$ by
setting first 
\begin{equation}
\left\{ 
\begin{array}{c}
\begin{array}{c}
\theta _{t}(f^{u})=1\ \text{ if }[u]\leq \lbrack t], \\ 
\theta _{t}(f^{u})=0\ \text{ if }[u]>[t],%
\end{array}
\\ 
\theta _{t}(e)=1,%
\end{array}%
\right.  \label{xa}
\end{equation}%
then extending it to the $\mathbb{C}$-vector subspace $V$ of $L^{\infty
}(\mu )$ generated by the family $(f^{t})_{t\in E}$ and the unit $e$ of~$A$
(constantly equal to $1$ on~$E$). As $f^{u}=f^{v}$ in $L^{\infty }(\mu )$
when $[u]=[v]$, every $g\in V$ can be expressed as 
\begin{equation}
g=\sum_{j=1}^{n}\alpha _{j}f^{t_{j}}+\alpha _{0}e,  \label{fonctionG}
\end{equation}%
where $\alpha _{j}\in \mathbb{C}$, $t_{j}\in E$ are such that $[t_{i}]\neq
\lbrack t_{j}]$ for $i\neq j$, and $f^{t_{j}}\neq e$, $f^{t_{j}}\neq 0$ for
all $1\leq j\leq n$ (as $L^{\infty }$ classes; condition~(\ref{y}) implies
that $f^{t}\neq 0$ when $t\neq +\infty $). Since $(E,\leq )$ is totally
ordered, we may assume that $[t_{1}]<[t_{2}]<\ldots <[t_{n}]$; since $%
f^{t_{1}}\neq e$ in~$A\subset L^{\infty }(\mu )$ and $e-f^{t_{1}}=\chi
_{F_{0}}$, $F_{0}\subset E$, we have $\mu (F_{0})>0$; it follows that $%
g=\alpha _{0}$ on~$F_{0}$ and $|\alpha _{0}|\leq \Vert g\Vert _{L^{\infty
}(\mu )}$. Next, $[t_{1}]<[t_{2}]$ implies that $f^{t_{1}}-f^{t_{2}}=\chi
_{F_{1}}$ with $F_{1}=(t_{1},t_{2}]$, $\mu (F_{1})>0$, and $g=\alpha
_{0}+\alpha _{1}$ on~$F_{1}$; it follows that $|\alpha _{0}+\alpha _{1}|\leq
\Vert g\Vert _{L^{\infty }(\mu )}$. We proceed in this way to $f^{t_{n}}\neq
0$, with the value $\sum_{j=0}^{n}\alpha _{j}$ on the set $F_{n}=\{t\in
E;t>t_{n}\}$ of positive measure. Finally, 
\begin{equation}
\sup_{0\leq k\leq n}\,\Bigl|\sum_{j=0}^{k}\alpha _{j}\Bigr|\leq \Vert g\Vert
_{L^{\infty }(\mu )}.  \label{normeG}
\end{equation}%
It follows from~(\ref{normeG}) that $e,f^{s_{1}},\ldots ,f^{s_{m}}$ are
linearly independent in~$A$ when $[s_{i}]\neq \lbrack s_{j}]$ for $i\neq j$
and $f^{s_{j}}\neq e$, $f^{s_{j}}\neq 0$ for all $1\leq j\leq m$. Indeed,
assuming as we may that $s_{1}<s_{2}<\ldots <s_{m}$, the equality 
\begin{equation*}
\sum_{j=1}^{m}\alpha _{j}f^{s_{j}}+\alpha _{0}e=0
\end{equation*}%
implies $\sum_{j=0}^{k}\alpha _{j}=0$ for $0\leq k\leq m$, hence $\alpha
_{j}=0$ for $j=0,\ldots ,m$. \vskip2pt plus 1pt minus 4pt

It follows that we may extend $\theta _{t}$ in a natural way as a linear
form on~$V$ by setting $\theta _{t}(g)=\sum_{j=1}^{n}\alpha _{j}\theta
_{t}(f^{t_{j}})+\alpha _{0}$ whenever $g$ is expressed as in~(\ref{fonctionG}%
). Furthermore, $\theta _{t}$ is a continuous function on~$V$ equipped with
the $L^{\infty }(\mu )$ norm: indeed, pick $g=\sum_{j=1}^{n}\alpha
_{j}f^{t_{j}}+\alpha _{0}e$ in the previous form and assume $%
t_{1}<t_{2}<\ldots <t_{n}$. We shall check that $\theta _{t}(g)$ is equal to
one of the sums $\sum_{j=0}^{k}\alpha _{j}$, $k\in \{0,\ldots ,n\}$,
depending on the position of~$t$: if $[t]<[t_{1}]$, we have that $\theta
_{t}(f^{t_{j}})=0$ for $1\leq j\leq n$, hence $\theta _{t}(g)=\alpha _{0}$;
if $[t_{k}]\leq \lbrack t]<[t_{k+1}]$ and $k\in \{1,\ldots ,n-1\}$, we see
that $\theta (t)(f^{t_{j}})=1$ when $j\leq k$ and $=0$ otherwise, thus $%
\theta _{t}(g)=\alpha _{0}+\alpha _{1}+\cdots +\alpha _{k}$; finally, when $%
[t]\geq \lbrack t_{n}]$, we have $\theta _{t}(f^{t_{j}})=1$ for every~$j$
and $\theta _{t}(g)=\alpha _{0}+\cdots +\alpha _{n}$. Therefore, we have by~(%
\ref{normeG}) that 
\begin{equation*}
|\theta _{t}(g)|\leq \Vert g\Vert _{L^{\infty }(\mu )}.
\end{equation*}%
Since the set $\{f^{u};u\in E\}\cup \{e\}$ forms a total set in~$A$, $\theta
_{t}$ extends (in a unique way) to a continuous linear form on $A=C(K)$.

We check finally that $\theta (t)$ is a character. First $\theta
_{t}(f^{u}f^{u})=\theta _{t}(f^{u})=\theta _{t}(f^{u})\theta _{t}(f^{u})$
because the values are $0$ or~$1$. Let $u,v$ in~$E$ with $[u]<[v]$; we have $%
f^{u}f^{v}=f^{v}$ and we can check that $\theta _{t}(f^{u}f^{v})=\theta
_{t}(f^{u})\theta _{t}(f^{v})$ by considering the cases $[t]<[u]$ (values $%
u\mapsto 0,v\mapsto 0$), $[u]\leq \lbrack t]<[v]$ (values $u\mapsto
1,v\mapsto 0$), and $[v]\leq \lbrack t]$ (values $u\mapsto 1,v\mapsto 1$).
This implies the result, by expanding products $g_{1}g_{2}$ with $g_{1}$, $%
g_{2}$ expressed as in~(\ref{fonctionG}). It follows that $\theta _{t}\in K$%
. $\blacksquare $

Let us set $h(\theta ^{\prime })=+\infty $ and $h(\theta ^{\prime \prime
})=-\infty $.

\begin{lemma}
\label{mq} Let $t\in E\diagdown \left\{ -\infty \right\} $. Then

a) $\sup \{u\in E;u\in \lbrack t]\}=b\in \lbrack t]$ and $\{u\in E;u>b\}$ is
not closed subset of $E$, if and only if there exists a character $\theta
_{t,E}^{{}}\in B_{2}^{{}}$ (we denote $\theta _{t,E}^{{}}=\theta _{t},$ When 
$E$ is fixed) such that $h(\theta _{t})=b\in \lbrack t]$.

b) $\inf \{u\in E;u\in \lbrack t]\}=a\in \lbrack t]$ and $\{u\in E;u<a\}$ is
not closed subset of $E$, if and only if there exists $\theta _{t,E}^{\prime
}=\theta _{t}^{\prime }\in B_{3}^{{}}$ such that $h(\theta _{t}^{\prime
})=a\in \lbrack t]$.

c) Assume $(E,\leq )$ complete. Then $\sup \{u\in E;u\in \lbrack t]\}=b\in
\lbrack t]$ and $\{u\in E;u>b\}$ is closed subset of $E$, if and only if
there exists $\theta _{t,E}^{\prime \prime }=\theta _{t}^{\prime \prime }\in
B_{1}^{{}}\setminus \left\{ \theta ^{\prime \prime }\right\} $ such that $%
h(\theta _{t}^{\prime \prime })=b\in \lbrack t]$.
\end{lemma}

\vskip 2pt plus 1pt minus 4pt \goodbreak

\noindent\textit{Proof}. \vskip 2pt plus 1pt minus 4pt

\textit{a)}. Assume that $\sup \{u\in E;u\in \lbrack t]\}=b\in \lbrack t]$
and $\{u\in E;u>b\}$ is not closed subset of $E$. Consider the character $%
\theta _{t}$ given by Lemma~\ref{charac}. We show that $\theta _{t}\in
B_{2}^{0}$ and $h(\theta _{t})=b$. Indeed, let $u\in E$. If $u>b$, then $%
[u]>[b]$ because $\sup \{v\in E;v\in \lbrack t]\}=b$. If $[u]>[b]$, it is
obvious that $u>b$. We conclude that $u>b$ if and only if $[u]>[b]$, or if
and only if $[u]>[t]$, since $[b]=[t]$.

Therefore, by the formula (\ref{xa}), we have $\psi (\theta _{t},u)=0$ if
and only if $u>b$. Thus $\{u\in E;\psi (\theta _{t},u)=0\}=\{u\in E;u>b\}$
is not closed subset of $E$ by the hypothesis. We conclude that $\theta
_{t}\in B_{2}^{{}}$. Finally by Lemma~\ref{kt}, 
\begin{equation*}
h(\theta _{t})=\sup \{u\in E;\psi (\theta _{t},u)=1\}=\sup \{u\in E;u\leq
b\}=b.
\end{equation*}

Conversely, assume that there is $\theta _{t}\in B_{2}^{{}}$ such that $%
h(\theta _{t})=b\in \lbrack t]$. By Lemma~\ref{kt}, we know that $\psi
(\theta _{t},v)=1$ when $v<b$ and $\psi (\theta _{t},v)=0$ when $v>b$. Thus
the set $E_{0}(\theta _{t})=\{u\in E;\psi (\theta _{t},u)=0\}$ contains $%
\{u\in E;u>b\}$ and is contained in $\{u\in E;u\geq b\}$; but $b\notin
E_{0}(\theta _{t})$, otherwise $E_{0}(\theta _{t})$ would be closed,
contrary to the assumption $\theta _{t}\in B_{2}^{{}}$. It follows therefore
that $E_{0}(\theta _{t})=\{u\in E;u>b\}$ is not closed subset of~$E$. It
remains to show that $\sup \{u\in E;u\in \lbrack t]\}=b$.

Let $u\in \lbrack t]=[b]$. We have then $\psi (\theta _{t},u)=\psi (\theta
_{t},b)=1$, and also $u\leq h(\theta _{t})=b$ by Lemma~\ref{kt}. It follows
that $b=\sup \{u\in E;u\in \lbrack t]\}$. $\blacksquare $ \medskip

\textit{b)}. The proof of \textit{b)} is similar to that of \textit{a)},
replacing in Lemma~\ref{charac} the formula \vskip2pt plus 1pt minus 1pt %
\vskip2pt plus 1pt minus 1pt

$\left\{ 
\begin{array}{c}
\begin{array}{c}
\theta _{t}(f^{u})=1\,\text{ if }[u]\leq \lbrack t], \\ 
\theta _{t}(f^{u})=0\,\text{ if }[u]>[t],%
\end{array}
\\ 
\theta _{t}(e)=1,%
\end{array}%
\right. $ \vskip2pt plus 1pt minus 1pt \vskip2pt plus 1pt minus 1pt

\noindent by the formula 
\begin{equation}
\left\{ 
\begin{array}{c}
\begin{array}{c}
\theta _{t}^{\prime }(f^{u})=1\,\text{ if }[u]<[t], \\ 
\theta _{t}^{\prime }(f^{u})=0\,\text{ if }[u]\geq \lbrack t],%
\end{array}
\\ 
\theta _{t}^{\prime }(e)=1.%
\end{array}%
\right.  \label{nt}
\end{equation}%
and replacing Lemma~\ref{kt} by Lemma~\ref{oo}. $\blacksquare $ \vskip2pt
plus 1pt minus 1pt \vskip2pt plus 1pt minus 1pt

\textit{c)}. Suppose that $\sup \{u\in E;u\in \lbrack t]\}=b\in \lbrack t]$
and $\{u\in E;u>b\}$ is a closed subset

\emph{Case 1}, t=+$\infty $ $($if $-\infty ,+\infty \in E):$ Note that by (%
\ref{y}) $\left[ t\right] =\left\{ +\infty \right\} $. Consider $\theta
_{t}^{\prime \prime }=\theta ^{\prime },$ it is clear that $\theta
_{t}^{\prime \prime }\in B_{1}^{{}}$ and $h(\theta ^{\prime \prime
}(t))=b=+\infty .$

\emph{Case 2}, t\TEXTsymbol{<}+$\infty :$ First note that $b\neq -\infty .$
Consider the character $\theta _{t}$ given by Lemma~\ref{charac}. Let us
show that $\theta _{t}\in B_{1}^{{}}$. Note that $u>b$ if and only if $%
[u]>[b]$, because $\sup \{v\in E;v\in \lbrack t]\}=b$. It implies that $%
\{u\in E;\psi (\theta _{t},u)=0\}=\{u\in E;[u]>[b]\}=\{u\in E;u>b\}$; this
is a closed subset of $E$ by hypothesis, as well as $\{u\in E;u\leq b\}$
(obviously), and the latter set is $\{u\in E;\psi (\theta _{t},u)=1\}$; it
follows that $\theta _{t}\in B_{1}^{{}}$. By Lemma~\ref{nf}, 
\begin{equation*}
h(\theta _{t})=\sup \{u\in E;\psi (\theta _{t},u)=1\}=\sup \{u\in E;u\leq
b\}=b.
\end{equation*}

We put $\theta _{t}^{\prime \prime }=\theta _{t}.$ It is obvious that $%
\theta _{t}^{\prime \prime }\neq \theta ^{\prime \prime }.$

Conversely, assume that there exists $\theta _{t}^{\prime \prime }\in
B_{1}^{{}}\setminus \left\{ \theta ^{\prime \prime }\right\} $ such that $%
h(\theta _{t}^{\prime \prime })=b\in \lbrack t]$.

\emph{Case 1, }$\theta _{t}^{\prime \prime }=\theta ^{\prime }:$ We note
that b=+$\infty ,$ then $\left\{ t\in E;t>b\right\} =\emptyset ,$ which is
closed subset of $E.$ On the other hand $\sup \{v\in E;v\in \lbrack
t]\}=b=+\infty .$

\emph{Case 2, }$\theta _{t}^{\prime \prime }\neq \theta ^{\prime }:$ Since $%
\psi (\theta _{t}^{\prime \prime },h(\theta _{t}^{\prime \prime }))=1$ by
Lemma~\ref{nf}, it follows that $\{u\in E;\psi (\theta _{t}^{\prime \prime
},u)=0\}=[\psi (\theta _{t}^{\prime \prime },.)]^{-1}\{0\}=\{u\in E;u>b\}$
is a closed subset of $E$, because $\theta _{t}^{\prime \prime }\in
B_{1}^{{}}$. An argument similar to that of \textit{a)} shows that $\sup
\{u\in E;u\in \lbrack t]\}=b$. $\blacksquare $ \vskip2pt plus 1pt minus 1pt

\begin{corollary}
\label{ds} Assume that the space $(E,\leq )$ is complete and that we have $%
\mu (]t\wedge t^{\prime },t\vee t^{\prime }])>0$ when $t\neq t^{\prime }$.
Then $h_{\mid _{B_{1}^{{}}\cup B_{2}^{{}}}}:B_{1}^{{}}\cup
B_{2}^{{}}\rightarrow \overline{E}$ is onto ($-\infty ,+\infty \notin E$).
\end{corollary}

\noindent \textit{Proof}. Note that the hypothesis implies that $[t]=\{t\}$
for every $t\in E$. Let $a\in E$. If $\{x\in E;x>a\}$ is not closed subset
of $E$, by Lemma~\ref{mq}-a), there exists $\theta \in B_{2}^{{}}$ such that 
$h(\theta )=a$. If $\{x\in E;x>a\}$ is closed subset of $E$, by Lemma~\ref%
{mq}-c), there exists $\theta \in B_{1}^{{}}$ such that $h(\theta )=a$. The
points $-\infty ,+\infty $ are images of $\theta ^{\prime },\theta ^{\prime
\prime }$. Thus $h_{\mid _{B_{1}^{{}}\cup B_{2}^{{}}}}$ is surjective. $%
\blacksquare $

\begin{proposition}
\label{omq} Assume $E$ satisfies $(\ast )$ and $\mu (]t\wedge t^{\prime
},t\vee t^{\prime }])>0$ for $t\neq t^{\prime }\in E$. Then $%
h_{j}^{{}}:B_{j}^{{}}\rightarrow E\setminus \{-\infty ,+\infty \}$ is onto, $%
j\in \{2,3\}$.
\end{proposition}

\noindent \textit{Proof}. It suffices to show by Lemma~\ref{mq} that the
sets $\{u\in E;u>t\}$ and $\{u\in E;u<t\}$ are not closed subsets of $E$,
for every point $t\in E\setminus \{-\infty ,+\infty \}$. Indeed, let $a,b\in
E$ be such that $t\in \,]a,b[$. Since $E$ satisfies the condition~$(\ast )$,
there exists $c\in E$ such that $a<c<t<b$, so $c\in \{u\in E;u<t\}$, thus $t$
is in the closure of $\{u\in E;u<t\}$. We conclude that $\{u\in E;u<t\}$ is
not closed subset in $E$. In the same way we can see that $\{u\in E;u>t\}$
is not a closed subset of $E$. $\blacksquare $

\section{Continuity of $h$}

In this part, we show that the function $h$ is continuous on the subset of $%
K $ where $h$ is defined.

Recall that $h(\theta ^{\prime })=+\infty $ and $h(\theta ^{\prime \prime
})=-\infty $. Let $K_{1}=(K\setminus B_{1}^{{}})\cup \{\theta ^{\prime
},\theta ^{\prime \prime }\}$. We denote again by $h$ the restriction of $h$
to $K_{1}$.

\begin{proposition}
\label{cc} Assume that $-\infty ,+\infty \notin E$. Then the map $%
h:(K_{1},\tau _{0})\rightarrow (\overline{E},\tau _{0})$ is continuous.
\end{proposition}

\noindent \textit{Proof}. The topology of $\overline{E}$ is generated by the
intervals of the form $\{u\in \overline{E};u<\alpha \}$ or of the form $%
\{u\in \overline{E};u>\alpha \}$, $\alpha \in E$. Let $\alpha \in E$ and let 
$\theta _{0}$ belong to the set $W=\{\theta \in K_{1};h(\theta )<\alpha \}$.
Let us show that there exists an open neighborhood of $\theta _{0}$
contained in~$W$. Note that $\theta _{0}\neq \theta ^{\prime }$.

By Lemma~\ref{oo}, we know when $\theta _{0}\neq \theta ^{\prime \prime }$
that $\psi (\theta _{0},\alpha )=0$ since $h(\theta _{0})<\alpha $ (and it
is also true for $\theta _{0}=\theta ^{\prime \prime }$), and we know that
if $\psi (\theta ,\alpha )=0$, then $h(\theta )\leq \alpha $. The set 
\begin{equation*}
U=\{\theta \in K_{1};\psi (\theta ,\alpha )=0\}=K_{1}\cap \{\theta \in
K;\theta (f^{\alpha })=0\}
\end{equation*}%
is open in~$K_{1}$, and $U$ is contained in $\{\theta \in K_{1};h(\theta
)\leq \alpha \}$. Consider $V=U\setminus \{\theta \in K_{1};h(\theta
)=\alpha \}$. Note that $V$ is an open set in $K_{1}$, because $\{\theta \in
K_{1};h(\theta )=\alpha \}$ is finite by Remark~\ref{cnn}, thus $V$ is an
open neighborhood of $\theta _{0}$ contained in~$W$.

Let $\alpha \in E$ and let $\theta _{0}$ belong to $W^{\prime }=\{\theta \in
K_{1};h(\theta )>\alpha \}$. Using now Lemma~\ref{nf}, we can show in a
similar way that the set 
\begin{equation*}
V^{\prime }=\{\theta \in K_{1};\psi (\theta ,\alpha )=1\}\setminus \{\theta
\in K_{1};h(\theta )=\alpha \}
\end{equation*}%
is an open neighborhood of $\theta _{0}$ and $V^{\prime }\subset W^{\prime }$%
. $\blacksquare $

\begin{remark}
\label{ff} The space $(E,\tau _{0})$ is connected if and only if $\ $($%
E,\leq )$ is a complete and $E$ satisfies the condition $(\ast )$ \cite[%
rem.~(d), p.~58]{Kell}.
\end{remark}

\begin{corollary}
\label{rq} Assume that $E$ is connected and $-\infty ,+\infty \notin E$.
Then $h:K\rightarrow (\overline{E},\tau _{0})$ is continuous.
\end{corollary}

\noindent \textit{Proof}. It suffices to see that $B_{1}^{{}}=\{\theta
^{\prime },\theta ^{\prime \prime }\},$ because $E$ is connected. $%
\blacksquare $ \vskip2pt plus 1pt minus 1pt

By an argument similar to that of Proposition~\ref{cc} and its corollary, we
can show the following result:

\begin{proposition}
\label{fd} Assume that $-\infty ,+\infty \in E$. Then $h:(K_{1},\tau
_{0})\rightarrow (E,\tau _{0})$ is continuous.
\end{proposition}

\begin{corollary}
\label{em} Assume that $E$ is connected and $-\infty ,+\infty \in E$. Then $%
h:K\rightarrow (E,\tau _{0})$ is continuous.
\end{corollary}

\begin{proposition}
\label{xh} Assume that $(E,\leq )$ is complete and $-\infty ,+\infty \notin
E $. Then $h:K\rightarrow \overline{E}$ is continuous.
\end{proposition}

\noindent \textit{Proof}. Let $\alpha \in E$ and let $\theta _{0}\in
W=h^{-1}(\{t\in \overline{E};t<\alpha \})=\{\theta \in K;h(\theta )<\alpha
\} $. We shall find an open neighborhood of $\theta _{0}$ in~$W$.

\emph{Case 1,} $\theta _{0}=\theta ^{\prime \prime }$: Consider $t_{0}\in E$
such that $t_{0}<\alpha $ and let $V=\{\theta \in K;\psi (\theta ,t_{0})=0\}$%
. Note that $\theta ^{\prime \prime }\in V$. It is obvious that $V\subset W$.

\emph{Case 2,} $\theta _{0}\in B_{1}^{{}}\setminus \{\theta ^{\prime \prime
}\}$: The set 
\begin{equation*}
\{u\in E;u<\alpha \}\cap \{t\in E;\psi (\theta _{0},t)=1\}
\end{equation*}%
is an open neighborhood of $h(\theta _{0})$, and because $\psi (\theta
_{0},h(\theta _{0}))=1$ by Lemma~\ref{nf}, it follows that there exists $%
a,b\in E$ ($a<b$) such that $h(\theta _{0})\in ]a,b[\subset \{u;u<\alpha
\}\cap \{t\in E;\psi (\theta _{0},t)=1\}$. Denote $V=\{\theta \in K;\psi
(\theta ,b)=0\}$; the set $V$ is an open neighborhood of $\theta _{0}$,
because $h(\theta _{0})<b$ and $\psi (\theta _{0},h(\theta _{0}))=1$. It
remains to show that $V\subset W$. Indeed, let $\theta \in V=\{\theta \in
K;\psi (\theta ,b)=0\}$ (hence $h(\theta )\leq b$). Let us show that $%
h(\theta )<b$. Assume that $h(\theta )=b$. Since $h(\theta _{0})<h(\theta
)=b $, $\theta _{0}<\theta $, there exists $\gamma \in E$ such that $\psi
(\theta _{0},\gamma )=0$ and $\psi (\theta ,\gamma )=1$, this implies that $%
a<h(\theta _{0})<\gamma <h(\theta )=b$ ($\gamma <h(\theta )$ because $\psi
(\theta ,b)=0$ and $\psi (\theta ,\gamma )=1$). We deduce that $\gamma \in
]a,b[\subset \{u\in E;u<\alpha \}\cap \{t\in E;\psi (\theta _{0},t)=1\}$,
which is impossible (note that $\psi (\theta _{0},\gamma )=0$), thus $%
h(\theta )<b$. \bigskip

Now let us show that $b\leq \alpha $. Assume that $b>\alpha $. Since $%
a<h(\theta _{0})<\alpha $, $\alpha \in \,]a,b[$, hence $\alpha \in \{u\in
E;u<\alpha \}$, this is impossible. It follows that $h(\theta )<b\leq \alpha 
$, \textit{i.e.} $\theta \in W=\{\varphi \in K;h(\varphi )<\alpha \})$. We
treat the case $\theta _{0}\in B_{2}^{{}}$ and the case $\theta _{0}\in
B_{3}^{{}}$, as in Proposition~\ref{cc}. $\blacksquare $

By an argument similar to that of Proposition~\ref{xh} we show:

\begin{proposition}
\label{fg} Assume that $(E,\leq )$ is complete and $-\infty ,+\infty \in E$.
Then $h:K\rightarrow E$ is continuous.
\end{proposition}

\bigskip

We define the measure $\mu ^{\prime }$ on $E$ by $\mu ^{\prime }(\{t\})=1$,
for every $t\in E$ ($\mu ^{\prime }$ is the \emph{counting measure}). Let~$%
A^{\prime }$ be the $C^{\ast }$-subalgebra generated by $\{f^{t},t\in E\}$
in $L^{\infty }(E,\mu ^{\prime })=\ell ^{\infty }(E)$ and let $K_{\mu
}^{\prime }(E)=K^{\prime }$ be the set of characters on $A^{\prime }$.
Denote 
\begin{equation*}
h^{\prime }(\theta )=\sup \{t\in E;\psi ^{\prime }(\theta ,t)=1\}
\end{equation*}%
where $\psi ^{\prime }(\theta ,t)=\theta (f^{t})$, $\theta \in K^{\prime }$, 
$t\in E\setminus \{-\infty ,+\infty \}$.

\section{Geometric properties of the compact sets $K$}

In this part, we shall see that the geometry of $K$ is related to the
geometry of $E.$

Recall that the lexicographic order on $E \times \{0, 1\}$ is defined in the
following way: for all $(r, i)$, $(s, j) \in E \times \{0, 1\}$,

$(r, i) < (s, j)$ if and only if $r < s$, or $r = s$ and $i < j$.

\noindent Note that if $E$ is totally ordered, then $E \times \{0, 1\} $ has
the same property. Recall that $\tau_1$ (resp. $\tau_2$) denotes the
topology on a totally ordered set generated by semi-open intervals $(a, b]$
(resp. $[a, b)$). We need the following classical lemmas:

\begin{lemma}
\label{hi} Let $E^{\prime }$ be a subset of $E$. Then $(E^{\prime },\tau
_{2})$ is homeomorphic to $(E^{\prime }\times \{1\},\tau _{0})$, viewed as
topological subspace of $(E\times \{0,1\},\tau _{0}^{{}})$ equipped with the
order topology $\tau _{0}^{{}}$ of the lexicographic order.
\end{lemma}

\begin{lemma}
\label{cw} Assume that $(E, \tau_0)$ is a strongly Lindel\"{o}f space. Then:

1) $(E, \tau_j)$ is a strongly Lindel\"{o}f space, $j \in \{1, 2\}$.

2) $\mathop{\rm Bor}\nolimits(E, \tau_0) = \mathop{\rm Bor}\nolimits(E,
\tau_j)$, $j \in \{1, 2\} $.

3) If $(E, \tau_0)$ is separable, then $(E, \tau_j)$ is separable, $j \in
\{1, 2\}$.

4) If every subset of $(E,\tau _{0})$ is separable, then every subset of $%
(E,\tau _{j})$ is separable, $j\in \{1,2\}$.
\end{lemma}

\begin{remark}
\label{nl} Let $(X,\tau )$ be a strongly Lindel\"{o}f space and let $Z$ be a
subspace of $(X,\tau )$. Then $(Z,\tau )$ is a strongly Lindel\"{o}f space.
\end{remark}

\noindent \textit{Proof of Lemma~\ref{hi}}. Consider $F:(E^{\prime },\tau
_{2})\rightarrow (E^{\prime }\times \{1\},\tau _{0})$ the map defined by $%
F(t)=(t,1)$, $t\in E^{\prime }$. We show that $F$ is a homeomorphism. Let $%
(a,j_{0})\in E\times \{0,1\}$ and let $V=\{(t,1)\in E^{\prime }\times
\{1\};(t,1)<(a,j_{0})\}$ be an open subset of $(E^{\prime }\times \{1\},\tau
_{0})$. We have 
\begin{equation*}
F^{-1}(V)=\{t\in E^{\prime };(t,1)<(a,j_{0})\}=\{t\in E^{\prime };t<a\}\in
(E^{\prime },\tau _{2}).
\end{equation*}%
Let $(a,j_{0})\in E\times \{0,1\}$ and let $W=\{(t,1)\in E^{\prime }\times
\{1\};(t,1)>(a,j_{0})\}$.

\emph{Case 1, } $j_{0}=1$: It is obvious that $F^{-1}(V)=\{t\in E^{\prime
};t>a\}\in (E_{1},\tau _{2})$.

\emph{Case 2, } $j_{0}=0$: We observe that $F^{-1}(V)=\{t\in E^{\prime
};t\geq a\}\in (E^{\prime },\tau _{2})$.

We conclude that if $V$ is an open subset of $(E^{\prime }\times \{1\},\tau
_{0})$, then $F^{-1}(V)$ is an open set of $(E^{\prime },\tau _{2})$ and
conversely. $\blacksquare $

\begin{remark}
\label{hs} By an argument similar to that of Lemma~\ref{hi}, we show that $%
(E^{\prime },\tau _{1})$ is homeomorphic to $(E^{\prime }\times \{0\},\tau
_{0})$.
\end{remark}

\noindent\textit{Proof of Lemma~\ref{cw}}.

1). We shall prove Lemma~\ref{cw} for $j=1$, the case $j=2$ can be shown by
a similar argument. The argument that allows us to show $1)$ is similar to 
\cite[p.~58, 59]{Kell}. \smallskip

Indeed, let $\{]a_i, b_i] ; a_i < b_i, i \in I\}$ be a family of $\tau_1$%
-open intervals in~$E$, let $U = \underset{i \in I}{\cup} ]a_i, b_i]$ and $V
= \underset{i \in I}{\cup} ]a_i, b_i[$. We will show that there exists a
countable subset $I^{\prime}$ of $I$ such that $U = \underset{i \in
I^{\prime}}{\cup} ]a_i, b_i]$. Denote by $J$ the subset of $I$ consisting of
those $j$ such that $b_j \notin V$.

Let us show that $\{b_j ; j \in J\}$ is countable. Let $b_j \neq
b_{j^{\prime}}$, $j, j^{\prime}\in J$. Assume that there exists $t \in \,
]a_j, b_j[ \, \cap \, ]a_{j^{\prime}}, b_{j^{\prime}}[$. If $b_j <
b_{j^{\prime}}$, then $a_{j^{\prime}} < t < b_j < b_{j^{\prime}}$, it
follows that $b_j \in \, ]a_{j^{\prime}}, b_{j^{\prime}}[ \subset V$, which
is impossible, hence $]a_j, b_j[ \, \cap \, ]a_{j^{\prime}}, b_{j^{\prime}}[
\, = \emptyset$, if $b_j \neq b_{j^{\prime}}$. By hypothesis $(E, \tau_0)$
is a strongly Lindel\"{o}f space, then there exists a countable subset $H$
of $J$, such that $\underset{i \in J}{\cup} ]a_i, b_i[ = \underset{i \in H}{%
\cup} ]a_i, b_i[$. On the other hand, if $b_j \neq b_{j^{\prime}}$, $j,
j^{\prime}\in J$, then $]a_j, b_j[ \, \cap \, ]a_{j^{\prime}},
b_{j^{\prime}}[$ is empty; we conclude that $\{b_j ; j \in J\} = \{b_j ; j
\in H\}$, \textit{i.e.} $\{b_j ; j \in J\}$ is countable.

Since $(E, \tau_0)$ is a strongly Lindel\"{o}f space, there exists a
countable subset $I_1$ of $I$, such that $V = \underset{i \in I_1}{\cup}
]a_i, b_i[$. Let $I_2 = I_1 \cup H$; this set is a countable subset of~$I$.
Let 
\begin{equation*}
W = \bigcup_{i \in I_2} \, ]a_i, b_i].
\end{equation*}
We check that $W = U$, this will end the proof. It is clear that $W \supset
V $ since $I_2 \supset I_1$. Furthermore, if $t \in U \setminus V$, then $t$
must be an endpoint $b_i$ of some $(a_i, b_i]$, for an $i \in I$. This means
that $i \in J$ and $b_i = t$ belongs to the family $\{b_j ; j \in H\}$,
hence to $\cup_{j \in H} \, ]a_j, b_j] \subset W$. The space $(E, \tau_1)$
is therefore a strongly Lindel\"{o}f space. $\blacksquare $ \vskip2pt plus
1pt minus 1pt

2). It is clear by 1) that $\mathop{\rm Bor}\nolimits(E, \tau_0) = %
\mathop{\rm Bor}\nolimits (E, \tau_1)$. $\blacksquare $ \vskip2pt plus 1pt
minus 1pt

3). Assume $(E, \tau_0)$ is separable. Let $(x_n)_{n \geq 0}$ be a dense
sequence in $(E, \tau_0)$ and $M = \{b \in E ; \, \exists a_b\in E \text{
such that } ]a_b, b] = \{b\} \}$.

\emph{Step 1: }Let us show that $M$ is countable.

By $1)$, $(E,\tau _{1})$ is a strongly Lindel\"{o}f space, it follows that $%
\underset{b\in M}{\cup }]a_{b},b]=\underset{b\in M}{\cup }\{b\}=M$ is
countable. $\blacksquare $

\emph{Step 2}: Let us show that $M_{2}=\{x_{n};n\in \mathbb{N}\}\cup M$ is
dense in $(E,\tau _{1})$.

Consider $]a,b]$ a nonempty open interval of $(E,\tau _{1})$. If $%
]a,b]=\{b\} $, then $b\in M$ and $]a,b]\cap M_{2}\neq \emptyset $. If $%
]a,b[\,\neq \emptyset $, then there exists $x_{n}\in \,]a,b[$, \textit{i.e.} 
$x_{n}\in \,]a,b]\cap M_{2}$. We conclude that $M_{2}$ is dense in $(E,\tau
_{1})$. $\blacksquare $ \vskip2pt plus 1pt minus 1pt

4). We need only observe that clearly, any subset $E^{\prime}$ of a strongly
Lindel\"of space~$E$ is strongly Lindel\"of for the induced topology, and
then apply~3) to~$E^{\prime}$. $\blacksquare $

\begin{remark}
\label{gg} Suppose that $(E, \tau_0)$ is a strongly Lindel\"{o}f space. Then 
$E\times \{0, 1\} $ has the same property.
\end{remark}

\noindent \textit{Proof}. By Lemma~\ref{cw}, the space $(E,\tau _{j})$ is a
strongly Lindel\"{o}f space, when $j\in \{1,2\}$. Using Lemma~\ref{hi} and
Remark~\ref{hs}, one obtains that $(E\times \{0\},\tau _{0})$ and $(E\times
\{1\},\tau _{0})$ are strongly Lindel\"{o}f spaces, therefore $E\times
\{0,1\}=E\times \{0\}\cup E\times \{1\}$ has the same property. $%
\blacksquare $

\begin{remark}
\label{wn} Suppose that $(E, \tau_0)$ is a separable strongly Lindel\"{o}f
space (resp. every subset of $(E, \tau_0)$ is separable). Then $E \times
\{0, 1\}$ is separable (resp. every subset of $E \times \{0, 1\}$ is
separable).
\end{remark}

\noindent \textit{Proof}. By Lemma~\ref{cw}, $(E,\tau _{j})$ is separable, $%
j\in \{1,2\}$ (resp. every subset of $(E,\tau _{j})$ is separable, $j\in
\{1,2\})$. On the other hand, by Lemma~\ref{hi} and Remark~\ref{hs}, $%
(E,\tau _{1})$ is homeomorphic to $E\times \{0\}$ and $(E,\tau _{2})$ is
homeomorphic to $E\times \{1\}$, it follows that $E\times \{0\},E\times
\{1\} $ are separables (resp. every subset of $E\times \{0\}$ or of $E\times
\{1\}$ is separable), hence $E\times \{0,1\}=E\times \{0\}\cup E\times \{1\}$
is separable (resp. every subset of $E\times \{0,1\}=E\times \{0\}\cup
E\times \{1\}$ is separable). $\blacksquare $

\begin{proposition}
\label{gf} i) Assume that every subset of $(E,\tau _{0})$ is separable. Then 
$B_{1}^{{}}$ is a strongly Lindel\"{o}f space.

ii) Assume that $(E,\tau _{0})$ is separable and that $(E,\leq )$ satisfies
the condition~$(\ast )$. Then $B_{1}^{{}}$ has a countable basis.
\end{proposition}

\noindent\textit{Proof}.

\textit{i)}. Let $B$ be a subset of $E$ and let $(u_{n})_{n\geq 0}$ be a
dense sequence in $(B,\tau _{0})$. For every $u\in B$ denote $V_{u}=\{\theta
\in B_{1}^{{}};\psi (\theta ,u)=1\}$. Let us show that $\underset{u\in B}{%
\cup }V_{u}=\underset{n\in \mathbb{N}}{\cup }V_{u_{n}}$.

Let $\theta \in \underset{u\in B}{\cup }V_{u}$. There exists $u\in B$ such
that $\theta \in V_{u}$, hence $u\in H_{\theta }=\{t\in B;\psi (\theta
,t)=1\}$. Thus there exists $n\in \mathbb{N}$ such that $u_{n}\in H_{\theta
} $, this implies that $\theta \in V_{u_{n}}$. It follows that $\underset{%
u\in B}{\cup }V_{u}=\underset{n\in \mathbb{N}}{\cup }V_{u_{n}}$.

If $V_{u}^{\prime }=\{\theta \in B_{1}^{{}};\psi (\theta ,u)=0\}$, $u\in B$,
by a similar argument, we show that $\underset{u\in B}{\cup }V_{u}^{\prime }=%
\underset{n\in \mathbb{N}}{\cup }V_{u_{n}}^{\prime }$. For every $u,v\in B$
let 
\begin{equation*}
W_{u,v}=\{\theta \in B_{1}^{{}};\psi (\theta ,u)=1\text{ and }\psi (\theta
,v)=0\}.
\end{equation*}%
Show that $\underset{u,v\in B}{\cup }W_{u,v}=\underset{m,n\in \mathbb{N}}{%
\cup }W_{u_{m},u_{n}}$.

Let $\theta \in W_{u,v}$. Note that $u\in H_{\theta }=\{t\in B;\psi (\theta
,t)=1\}$ and $v\in H_{\theta }^{\prime }=\{t\in B;\psi (\theta ,t)=0\}$,
hence there $u_{m}\in H_{\theta }$ and $u_{n}\in H_{\theta }^{^{\prime }}$.
It is obvious that $\theta \in W_{u_{m},u_{n}}$, this implies that $\underset%
{u,v\in B}{\cup }W_{u,v}=\underset{m,n\in \mathbb{N}}{\cup }W_{u_{m},u_{n}}$%
. $\blacksquare $ \medskip

\textit{ii)}. Let $(a_{n})_{n\geq 0}$ be a dense sequence in $(E,\tau _{0})$%
. For every $m,n\in \mathbb{N}$ consider the sets $W_{m,n}=\{\theta \in
B_{1}^{{}};\psi (\theta ,a_{n})=0\text{ and }\psi (\theta ,a_{m})=1\}$, $%
V_{m}=\{\theta \in B_{1}^{{}};\psi (\theta ,a_{m})=1\}$, $V_{n}^{\prime
}=\{\theta \in B_{1}^{{}};\psi (\theta ,a_{n})=0\}$, and $D=\{(m,n)\in 
\mathbb{N}^{2};W_{m,n}\neq \emptyset \}$. Let us show that 
\begin{equation*}
\{W_{m,n};(m,n)\in D\}\cup \{V_{m};m\in \mathbb{N}\}\cup \{V_{n}^{\prime
};n\in \mathbb{N}\}
\end{equation*}%
forms a basis of $B_{1}^{{}}$.

In the proof of Proposition~\ref{omq}, we proved that $u$ is in the closure
of $\{x\in E;x<u\}$. Let now $t,u\in E$ and 
\begin{equation*}
W_{t,u}=\{\theta \in B_{1}^{{}};\psi (\theta ,u)=0\text{ and \ }\psi (\theta
,t)=1\}\neq \emptyset .
\end{equation*}%
Pick $\varphi \in W_{t,u}$. The set $\{x\in E;\psi (\varphi ,x)=0\}$ is an
open neighborhood of $u$, hence $\{x\in E;\psi (\varphi ,x)=0\}\cap \{x\in
E;x<u\}$ is a nonempty open subset of $E$, hence there exists $n\in \mathbb{N%
}$ such that $\psi (\varphi ,a_{n})=0$ and $a_{n}<u$. By a similar argument,
one shows that there exists $m\in \mathbb{N}$ such that $\psi (\varphi
,a_{m})=1$ and $a_{m}>t$. This yields that $\varphi \in W_{m,n}\subset
W_{t,u}$.

By a reasoning similar to the previous one, we show that for every $\varphi
\in \{\theta \in B_{1}^{{}};\psi (\theta ,u)=0\}$ (resp. $\varphi \in
\{\theta \in B_{1}^{{}};\psi (\theta ,t)=1\})$, there exists $n\in \mathbb{N}
$ such that $\varphi \in V_{n}^{\prime }\subset \{\theta \in B_{1}^{{}};\psi
(\theta ,u)=0\}$ (resp. there exists $m\in \mathbb{N}$ such that $\varphi
\in V_{m}\subset \{\theta \in B_{1}^{{}};\psi (\theta ,t)=1\})$. $%
\blacksquare $

\begin{proposition}
\label{ee} Assume that $(E, \tau_0)$ is a strongly Lindel\"{o}f space. Then

I) ($B_{j}^{{}},\tau _{0})$ is a strongly Lindel\"{o}f space ($(B_{j},\tau
_{0})$ viewed as topological subspace of $(K,\tau _{0}))$, $j\in \left\{
1,2\right\} .$

II) Assume that every subset of $(E,\tau _{0})$ is separable. Then every
subset of ($B_{j}^{{}},\tau _{0})$ is separable, $j\in \left\{ 1,2\right\} .$
\end{proposition}

\noindent \textit{Proof}. We shall prove I) and II) for $j=3,$ the proof for 
$j=2$ is similar.

I).

\emph{Step 1:} Let us show that the map $h_{3}^{{}}:(B_{3}^{{}},\tau
_{0})\rightarrow (E,\tau _{1})$ is continuous.

Since $h_{3}^{{}}:(B_{3}^{{}},\tau _{0})\rightarrow (E,\tau _{0})$ is
continuous (by Proposition~\ref{cc}), it suffices to show that $%
(h_{3}^{{}})^{-1}(V)$ is an open subset of $B_{3}^{{}}$, when $V=\{t\in
E;t\leq \alpha \}$, $\alpha \in E$. Let $\alpha \in E$ and let $V=\{t\in
E;t\leq \alpha \}$. Consider $\theta _{0}\in (h_{3}^{{}})^{-1}(V)$ and $%
V_{1}=\{\theta \in B_{3}^{{}};\psi (\theta ,h_{3}^{{}}(\theta _{0}))=0\}$.
It is clear that $V_{1}$ is an open neighborhood of $\theta _{0}$. On the
other hand, if $\theta \in V_{1}$, $h(\theta )\leq h(\theta _{0})\leq \alpha 
$, hence $\theta \in (h_{3})^{-1}(V)$. This implies that $V_{1}\subset
(h_{3})^{-1}(V)$. Thus $(h_{3}^{{}})^{-1}(V)$ is an open subset of $%
(B_{3}^{{}},\tau _{0})$. $\blacksquare $\vskip2pt plus 1pt minus 4pt

\emph{Step 2}: Let us show that the map $h_{3}^{{}}:(B_{3}^{{}},\tau
_{0})\rightarrow (h_{3}^{{}}(B_{3}^{{}}),\tau _{1})$ is a homeo\-mor\-phism.

Consider $t\in E$ and $V_{t}=\{\theta \in B_{3}^{{}};\psi (\theta ,t)=1\}$.
By defi\-ni\-tion of $B_{3}^{{}}$, $V_{t}=\{\theta \in
B_{3}^{{}};h_{3}^{{}}(\theta )>t\}$, it follows that $h_{3}^{{}}(V_{t})=\{u%
\in h_{3}^{{}}(B_{3}^{{}});u>t\}$ is an open subset of ($%
h_{3}^{{}}(B_{3}^{{}}),\tau _{1})$. In the same way, if we let $%
W_{t}=\{\theta \in B_{3}^{{}};\psi (\theta ,t)=0\}=\{\theta \in
B_{3}^{{}};h_{3}^{{}}(\theta )\leq t\}$, we see that $h_{3}^{{}}(W_{t})=\{u%
\in h_{3}^{{}}(B_{3}^{{}});u\leq t\}$ is an open set of $%
(h_{3}^{{}}(B_{3}^{{}}),\tau _{1})$. Thus $h_{3}^{{}}$ is a homeomorphism.

By Lemma~\ref{cw} and Remark~\ref{nl} $(h_{3}^{{}}(B_{3}^{{}}),\tau _{1})$
is a strongly Lindel\"{o}f space, hence $(B_{3}^{{}},\tau _{0})$ has the
same property. $\blacksquare $

II). By Lemma~\ref{cw}, every subset of $(h_{3}^{{}}(B_{3}^{{}}),\tau _{1})$
is separable, using step~2, one obtains that every subset of $%
(B_{3}^{{}},\tau _{0})$ is separable. $\blacksquare $

\begin{corollary}
\label{vg} Assume that $(E,\tau _{0})$ is a strongly Lindel\"{o}f space, $%
(E,\leq )$ is complete and every subset of $(E,\tau _{0})$ is separable.
Then $K$ is separable.
\end{corollary}

\noindent \textit{Proof}. First et us show that $(B_{1}^{{}},\tau _{0})$ is
separable. By an argument similar to that of Proposition~\ref{ee}-II), one
shows that the map $h_{1}^{{}}:(B_{1}^{{}},\tau _{0})\rightarrow
(h(B_{1}^{{}}),\tau _{2})$, $\theta \rightarrow h(\theta )$ is a
homeo\-morphism. By Lemma~\ref{cw}, $(h(B_{1}^{{}}),\tau _{2})$ is
separable, hence $(B_{1}^{{}},\tau _{0})$ is separable. On the other hand,
by Proposition~\ref{ee}, $(B_{j}^{{}},\tau _{0})$ is separable, $j=0,1$.
Since $K=B_{1}^{{}}\cup B_{2}^{{}}\cup B$, $K$ is separable. $\blacksquare $

\begin{corollary}
\label{mm} Assume that $(E,\tau _{0})$ is a separable space and satisfies
the condition~$(\ast )$. Then $K$ is a separable strongly Lindel\"{o}f space.
\end{corollary}

\noindent \textit{Proof}. By Proposition~\ref{gf}, $B_{1}^{{}}$ has a
countable basis. On the other hand by Remark~\ref{gd}, $(E,\tau _{0})$ has a
countable basis, this implies that $(E,\tau _{0})$ is a strongly Lindel\"{o}%
f space and every subset of $(E,\tau _{0})$ is separable, therefore by
Proposition~\ref{ee}, $B_{2}^{{}},B_{3}^{{}}$ are separables strongly Lindel%
\"{o}f spaces. Thus $K=B_{1}^{{}}\cup B_{2}^{{}}\cup B_{3}^{{}}$ is a
separable strongly Lindel\"{o}f space. $\blacksquare $

\begin{remark}
\label{NO}Let ($L,<)$ non-empty totally ordered set; assume that the order
topology of $L$ is compact and seprable. By \cite{Osta} every subset of $L$
is separable and the order topology of $L$ is strongly Lindel\"{o}f.
\end{remark}

\begin{corollary}
\label{zo} Assume that $(E,\tau _{0})$ is separable and connected. Then $K$
is sepa\-rable space.
\end{corollary}

Let $P$ be a Polish space. We denote by $B_1(P)$ the space of first class
functions (with real values or complex values) on $P$.

We denote by $\tau_p$ the topology of pointwise convergence on $B_1(P)$.

\begin{remark}
\label{az} Let $(E,\tau _{0})$ be a separable strongly Lindel\"{o}f space.
Suppose that the map $h_{j}^{{}}:B_{j}^{{}}\rightarrow E\setminus \{-\infty
,+\infty \}$ is onto, for $j\in \{2,3\}$. Then $K\setminus B_{1}^{{}}$ is
separable.
\end{remark}

\noindent \textit{Proof}. Note first that $(E\setminus \{-\infty ,+\infty
\},\tau _{0})$ is a separable strongly Lindel\"{o}f space. By Lemma~\ref{cw}%
, $(E\setminus \{-\infty ,+\infty \},\tau _{1})$ is a separable space. On
the other hand by step~2 of Proposition~\ref{ee}, the map $%
h_{3}^{{}}:(B_{3}^{{}},\tau _{0})\rightarrow (E\setminus \{-\infty ,+\infty
\},\tau _{1})$ is a homeomorphism, therefore $(B_{3}^{{}},\tau _{0})$ is
separable. By a similar argument, one shows that $(B_{2}^{{}},\tau _{0})$ is
separable.

As $K-B_{1}^{{}}=B_{2}^{{}}\cup B_{3}^{{}}$, $K-B_{1}^{{}}$ is separable. $%
\blacksquare $

\begin{definition}
\label{ad} Let $L$ be a compact Hausdorff space. We say that $L$ is a
Rosenthal compact set, if there exists a Polish space $P$ such that $L$
continuously embedded in $(B_{1}(P),\tau _{p})$ \cite{Gode}.
\end{definition}

\begin{definition}
\label{nn} Let $X$ be a Hausdorff space. $X$ is said to be an analytic set,
if $X$ is a continuous image of a Polish space \cite[chap.~II, p.~96]{Sch(1)}%
.
\end{definition}

\begin{proposition}
\label{ls} Suppose that $(E,\tau _{0})$ is an analytic set and suppose that $%
-\infty ,+\infty \notin E$. Then $K$ is a Rosenthal compact set.
\end{proposition}

\noindent \textit{Proof}. Since $(E,\tau _{0})$ is analytic, there exists a
Polish space $P$ and a surjective continuous map $\phi :P\rightarrow (E,\tau
_{0})$. Let $\theta \in K$. Consider the map $\sigma _{\theta }:P\rightarrow
\{0,1\}$, $x\in P\rightarrow \psi (\theta ,\phi (x))$. Note at first that if 
$\theta \in B_{1}^{{}}$, $\sigma _{\theta }$ is a continuous function. Let $%
\theta \in B_{2}^{{}}$. Prove that $\sigma _{\theta }$ is a first class
function.

It is clear that 
\begin{equation*}
(\sigma _{\theta })^{-1}\{1\}=\{x\in P;\psi (\theta ,\phi (x))=1\}=\{x\in
P;\phi (x)\leq h(\theta )\}
\end{equation*}%
which is a closed subset of $P$. On the other hand 
\begin{equation*}
(\sigma _{\theta })^{-1}\{0\}=\{x\in P;\psi (\theta ,\phi (x))=0\}=\{x\in
P;\phi (x)>h(\theta )\}
\end{equation*}%
is an open subset of $P$. Therefore the sets $(\sigma _{\theta })^{-1}\{1\}$%
, $(\sigma _{\theta })^{-1}\{0\}$ are $G_{\delta }$ subsets of $P$. It
follows that $\sigma _{\theta }$ is a first class function. If $\theta \in
B_{3}^{{}}$, by a similar argument, one shows that $\sigma _{\theta }$ is a
first class function. One defines the map $\sigma :K\rightarrow B_{1}(P)$,
by $\sigma (\theta )=\sigma _{\theta }$, $\theta \in K$. It is obvious that $%
\sigma $ is into and continuous, hence $K$ is a Rosenthal compact set. $%
\blacksquare $

\begin{proposition}
\label{wb} Suppose that $(E,\leq )$ is complete and $(E,\tau _{0})$ is a
strongly Lindel\"{o}f space. Then $K$ is a strongly Lindel\"{o}f space.
\end{proposition}

\noindent \textit{Proof}. We can suppose that $-\infty ,+\infty \notin E$.
Denote $K_{2}=B_{1}^{{}}\cup B_{2}^{{}}-\{\theta ^{\prime },\theta ^{\prime
\prime }\}$ and consider the map $h_{\mid _{K_{2}}}:(K_{2},\tau
_{0})\rightarrow (h(K_{2}),\tau _{2})$ (note that $\sup \{t;\psi (\theta
,t)=1\}$ exists for all $\theta \in K_{2}$ by Lemma~\ref{kt} and Remark~\ref%
{lmv}). By hypothesis, $(E,\tau _{0})$ is a strongly Lindel\"{o}f space,
hence by Lemma~\ref{cw}, $(E,\tau _{2})$ has the same property, therefore $%
(h(K_{2}),\tau _{2})$ is a strongly Lindel\"{o}f space. Thus $(K_{2},\tau
_{0})$ is a strongly Lindel\"{o}f space, because $h_{\mid _{K_{2}}}$ is a
homeomorphism (By an argument similar to that of Proposition~\ref{ee}-I)).
Observe finally that $K=K_{2}\cup B_{3}^{{}}\cup \{\theta ^{\prime },\theta
^{\prime \prime }\}$, $K_{2}$ is a strongly Lindel\"{o}f space and $%
B_{3}^{{}}$ is a strongly Lindel\"{o}f space by Proposition~\ref{ee}, hence $%
K$ is a strongly Lindel\"{o}f space. $\blacksquare $

Let $L$ be a Hausdorff compact space. Denote $M_{{}}^+(L)$ the Radon
probability measures space on $L$. \bigskip

For every Banach space $X$, we denote by $B_{X}$ the closed unit ball of~$X$.

\begin{definition}
\label{ih} \cite[chap.~IX.8]{Bourb} A topological space $(X, \tau)$ is said
to be a completely regular space, if $(X, \tau)$ can be embedded in a
compact Hausdorff space.
\end{definition}

It is clear that every topological subspace of a completely regular space is
a completely regular space.

Recall that $\mu ^{\prime }$ is the counting measure defined on $E$ and $%
h^{\prime }(\theta )=\sup \{t\in E;\psi ^{\prime }(\theta ,t)=1\}$, $\theta
\in K^{\prime }\setminus B_{1}^{\prime }$, where $\psi ^{\prime }(\theta
,t)=\theta (f^{t})$, $t\in E$. Denote by $B_{1}^{\prime }$ the set of $%
\theta \in K^{\prime }$ such that the map $t\in E\rightarrow \psi ^{\prime
}\theta ,t)\in \{0,1\}$ is continuous, $B_{2}^{\prime }=\{\theta \in
K^{\prime }\setminus B_{1}^{\prime };\psi ^{\prime }(\theta ,h^{\prime
}(\theta ))=1\}$ and by $B_{3}^{\prime }=\{\theta \in K^{\prime }\setminus
B_{1}^{\prime };\psi ^{\prime }(\theta ,h_{{}}^{\prime }(\theta ))=0\}$.

\begin{remark}
\label{cs} Suppose that $(E, \leq)$ is complete ($-\infty, +\infty \notin E$%
). Then $(E, \tau_j)$ is a completely regular, $j \in \{1, 2\}$.
\end{remark}

\noindent \textit{Proof}. We shall show this remark for $j=2$. Consider $\mu
=\mu ^{\prime }$. By Corollary~\ref{ds}, $h^{\prime }(K_{2})=E$. On the
other hand, in the proof of Proposition~\ref{wb}, we have shown that $%
h_{{}}^{\prime }:B_{2}^{\prime }\cup B_{1}^{\prime }\setminus \{\theta
^{\prime },\theta ^{\prime \prime }\}=K_{2}\rightarrow (E,\tau _{2})$ is \ a
homeomorphism. Thus $(E,\tau _{2})$ is a completely regular space. $%
\blacksquare $

Let 
\begin{equation*}
E_{1}=\{a\in E;\{x\in E;x>a\}\text{ is not closed subset of }E\}
\end{equation*}%
and 
\begin{equation*}
E_{2}=\{a\in E;\{x\in E;x<a\}\text{ is not closed subset of }E\}.
\end{equation*}

\begin{remark}
\label{wj} The space $(E_{2},\tau _{1})$ is completely regular space.
\end{remark}

\noindent \textit{Proof}. Consider $\mu =\mu ^{\prime }$. By Lemma~\ref{mq}%
-b), the restriction of $h_{{}}^{\prime }$ (in Remark~\ref{cs}) to $%
B_{3}^{\prime }$ is onto on $E_{2}$ and by step~2 of Proposition~\ref{ee},
this map is homeomorphism. On the other hand, $B_{3}^{\prime }$ is a
completely regular space, hence $(E_{2},\tau _{1})$ has the same property. $%
\blacksquare $

\begin{remark}
\label{db} Suppose that $(E,\leq )$ satisfies the condition~$(\ast )$. Since
the point $a$ is in the closure of the set $\{x\in E;x<a\}$, for every $a\in
E,$ the set $\{x\in E;x<a\}$ is not closed. We conclude that $E=E_{2}$.
\end{remark}

\begin{proposition}
\label{ai} Suppose that $E$ is a Polish space and $\mu (\{a\})=0$ for every $%
a\in E$ (i.~e., $\mu $ is diffuse on the Borel $\sigma $-algebra of $E$).
Then for every $\xi \in C(K)^{\ast }$ the map $t\in (E,\tau _{0})\rightarrow
\xi (f^{t})$ is a first-class function.
\end{proposition}

\noindent\textit{Proof}. It is suffices to show Proposition~\ref{ai} when $%
\xi \in B_{C(K)^*}$.

\emph{Step 1}: Let us show that for every $g\in L^1(E, \mu)$ the map $t\in
E\rightarrow (g,f^t)=\mathop{\displaystyle \int} \limits_{E}g(x)f^t(x)d\mu
(x)$ is continuous.

Let $g\in L^1(E, \mu)$ and let $(t_n)_{n\geq 0}$ be a sequence in $E$ such
that $t_n\rightarrow t$ in $(E, \tau_0)$. Consider $y\in E$ such that $y\neq
t$. We shall show that $f^{t_n}(y) \rightarrow_{n \rightarrow +\infty}
f^t(y) $.

We show that there does not exist two subsequences $(t_{m_k})_{k\geq 0}$ and 
$(t_{n_k})_{k\geq 0}$ such that $f^{t_{m_k}}(y)=1$ and $f^{t_{n_k}}(y)=0 $
for every $k\in \mathbb{N}$. Assume there are such subsequences. We have
then $t_{n_k}$ $\geq y>t_{m_k}$ for every $k \in \mathbb{N}$, by passing to
the limit we obtain $y = t$, which is impossible. \smallskip

Thus we distinguish two cases:

\emph{Case 1}, there exists $n_0\in \mathbb{N}$ such that $f^{t_n}(y)=1$ for
every $n\geq n_0$: Note that $y>t_n$, $\forall n\geq n_0$, hence $y \geq t$.
Since $y\neq t$, $y>t$, i.e. $f^{t_n}(y)\rightarrow _{n\rightarrow +\infty
}f^t(y)$.

\emph{Case 2}, there exists $n_0$ such that $f^{t_n}(y)=0$ for every $n \geq
n_0$: This implies that $y\leq t_n$ for every $n\geq n_0$, by passing to the
limit, we obtain $y\leq t$. Thus $f^{t_n}(y)\underset{ n\rightarrow \infty }{%
\rightarrow }f^t(y)=0$.

By hypothesis on $\mu $, $(g,f^{t_{n}})=\mathop{\displaystyle \int}%
\limits_{E\setminus \{t\}}g(y)f^{t_{n}}(y)\,d\mu (y)$. By Lebesgue's
dominated convergence theorem, we have 
\begin{eqnarray*}
(g,f^{t_{n}}) &=&\int_{E\setminus \{t\}}g(y)f^{t_{n}}(y)\,d\mu (y) \\
&&\underset{n\rightarrow \infty }{\rightarrow }\int_{E\setminus
\{t\}}g(y)f^{t}(y)\,d\mu (y) \\
&=&(g,f^{t}).\text{ }\blacksquare
\end{eqnarray*}

\emph{Step 2}: Let $\xi \in M^{+}(K)$. Let us show that the map $t\in
E\rightarrow \xi (f^{t})\in \mathbb{R}$ is a first class function.

By Hahn-Banach theorem there exists $\xi ^{\prime }\in B_{L^{\infty }(\mu
)^{\ast }}$ which extends $\xi $ with $\Vert \xi \Vert =\Vert \xi ^{\prime
}\Vert $. Since $B_{L^{1}(\mu )}$ is $w^{\ast }$-dense in $B_{(L^{\infty
}(\mu ))^{\ast }}=B_{(L^{1}(\mu ))^{\ast \ast }}$, there is a net $%
(g_{\alpha })_{\alpha \in I}$ in the unit ball of $L^{1}(\mu )$ such that $%
g_{\alpha }\underset{\mathcal{U}}{\rightarrow }\xi ^{\prime }$ for the
topology $w^{\ast }$ in $[L^{\infty }(\mu )]^{\ast }$. There exist $\xi
_{1},\xi _{2}\in M^{+}(K)$ such that $g_{\alpha }^{+}\rightarrow \xi _{1}\in
M^{+}(K)$ and $g_{\alpha }^{-}\rightarrow \xi _{2}\in M^{+}(K)$, $\sigma
\lbrack C(K)^{\ast },C(K)]$, hence $\xi =\xi _{1}-\xi _{2}$, because $%
g_{\alpha }=g_{\alpha }^{+}-g_{\alpha }^{-}$, for all $\alpha $ and $%
M^{+}(K) $ is $\sigma (C(K)^{\ast },C(K))$ compact. By using Proposition~\ref%
{ls}, $K$ is a Rosenthal compact set. On the other hand, by \cite{Gode}, $%
M_{{}}^{+}(K) $ is a Rosenthal compact set, and by \cite{Bour-Frem-Tal} $%
M_{{}}^{+}(K)$ is angelic. Therefore, there exist two subsequences $%
(g_{j_{k}}^{+})_{k\geq 0}$, $(g_{i_{k}}^{-})_{k\geq 0}$, such that $%
g_{j_{k}}^{+}\underset{k\rightarrow +\infty }{\rightarrow }\xi _{1}$ and $%
g_{i_{k}}^{-}\underset{k\rightarrow +\infty }{\rightarrow }\xi _{2}$, for
the topology $w^{\ast }$, i.e. $h_{k}=g_{j_{k}}^{+}-g_{i_{k}}^{-}\underset{%
k\rightarrow +\infty }{\rightarrow }\xi _{1}-\xi _{2}=\xi $ for the topology 
$w^{\ast }$. Thus, $(h_{k},f^{t})\underset{k\rightarrow +\infty }{%
\rightarrow }\xi (f^{t})$ for every $t\in E$.

Finally, by step~1, the map $t\in E\rightarrow (h_{k},f^{t})$ is continuous
for every~$k$, which implies that $t\in E\rightarrow \xi (f^{t})$ is a first
class function. $\blacksquare $

If $\xi $ in the closed unit ball of $[C(K)]^{\ast }$, it suffices to see
that there exists $u^{+},u^{-},v^{+},v^{-}\in M^{+}(K)$ such that $\xi
=(u^{+}-u^{-})+i(v^{+}-v^{-})$, by using step~2, we obtain that the map $%
t\in E\rightarrow \xi (f^{t})$ is a first class function. $\blacksquare $

Denote $J=[E\cup \{+\infty \}]\times \{1\}\cup \lbrack E\cup \{-\infty
\}]\times \{0\}$ ($-\infty ,+\infty \notin E)$. As $J$ is a subset of $%
\overline{E}\times \{0,1\}$ and $\overline{E}\times \{0,1\}$ is totally
ordered (for the lexicographic order), $J$ is totally ordered. We denote
again by $\tau _{0}$ the order topology of $J$. \bigskip

Let $J_{1}$ be a subset of $J$. Recall that the notation $(J_{1},\tau _{0})$
means that the topology $\tau _{0}$ on $J_{1}$ is the topology formed of
relative open subsets of $J_{1}$ in $(J,\tau _{0})$. Let 
\begin{equation*}
J_{1}=[h(B_{2}^{{}}\cup B_{1}^{{}}\setminus \{\theta ^{\prime \prime
}\})\times \{1\}]\cup \lbrack h(B_{3}^{{}}\cup \{\theta ^{\prime \prime
}\})\times \{0\}].
\end{equation*}

In following Theorem $J_{1}$ is equipped by the relative topology $\tau
_{0}. $

\begin{theorem}
\label{vb} Assume that $(E,\leq )$ is complete and $-\infty $, $+\infty
\notin E$. Then $K$ is topological homeomorphic and order isomorphic to $%
J_{1}.$
\end{theorem}

\noindent \textit{Proof}.

One defines $T:K\rightarrow (J_{1},\tau _{0})$, by 
\begin{equation*}
T(\theta )=\left\{ 
\begin{array}{c}
(h(\theta ),1),\text{ if }\theta \in B_{2}^{{}}\cup (B_{1}^{{}}\setminus
\{\theta ^{\prime \prime }\}) \\ 
(h(\theta ),0),\text{ if }\theta \in B_{3}^{{}}\cup \{\theta ^{\prime \prime
}\}%
\end{array}%
\right. ,\quad \theta \in K.
\end{equation*}

Let us show that $T$ is a homeomorphis. Since the restriction of $h$ to $%
B_{j}^{{}}$ is injective, $j\in \{1,2,3\}$, $T$ is injective. So suffices to
show that $T$ is continuous, because $K$ is a Hausdorff compact space.

Let $(\alpha ,\beta )\in J$ and let $V=\{(x,y)\in J_{1};(x,y)<(\alpha ,\beta
)\}$ be an open subset of~$J_{1}$. We prove that $T^{-1}(V)$ is an open
subset of $K$. We distinguish two cases:

\emph{Case 1}, $\beta =0$: It is obvious that $T^{-1}(V)=(h)^{-1}(\left\{
t\in \overline{E};t<\alpha \right\} )$ which is an open subset of $K$, by
Proposition~\ref{xh}. \vskip2pt plus 1pt minus 4pt

\emph{Case 2}, $\beta =1$: Let $\theta \in T^{-1}(V)$. We shall find an open
neighborhood $C$ of $\theta $ in $T^{-1}(V)$. We distinguish again two
cases: \vskip 2pt plus 1pt minus 4pt

\emph{Case 2-a}, $h(\theta )<\alpha $: Put $C=h^{-1}(\left\{ t\in \overline{E%
};t<\alpha \right\} )$ (note that $\theta \in C)$. if $\varphi \in C$, $%
T(\varphi )=(h(\varphi ),j)<(\alpha ,1)$, ( $j\in \left\{ 0,1\right\} $),
hence $C\subset T^{-1}(V)$. On the other hand, by Proposition~\ref{xh}, $C$
is an open subset of $K$. \vskip2pt plus 1pt minus 4pt

\emph{Case 2-b}, $h(\theta )=\alpha >-\infty $: Note that $\theta \neq
\theta ^{\prime }$, because $\theta \in T^{-1}(V)$. Denote $C=\left\{
\varphi \in K;\psi (\varphi ,h(\theta ))=0\right\} $. $C$ is an open subset
of $K$. \vskip2pt plus 1pt minus 4pt

Let us show that $C$ contains $\theta$. \vskip 2pt plus 1pt minus 4pt

Since $\theta \in T^{-1}(V)$, $T(\theta )=(h(\theta ),j)<(h(\theta ),1)$,
hence $j=0$. By definition of $T$ $\psi (\theta ,h(\theta ))=0$, i.e. $%
\theta \in C$.

It remains to show that $C\subset T^{-1}(V)$. Let in fact $\varphi \in C$;
this means that $\psi (\varphi ,h(\theta ))=0$, we conclude that $h(\varphi
)\leq h(\theta )$, if $h(\varphi )=h(\theta )=\alpha $, then $\psi (\varphi
,h(\varphi ))=\psi (\varphi ,h(\theta ))=0$. It follows that $T(\varphi
)=(h(\varphi ),0)<(\alpha ,1)$. Thus $\varphi \in T^{-1}(V)$.

If $h(\varphi )<h(\theta )=\alpha $, then $T(\varphi )=(h(\varphi
),j)<(\alpha ,1)$ ($j\in \left\{ 0,1\right\} )$, hence $\varphi \in
T^{-1}(V) $. Therefore $C\subset T^{-1}(V)$.

Let now $(\alpha ,\beta )\in J$ and let $V=\left\{ (x,y)\in
J_{1};(x,y)>(\alpha ,\beta )\right\} $ be an open subset of $J_{1}$. Show
that $T^{-1}(V)$ is an open subset of $K$.

\emph{Case 1}, $\beta =1$: By Proposition~\ref{xh}, $T^{-1}(V)=h^{-1}(\left%
\{ t\in \overline{E};\text{ }t>\alpha \right\} )$ is an open subset of $K$.

\emph{Case 2}, $\beta =0$: Consider $\theta \in T^{-1}(V)$. We distinguish
two cases:

\emph{Case 2-c, }$h(\theta )>\alpha $: The set $C=h^{-1}(\left\{ t\in 
\overline{E};\text{ }t>\alpha \right\} )$ is an open subset of $K$, by
Proposition~\ref{xh} which contains $\theta $. The condition $h(\varphi
)>\alpha $ implies that $T(\varphi )=(h(\varphi ),j)>(\alpha ,0)$, hence $%
C\subset T^{-1}(V)$.

\emph{Case 2-d}, $h(\theta )=\alpha <+\infty $: Note that $\theta \neq
\theta ^{\prime \prime }$, because $\theta \in T^{-1}(V))$. Choose $%
C=\left\{ \varphi \in K;\psi (\varphi ,h(\theta ))=1\right\} $. The
hypothesis $\theta \in T^{-1}(V)$ implies that $T(\theta )=(h(\theta
),j)>(h(\theta ),0)$, hence $j=1$, this means that $\psi (\theta ,h(\theta
))=1$. It follows that $C$ contains $\theta $. Show that $C\subset T^{-1}(V)$%
. Let $\varphi \in C$. It is clear that $h(\varphi )\geq h(\theta )$, if $%
h(\theta )=h(\varphi )$, then $\psi (\varphi ,h(\varphi ))=\psi (\varphi
,h(\theta ))=1$, hence $T(\varphi )=((h(\varphi ),1)>(\alpha ,0)$, i.e. $%
\varphi \in T^{-1}(V)$. Finally if $h(\theta )<h(\varphi )$, then $T(\varphi
)=(h(\varphi ),j)>(h(\theta ),0)=(\alpha ,0)$, this implies that $\varphi
\in T^{-1}(V)$. Thus $C\subset T^{-1}(V)$.

Let us show that $T$ is an order isomorphism.

Let $\theta ,\varphi \in K$ such that $\theta <\varphi $. There exists $t\in
E$ such that $\psi (\theta ,t)=0$ and $\psi (\varphi ,t)=1$. This implies
that 
\begin{equation}
h(\varphi )\geq t\geq h(\theta ).  \label{g}
\end{equation}%
Obviously if $h(\theta )<h(\varphi )$, then $T(\theta )=(h(\theta
),j_{1})<(h(\varphi ),j_{2})=T(\varphi )$. If $h(\theta )=h(\varphi )$, by (%
\ref{g}) $h(\theta )=h(\varphi )=t$, hence $\theta \in B_{3}^{0}$ and $%
\varphi \in B_{2}^{{}}$, this implies that $T(\theta )=(h(\theta
),0)<(h(\varphi ),1)=T(\varphi )$.

Conversely, suppose that $T(\theta )<T(\varphi )$. If $h(\theta )<h(\varphi
) $, then $\theta <\varphi $, if $h(\theta )=h(\varphi )$, since $T(\theta
)=(h(\theta ),j_{1})<T(\varphi )=(h(\varphi ),j_{2})$, $j_{1}=0$ and $%
j_{2}=1 $. Thus $\theta \in B_{3}^{{}}$ and $\varphi \in B_{2}^{{}}$, hence $%
\theta <\varphi $. $\blacksquare $

\begin{corollary}
\label{wi} Suppose that $E$ is connected and $\mu (\left] t\wedge t^{\prime
},t\vee t^{\prime }\right] )>0$, for $t\neq t^{\prime }$. Then $K$ is
topological homeomorphic and order isomorphic to 
\begin{equation*}
J=(\left[ E\cup \left\{ +\infty \right\} \right] \times \left\{ 1\right\}
\cup \left[ E\cup \left\{ -\infty \right\} \right] \times \left\{ 0\right\}
,\tau _{0}).
\end{equation*}
\end{corollary}

\noindent \textit{Proof}. Since $E$ is connected, by Remark~\ref{ff}, $%
(E,\leq )$ is complete. On the other hand, by Proposition~\ref{omq}, $%
h_{j}^{{}}$ is surjective onto $E\setminus \{-\infty ,+\infty \}$, $j\in
\{2,3\}$, it follows that $T$ (defined in Theorem~\ref{fg}) is surjective
(note that $B_{1}^{{}}=\{\theta ^{\prime },\theta ^{\prime \prime }\}$). $%
\blacksquare $

\begin{corollary}
\label{cg} Let $E$ be an abelian group (-$\infty ,+\infty \notin E$). Assume
that $E$ separable and connected. Then $K$ is topological homeomorphic and
order isomorphic to 
\begin{equation*}
(\left[ E\cup \left\{ +\infty \right\} \right] \times \left\{ 1\right\} \cup %
\left[ E\cup \left\{ -\infty \right\} \right] \times \left\{ 0\right\} ,\tau
_{0}).
\end{equation*}
\end{corollary}

\noindent \textit{Proof}. The hypothesis that $E$ is connected implies that (%
$E,\leq )$ is complete. On the other hand, $E$ satisfies the condition~$%
(\ast )$, the example (\ref{yu}) shows that the measure $\mu =m$ verifies
the hypotheis of Corollary~\ref{wi}, hence $K$ is isomorphic to $(\left[
E\cup \{+\infty \}\right] \times \{1\}\cup \left[ E\cup \{-\infty \}\right]
\times \{0\},\tau _{0})$. $\blacksquare $

\begin{lemma}
\label{jf} Assume that $(E, \leq)$ is complete ($E$ is a separable strongly
Lindel\"{o}f space). Then:

1) $J = (\left[ (E \cup \{+\infty\}) \times \{1\} \right] \cup \left[ (E\cup
\{-\infty\}) \times \{0\} \right], \tau_0)$ is compact.

2) $(J, \tau_0)$ is a separable space.
\end{lemma}

\noindent\textit{Proof}.

1). Consider $\mu =\mu ^{\prime }$ the counting measure. By using Theorem~%
\ref{vb}, we obtain that the set $\left[ (E\cup \{+\infty \})\times \{1\}%
\right] \cup \left[ (h^{\prime }(B_{3}^{\prime })\cup \{-\infty \})\times
\{0\}\right] $ is compact in $(J,\tau _{0})$, because this space is
homeomorphic to $K^{\prime }$. By a similar argument (by changing the
statement of Theorem~\ref{vb}) we may show that the set $\left[ (h^{\prime
}(B_{2}^{\prime })\cup \{+\infty \})\times \{1\}\right] \cup \left[ (E\cup
\{-\infty \})\times \{0\}\right] $ is compact in $(J,\tau _{0})$, thus 
\begin{equation*}
\begin{array}{c}
\left[ (E\cup \{+\infty \})\times \{1\}\right] \cup \left[ (h^{\prime
}(B_{3}^{\prime })\cup \{-\infty \})\times \{0\}\right] \\ 
\cup \left[ (h^{\prime }(B_{2}^{\prime })\cup \{+\infty \})\times \{1\}%
\right] \cup \left[ (E\cup \{-\infty \})\times \{0\}\right] =J%
\end{array}%
\end{equation*}%
is compact. $\blacksquare $

2). Since $E$ is a separable strongly Lindel\"{o}f space, $(\overline{E}%
,\tau _{0})$ has the same property. By Lemma~\ref{cw}, $(\overline{E},\tau
_{2})$ is a separable strongly Lindel\"{o}f space. On the other hand Lemma~%
\ref{hi} shows us that $(\overline{E}\times \{1\},\tau _{0})$ is a separable
subspace of $(\overline{E}\times \{0,1\},\tau _{0})$. Note now that $\left[
E\cup \{+\infty \}\right] \times \{1\}=\overline{E}\times \{1\}\setminus
\{(-\infty ,1)\}$, hence $(\left[ E\cup \{+\infty \}\right] \times
\{1\},\tau _{0})$ is separable. By a similar argument, one shows that $(%
\left[ E\cup \left\{ -\infty \right\} \right] \times \left\{ 0\right\} ,\tau
_{0})$ is a separable of $\overline{E}\times \left\{ 0,1\right\} $. It
follows that $\left[ (E\cup \left\{ +\infty \right\} )\times \left\{
1\right\} \right] \cup \left[ (E\cup \left\{ -\infty \right\} )\times
\left\{ 0\right\} \right] =J$ is a separable space of $\overline{E}\times
\left\{ 0,1\right\} $. Let $i:(J,\widetilde{\tau }_{0})\rightarrow (J,\tau
_{0})$ be the identity map of $J$ ($(J,\widetilde{\tau }_{0})$, where is a
topological subspace of $\overline{E}\times \left\{ 0,1\right\} )$. It is
clear that $i$ is continuous, thus $(J,\tau _{0})$ is separable. $%
\blacksquare $ \bigskip

By an argument similar to that of Theorem~\ref{vb}, one shows the following
theorem

\begin{theorem}
\label{ct} Suppose that $(E,\leq )$ is complete and $-\infty $, $+\infty \in
E$. Then $K$ is topological homeomorphic and order isomorphic to 
\begin{equation*}
(\left[ (h(B_{2}^{{}}\cup B_{1}^{{}}\setminus \left\{ \theta ^{\prime \prime
}\right\} )\times \left\{ 1\right\} \right] \cup \left[ (h(B_{3}^{{}}\cup
\left\{ \theta ^{\prime \prime }\right\} )\times \left\{ 0\right\} \right]
,\tau _{0}).
\end{equation*}
\end{theorem}

By an argument similar to that of Corollary~\ref{wi}, one shows the
following corollary:

\begin{corollary}
\label{vv} Suppose that $E$ is connected, that $-\infty ,+\infty \in E$ and
that $\mu (\left] t\wedge t^{\prime },t\vee t^{\prime }\right] )>0$, for $%
t\neq t^{\prime }$. Then $K$ is topological homeomorphic and order
isomorphic to $(\left[ E\setminus \{-\infty \}\right] \times \{1\}\cup \left[
E\setminus \{+\infty \}\right] \times \{0\},\tau _{0})$.
\end{corollary}

\begin{theorem}
\begin{remark}
\label{fr} Let $(E,\tau _{0})$ be a separable strongly Lindel\"{o}f space.
Suppose that $(E,\leq )$ is complete and $-\infty ,+\infty \in E$. Then $K$
is a seprable..
\end{remark}
\end{theorem}

\section{Construction of a totally ordered compact set homeomorphic to $K_{%
\protect\mu }(E)$}

This part is devoted to find sufficient conditions for $K$ to be isomorphic
to $E.$

Recall that $E_{1}=\{t\in E;\{x\in E;x>t\}\text{ is not closed in }E\}$,
that $E_{2}=\{t\in E;\{x\in E;x<t\}\text{ is not closed in }E\}$, $h(\theta
^{\prime })=+\infty $ and $h(\theta ^{\prime \prime })=-\infty $. \bigskip

\begin{theorem}
\label{vl} Assume that $(E,\leq )$ is complete, that $-\infty ,+\infty \in E$%
, $E\setminus \left\{ +\infty ,-\infty \right\} =E_{1}\cup E_{2}$ and for
every $t\in E_{1}\cup E_{2}$, there exists $b\in E_{1}$ and $a\in E_{2}$
such that $[t]=\{a,b\}$, $b>a$. Then the map $h:K\rightarrow E$ is a
homeo\-morphism.
\end{theorem}

\noindent\textit{Proof}.

\emph{Step 1}: Show that $h:K\rightarrow E$ is onto.

Let $b\in E_{1}$. By hypothesis, there exists $a\in E_{2}$ such that $%
[b]=\{b,a\}$ and $a<b$. By Lemma~\ref{mq}, there is $\theta _{b}\in
B_{2}^{{}}$ such that $h(\theta _{b})=b$. In the same way, for every $%
a^{\prime }\in E_{2}$, there is $\theta _{a^{\prime }}^{\prime }\in
B_{3}^{{}}$ such that $h(\theta _{a^{\prime }}^{\prime })=a^{\prime }$.
\bigskip $\blacksquare $

\emph{Step 2}: Let us show that the map $h:K\rightarrow E$ is into.

Suppose there exist $\theta _{1},\theta _{2}\in K$ such that $h(\theta
_{1})=h(\theta _{2})=t_{0}$. Since the restriction of $h$ to $B_{1}^{{}}\cup
B_{2}^{{}}$ (resp. to $B_{3}^{{}})$ is injective, by step~1 of Proposition~%
\ref{wb} (resp. by Remark~\ref{cnn})), it is enough to treat the following
cases:

\emph{Case 1}, $\theta _{1}\in B_{2}^{{}}$ and $\theta _{2}\in B_{3}^{{}}$:

\emph{Case 1-a}, $t_{0}\in E_{1}$: By hypothesis, there exists $a\in E_{2}$
such that $[t_{0}]=\{a,b\}$ and $b=t_{0}>a$. On the other hand, since $a\in
\lbrack t_{0}]=[b]$ and $\theta _{2}\in B_{3}^{{}},$ $\psi (\theta
_{2},a)=\psi (\theta _{2},b)=\psi (\theta _{2},h(\theta _{2}))=0$, \ this
implies that $a\geq h(\theta _{2})=b,$ it is impossible because $b=t_{0}>a$.

\emph{Case 1-b}, $t_{0}=a\in E_{2}$: By hypothesis, there exists $b\in E_{1}$
such that $[t_{0}]=\{a,b\}$ and $t=t_{0}<b$. Observe that $\psi (\theta
_{1},b)=\psi (\theta _{1},a)=\psi (\theta _{1},h(\theta _{1}))=1$. Thus $%
b\leq h(\theta _{1})=t_{0},$ which is impossible.

\emph{Case 2}, $\theta _{1}\in B_{1}^{{}}\setminus \left\{ \theta ^{\prime
},\theta ^{\prime \prime }\right\} $ and $\theta _{2}\in B_{3}^{{}}$:

\emph{Case 2-c}, $t_{0}=b\in E_{1}$: By hypothesis, there exists $a\in E_{2}$
such that $[t_{0}]=\{a,b\}$ and $a<b=t_{0}$. Note that $\psi (\theta
_{2},a)=\psi (\theta _{2},b)=0,$ hence $a\geq h(\theta _{2})=t_{0}=b,$ this
is impossible.

\emph{Case 2-d,} $t_{0}=a\in E_{2}$: There exists $b\in E_{1}$ such that $%
[t_{0}]=\{a,b\}$ and $t_{0}=a<b$. \ In this case, we have $\psi (\theta
_{1},b)=\psi (\theta _{1},a)=1,$ by lemma \ref{nf}, hence $b\leq h(\theta
_{1})=a,$ which is impossible.

Finally, by Proposition~\ref{fg}, $h$ is continuous, therefore $h$ is a
homeomorphism. $\blacksquare $

\begin{remark}
\label{cx} Under the hypothesis of Theorem~\ref{vl}, we have $%
B_{1}^{{}}=\{\theta ^{\prime },\theta ^{\prime \prime }\}$.
\end{remark}

\noindent \textit{Proof}. Suppose that $\theta \in B_{1}^{{}}\setminus
\{\theta ^{\prime },\theta ^{\prime \prime }\}$. Note that $h(\theta
)=t_{0}\in E\setminus \{-\infty ,+\infty \}$. In the proof of Theorem~\ref%
{vl}-step~1, we showed that there exists $\varphi \in B_{2}^{{}}\cup
B_{3}^{{}}$ such that $h(\varphi )=h(\theta )=t_{0}$. We deduce that $\theta
=\varphi $, which is impossible. $\blacksquare $

\begin{remark}
\label{lu} Under the hypothesis of Theorem~\ref{vl}, by step~1, one notes
that if $\theta \in B_{2}^{0}$ (resp. $\theta \in B_{3}^{{}}$) $h(\theta
)\in E_{1}$ (resp. $h(\theta )\in E_{2}$).
\end{remark}

Consider $E=[0,1]\times \{0,1\}$ equipped with the lexicographic order and $%
\Phi :\left[ 0,1\right] \rightarrow \lbrack 0,1]\times \{0,1\}$ the map
defined by $\Phi (t)=(t,1)$, $t\in \lbrack 0,1]$. $\Phi $ is Borel. Indeed,
let $F:([0,1],\tau _{2})\rightarrow ([0,1]\times \{1\},\tau _{0})$ the map
defined by $F(t)=(t,1)$. In the proof of Lemma~\ref{hi} we showed that $F$
is a homeomorphism. On the other hand by Lemma~\ref{cw}, $\mathop{\rm Bor}%
\nolimits([0,1],\tau _{2})=\mathop{\rm Bor}\nolimits([0,1],\tau _{0})$,
hence $\Phi $ is Borel.

Put $\mu =\Phi (m_{1})$, where $m_{1}$ the Lebesgue measure on $[0,1]$.

\begin{corollary}
\label{kv} $h:K\rightarrow (E=\left[ 0,1\right] \times \{0,1\},\mu )$ is a
homeomorphism.
\end{corollary}

\noindent \textit{Proof}. Put $E_{1}=[0,1]\times \{1\}\setminus \left\{
(1,1)\right\} $ and $E_{2}=[0,1]\times \{0\}\setminus \left\{ (0,0)\right\} $%
. Let $b=(t_{0},1)\in E_{1}$. Suppose that $b_{1}=(t_{1},j_{1})\in E\times
\left\{ 0,1\right\} $ and $t_{0}>t_{1}$, $j_{1}\in \left\{ 0,1\right\} $.
One has 
\begin{align*}
\mu (\left] b_{1},b\right] )& =\mu (\left\{ (t,j);\ (t_{1},j_{1})<(t,j)\leq
(t_{0},1\right\} ) \\
& =m_{1}(\left\{ t\in E;(t_{1},j_{1})<(t,1)\leq (t_{0},1)\right\} ) \\
& \geq m_{1}(\left\{ t\in E;t_{1}<t<t_{0}\right\} )>0.
\end{align*}%
If $b_{1}=(t_{1},j_{1})$ and $t_{0}<t_{1}$, by a similar argument, one shows
that $\mu (\left] b,b_{1}\right] )>0$.

Consider $a=(t_{0},0)$. It is clear that $\mu (\left] a,b\right] )=0$, hence 
$\left[ b\right] =\left\{ a,b\right\} $ (by above).

It remains to show that $V=\left\{ (t,j)\in E\times \left\{ 0,1\right\}
;((t,j)>(t_{0},1)\right\} $ is not closed in $E$ (by a similar argument, one
shows that $V=\left\{ (t,j)\in E\times \left\{ 0,1\right\}
;((t,j)<(t_{0},0)\right\} $ is not closed in $E$).

Let $(t_{n})_{n\geq 0}$ be a sequence in $\left] 0,1\right[ $ such that for
every $n$ $\in \mathbb{N}$ $t_{n}>t_{0}$ and $t_{n}\underset{n\rightarrow
+\infty }{\rightarrow }t_{0}$. Since $(t_{n},1)\in V$ for every $n\in 
\mathbb{N}$, it is enough to show that $(t_{n},1)\rightarrow _{n\rightarrow
+\infty }(t_{0},1)\notin V$.

Let $W=\left\{ (u,j)\in E\times \left\{ 0,1\right\} ;
(u,j)<(u_0,j_0)\right\} $ be an open neighborhood of $(t_0,1)$. Since $%
t_0<u_0$ (because $(t_0,1)\in W)$, there exists $n_0\in \mathbb{N}$ such
that $t_n<u_0$ for every $n\geq n_0$, it follows that $(t_n,1)\in W$.

Consider $W=\left\{ (u,j)\in E\times \left\{ 0,1\right\}
;(u,j)>(u_{0},j_{0})\right\} $ an open neighborhood of $(t_{0},1)$.
Obviously $t_{0}\geq u_{0}$, hence $t_{n}>t_{0}\geq u_{0}$. It follows that $%
(t_{n},1)\in W$, for every $n$. By Theorem~\ref{vl}, $h$ is a homeomorphism. 
$\blacksquare $ \bigskip

In the following of this part, one supposes that $(E,\tau _{0})$ is a
connected strongly Lindel\"{o}f space, $-\infty ,+\infty \in E$ and $%
[t]=\{t\}$ for all $t\in E$. One defines the map $\Pi _{1}:E\rightarrow K$
by 
\begin{equation*}
\Pi _{1}(t)=\left\{ 
\begin{array}{c}
\theta _{t},\text{ if }t\in \,]-\infty ,+\infty \lbrack \\ 
\theta ^{\prime },\text{ if }t=+\infty \\ 
\theta ^{\prime \prime },\text{ if }t=-\infty .%
\end{array}%
\right.
\end{equation*}%
By Lemma~\ref{mq}, $\theta _{t}\in B_{2}^{{}}$, and $h(\theta _{t})=t$, $%
\forall t\in E\setminus \{-\infty ,+\infty \}$.

\begin{lemma}
\label{er} The map $\Pi_1$ is Borel.
\end{lemma}

\noindent \textit{Proof}. Consider $V_{a}=\left\{ \theta \in K;\psi (\theta
,a)=1\right\} $ and $V_{a}^{\prime }=\{\theta \in K;\psi (\theta ,a)=0\}$.
It suffices to show ($\Pi _{1})^{-1}(V_{a})$ and $(\Pi
_{1})^{-1}(V_{a}^{\prime })$ are Borel subsets of $E$ for every $a\in
E\setminus \left\{ -\infty ,+\infty \right\} $, because by Proposition~\ref%
{wb} $K$ is a strongly Lindel\"{o}f space (note that $(E,\tau _{0})$ is
connected implies that ($E,\leq )$ is complete and satisfies the
condition~(*) (cf.~\cite[rem.~(d), p.~58]{Kell}). \vskip2pt plus 1pt minus
4pt

Pick $a\in E\setminus \left\{ -\infty ,+\infty \right\} $. One has 
\begin{eqnarray*}
(\Pi _{1})^{-1}(V_{a}) &=&\left\{ t\in E\setminus \left\{ -\infty ,+\infty
\right\} ;\psi (\theta _{t},a)=1\right\} \cup \left\{ +\infty \right\} = \\
\left\{ t\in E\setminus \left\{ -\infty ,+\infty \right\} ;a\leq h(\theta
_{t})\right\} \cup \left\{ +\infty \right\} &=&\left\{ t\in E\setminus
\left\{ -\infty ,+\infty \right\} ;a\leq t\right\} \cup \left\{ +\infty
\right\} .
\end{eqnarray*}%
Thus $(\Pi _{1})^{-1}(V_{a}^{{}})$ is Borel. By a similar argument, one
shows that $(\Pi _{1})^{-1}(V_{a}^{\prime })$ is Borel. $\blacksquare $ %
\vskip2pt plus 1pt minus 4pt

Put $\mu_1=\Pi_1(\mu)$. Note that $\mu_1$ is a measure on $K$. \vskip 2pt
plus 1pt minus 4pt \vskip 2pt plus 1pt minus 4pt

\begin{lemma}
\label{qc} Let $\theta _{0}\in B_{2}^{{}}$ and $\theta _{1}\in B_{3}^{{}}$.
Then:

I) The sets $\left\{ \varphi \in K;\varphi >\theta _{0}\right\} $, $\left\{
\varphi \in K;\text{ }\varphi <\theta _{1}\right\} $ are not closed subsets
in $K$.

II) For every $\varphi _{0}\in B_{2}^{{}}$, there exists $\varphi _{1}\in
B_{3}^{{}}$ such that $[\varphi _{0}]=\{\varphi _{0},\varphi _{1}\}$ and $%
h(\varphi _{0})=h(\varphi _{1})$.
\end{lemma}

\vskip 2pt plus 1pt minus 4pt

\noindent\textit{Proof}. \vskip 2pt plus 1pt minus 4pt

I). We shall show that the set $\left\{ \varphi \in K;\varphi >\theta
_{0}\right\} $ is not closed, the proof is similar for $\left\{ \varphi \in
K;\varphi <\theta _{1}\right\} $. \vskip2pt plus 1pt minus 4pt

\emph{Step 1:} Let $\varphi \in K$ such that $\varphi >\theta _{0}$. Show
that $h(\varphi )>h(\theta _{0})$.

Suppose that $h(\varphi )=h(\theta _{0})$. Since the restriction of $h$ to $%
B_{j}^{{}}$ is injective, $j\in \{2,3\}$, $\varphi \in B_{3}^{{}}$, this is
impossible because $\psi (\varphi ,h(\varphi ))=0$, $\psi (\theta
_{0},h(\theta _{0}))=1$ and $\varphi >\theta _{0}$. It follows that $%
h(\varphi )>h(\theta _{0})$. $\blacksquare $\vskip2pt plus 1pt minus 4pt

\emph{Step 2}: Show that $\left\{ \varphi \in K;\varphi >\theta _{0}\right\} 
$ is not closed in $K$.

By the step~1, $\left\{ \varphi \in K;\varphi >\theta _{0}\right\} =\left\{
\varphi \in K;h(\varphi )>h(\theta _{0})\right\} $. Sup\-pose that this
subset is closed in $K$. We shall show that the set \-$\left\{ t\in
E;t>h(\theta _{0})\right\} $ is closed, in this case the proof will be
finished because $E$ is connected, hence $\left\{ t\in E;t>h(\theta
_{0})\right\} $ is not closed. $\blacksquare $\vskip2pt plus 1pt minus 4pt

By Proposition~\ref{fg}, the map $h:K\rightarrow E$ is continuous and the
set $\left\{ \varphi \in K;h(\varphi )>h(\theta _{0})\right\} $ is a compact
subset of $K$, it follows that $h\left[ \left\{ \varphi \in K;h(\varphi
)>h(\theta _{0})\right\} \right] =\left\{ h(\varphi );h(\varphi )>h(\theta
_{0}),\varphi \in K\right\} $ is compact in $E$. On the other hand, by
Proposition~\ref{omq}, $h$ is surjective, this implies that $\left\{
h(\varphi );h(\varphi )>h(\theta _{0})\right\} =\left\{ t\in E;t>h(\theta
_{0})\right\} $ is compact, hence it is closed in $E$. $\blacksquare $ \vskip%
2pt plus 1pt minus 4pt

II). Let $\varphi _{0}\in B_{2}^{{}}$. By Lemma~\ref{mq}, there exists $%
\varphi _{1}\in B_{3}^{{}}$ such that $h(\varphi _{0})=h(\varphi _{1})$, it
is clear that $\varphi _{0}>\varphi _{1}$ (note that $\psi (\varphi
_{0},h(\varphi _{0}))=1$ and $\psi (\varphi _{1},h(\varphi _{1}))=0)$. On
the other hand 
\begin{align*}
\mu _{1}(\left] \varphi _{1},\varphi _{0}\right] )& =\mu (\left\{ t\in
E;\varphi _{1}<\Pi _{1}(t)\leq \varphi _{0}\right\} ) \\
& \leq \mu (\left\{ t\in E;h(\varphi _{1})\leq h(\Pi _{1}(t))\leq h(\varphi
_{0})\right\} ) \\
& =\mu (\left\{ t\in E;h(\varphi _{1})\leq t\leq h(\varphi _{0})\right\} )=0.
\end{align*}%
Thus $\varphi _{1}\in \lbrack \varphi _{0}]$. \vskip2pt plus 1pt minus 4pt

Let $\varphi \in K\setminus \{\theta ^{\prime },\theta ^{\prime \prime
},\varphi _{0},\varphi _{1}\}$. Let us show that $\varphi \notin \lbrack
\varphi _{0}]$.

\emph{Case 1}, $\varphi >\varphi _{0}$: Since the restriction of $h$ to $%
B_{j}$ is injective, $j=2,3$, $h(\varphi )>h(\varphi _{0})$. Observe that 
\begin{align}
\mu _{1}(\left] \varphi _{0},\varphi \right] )& =\mu (\left\{ t\in E;\varphi
_{0}<\Pi _{1}(t)\leq \varphi \right\} )  \notag \\
& \geq \mu (\left\{ t\in E;h(\varphi _{0})<h(\Pi _{1}(t))<h(\varphi
)\right\} )  \label{hhi} \\
& =\mu (\left\{ t\in E;h(\varphi _{0})<t<h(\varphi )\right\} )>0,  \notag
\end{align}%
(because $\left[ t\right] =\left\{ t\right\} ,$ for all $t\in E)$ hence $%
\varphi \notin \left[ \varphi _{0}\right] $.

\emph{Case 2,} $\varphi <\varphi _{0}$: \ Let us to show that $h(\varphi
)<h_{1}^{{}}(\varphi _{0}).$ Assume that $h(\varphi )=h_{1}^{{}}(\varphi
_{0}),$ hence $h(\varphi )=h(\varphi _{0})=h(\varphi _{1}),$ since the
restriction of $h$ to $B_{j}^{{}}$ is injective, $j=2,3$ and $\varphi \in
B_{0}\cup B_{1}$, either $\varphi =\varphi _{0}$ or $\varphi =\varphi _{1},$
which is impossible. Thus $h(\varphi )<h(\varphi _{0}).$ By a similar
argument, one shows that $\mu (\left] \varphi ,\varphi _{0}\right] )>0$. We
deduce that $\left[ \varphi _{0}\right] =\left\{ \varphi _{0},\varphi
_{1}\right\} $. $\blacksquare $

Recall that $\psi _{\mu _{1},K}(\eta ,\varphi )=f^{\varphi }(\eta ),$ where $%
f^{\varphi }$ is the characteristic function of $\{\theta \in K;\theta
>\varphi \}.$ Let $A_{1}$ be the $C^{\ast }$-subalgebra (with unit)
generated by the family $\{f^{\varphi },\varphi \in K\}$ in $L^{\infty }(\mu
_{1})$. Note that $A_{1}=C(K_{\mu _{1}}(K)).$ By lemmas \ref{kt} and \ref{nf}%
, $h_{K}(\eta )=$sup$\left\{ \varphi \in K;\text{ }\psi _{\mu _{1},K}(\eta
,\varphi )=1\right\} $ exists for every $\eta \in K(K)$ ($\theta
_{K}^{\prime }=$max($K_{\mu _{1}}(K))$ and $\theta _{K}^{\prime \prime }=$%
min($K_{\mu _{1}}(K))$.

\begin{corollary}
\label{gn} The map $h_{\mu _{1},K}:K_{\mu _{1}}(K)\rightarrow K_{\mu }(E)=K$
is a homeomorphism.
\end{corollary}

\noindent \textit{Proof}. Consider the reference set $\widetilde{E}=K$, $%
E_{1}=B_{2}^{\mu ,E}=B_{2}$, $E_{2}=B_{3}^{\mu ,E}=B_{3}$. By Lemma~\ref{qc}
and Theorem~\ref{vl} we obtain that $h_{\mu _{1},K}:K_{\mu
_{1}}(K)\rightarrow K$ is a homeomorphism (note that $\mu _{1}$ satisfies
the condition (\ref{y})). $\blacksquare $ \smallskip

Let $\Pi _{2}:K\rightarrow K_{\mu _{1}}(K)$ the map defined by 
\begin{equation*}
\Pi _{2}(\varphi )=\left\{ 
\begin{array}{c}
\theta _{\varphi ,K}\in B_{2}^{\mu _{1},K},\text{ if }\varphi \in
E_{1}=B_{2}^{\mu ,E}=B_{2} \\ 
\theta _{\varphi ,K}^{\prime }\in B_{3}^{\mu _{1},K},\text{ if }\varphi \in
E_{2}=B_{3}^{\mu ,E}=B_{3} \\ 
\Pi _{2}(\theta _{E}^{\prime })=\theta _{K}^{\prime }\text{ and }\Pi
_{2}(\theta _{E}^{\prime \prime })=\theta _{K}^{\prime \prime },%
\end{array}%
\right.
\end{equation*}%
Where $\theta _{\varphi ,K},\theta _{\varphi ,K}^{\prime }$ the elements
founded by lemma \ref{mq}. Note that $\Pi _{2}(B_{j}^{{}})=B_{j}^{\mu
_{1},K} $, $j\in \{2,3\}$.

\begin{lemma}
\label{mi} The map $\Pi _{2}:K\rightarrow K_{\mu _{1}}(K)$ is Borel.
\end{lemma}

\noindent \textit{Proof}. By Proposition~\ref{wb}, $K_{\mu _{1}}(K)$ is a
strongly Lindel\"{o}f space, it suffices then to show that $(\Pi
_{2})^{-1}(V_{\theta }),(\Pi _{2})^{-1}(V_{\theta }^{\prime })$ are Borel
subsets of $K$, where $V_{\theta }=\left\{ \eta \in K_{\mu _{1}}(K);\psi
_{\mu _{1},K}(\eta ,\theta )=1\right\} $ and $V_{\theta }^{\prime }=\left\{
\eta \in K_{\mu _{1}}(K);\psi _{\mu _{1},K}(\eta ,\theta )=0\right\} $. We
can suppose that $\theta \in B_{2}^{{}}$. We have 
\begin{align*}
(\Pi _{2})^{-1}(V_{\theta })& =\left\{ \varphi \in B_{2}^{{}};\psi _{\mu
_{1},K}(\theta _{\varphi ,K},\theta )=1\right\} \\
& \mskip30mu\cup \left\{ \varphi \in B_{3}^{{}};\psi _{\mu _{1},K}(\theta
_{\varphi ,K}^{\prime },\theta )=1\right\} \cup \left\{ \theta ^{\prime
}\right\} \\
& =\left\{ \varphi \in B_{2}^{{}};\theta \leq h_{\mu _{1},K}(\theta
_{\varphi ,K}))\right\} \\
& \mskip30mu\cup \left\{ \varphi \in B_{3}^{{}};\theta <h_{\mu
_{1},K}(\theta _{\varphi ,K}^{\prime })\right\} \cup \left\{ \theta ^{\prime
}\right\} \\
& =\left\{ \varphi \in B_{2}^{{}};\theta \leq \varphi \right\} \cup \left\{
\varphi \in B_{3}^{{}};\theta <\varphi \right\} \cup \left\{ \theta ^{\prime
}\right\} \\
& =\left\{ \varphi \in K;\theta \leq \varphi \right\} .
\end{align*}%
Thus $(\Pi _{2})^{-1}(V_{\theta })$ is a Borel subset of $K$. By a similar
argument, one shows that $(\Pi _{2})^{-1}(V_{\theta }^{\prime })$ is a Borel
subset of $K$. $\blacksquare $

Let $\mu_2=\Pi_2(\mu_1)$.

\begin{lemma}
\label{nb} Suppose that $\mu (]t\wedge t^{\prime },t\vee t^{\prime }[)>0$,
for every $t\neq t^{\prime }$. Let $\eta _{0}\in B_{2}^{K}$ and $\eta
_{1}\in B_{3}^{K}$. Then:

I) The sets $\left\{ \eta \in K_{\mu _{1}}(K);\eta >\eta _{0}\right\} $, $%
\left\{ \eta \in K_{\mu _{1}}(K);\eta <\eta _{1}\right\} $ are not closed
subsets of $K_{\mu _{1}}(K)$.

II) For every $\lambda _{0}\in B_{2}^{\mu _{1},K}$, there exists $\lambda
_{1}\in B_{3}^{\mu _{1},K}$ such that $\left[ \lambda _{0}\right] =\left\{
\lambda _{0},\lambda _{1}\right\} $ and $h(h_{\mu _{1},K}(\lambda
_{0}))=h(h_{\mu _{1},K}(\lambda _{1}))$.
\end{lemma}

\noindent\textit{Proof}.

I). Let us show that $\left\{ \eta \in K_{\mu _{1}}(K);\eta >\eta
_{0}\right\} $ is not closed. Suppose that there exists $\eta \in K_{\mu
_{1}}(K)$ such that $\eta >\eta _{0}$. There exists $\theta \in K$ such that 
$\psi _{\mu _{1},K}(\eta ,\theta )=f^{\theta }(\eta )=1$ and $\psi _{\mu
_{1},K}(\eta _{0},\theta )=f^{\theta }(\eta _{0})=0$. Since $\eta _{0}\in
B_{2}^{\mu _{1},K}$, $h_{\mu _{1},K}(\eta )\geq \theta >h_{\mu _{1},K}(\eta
_{0})$. It follows that $h_{K}^{{}}(\eta )>h_{K}^{{}}(\eta _{0})$. We
conclude that $\left\{ \eta \in K(K);\eta >\eta _{0}\right\} =\left\{ \eta
\in K(K);h_{\mu _{1},K}(\eta )>h_{\mu _{1},K}(\eta _{0})\right\} $.

Assume now that $\left\{ \eta \in K_{\mu _{1}}(K);\eta >\eta _{0}\right\}
=\left\{ \eta \in K_{\mu _{1}}(K);h_{\mu _{1},K}^{{}}(\eta )>h_{\mu
_{1},K}(\eta _{0})\right\} $ is closed. The map $h_{\mu _{1},K}^{{}}:K_{\mu
_{1}}(K)\rightarrow K$ is a homeomorphism by Corollary~\ref{gn}, it follows
that 
\begin{equation*}
h_{\mu _{1},K}^{{}}\left[ \left\{ \eta \in K_{\mu _{1}}(K);h_{\mu
_{1},K}(\eta )>h_{\mu _{1},K}(\eta _{0})\right\} \right] =\left\{ h_{\mu
_{1},K}(\eta );h_{K}(\eta )>h_{\mu _{1},K}(\eta _{0}),\eta \in K_{\mu
_{1}}(K)\right\}
\end{equation*}%
is a compact subset of $K$. Since $h_{\mu _{1},K}$ is onto, the set 
\begin{equation*}
\left\{ h_{\mu _{1},K}(\eta )\in K;h_{\mu _{1},K}(\eta )>h_{\mu _{1},K}(\eta
_{0}),\eta \in K_{\mu _{1}}(K)\right\} =\left\{ \varphi \in K;\varphi
>h_{\mu _{1},K}(\eta _{0})\right\}
\end{equation*}%
is closed in $K$. On the other hand, by Remark~\ref{lu}, $h_{\mu
_{1},K}(\eta _{0})\in B_{2}^{{}}$, using Lemma~\ref{qc}-I) we obtain that $%
\left\{ \varphi \in K;\varphi >h_{\mu _{1},K}(\eta _{0})\right\} $ is not
closed, which is impossible. Thus $\left\{ \eta \in K_{\mu _{1}}(K);\eta
>\eta _{0}\right\} $ is not closed. By a similar argument, one shows that
the set $\left\{ \eta \in K_{\mu _{1}}(K);\eta <\eta _{1}\right\} $ is not
closed. $\blacksquare $ \vskip2pt plus 1pt minus 4pt \vskip2pt plus 1pt
minus 4pt

II). Let $\lambda _{0}\in B_{2}^{\mu _{1},K}$. Since $h_{3}^{{}}:B_{3}^{{}}%
\rightarrow E\setminus \{-\infty ,+\infty \}$ is onto, there exists $\theta
_{1}\in B_{3}^{{}}$ such that%
\begin{equation}
h(\theta _{1})=h(h_{\mu _{1}K}(\lambda _{0})).  \label{L}
\end{equation}%
On the other hand, by Corollary~\ref{gn} and Remark~\ref{lu}, there exists $%
\lambda _{1}\in B_{3}^{\mu _{1},K}$ such that $h_{\mu _{1},K}(\lambda
_{1})=\theta _{1},$ this implies by \ref{L} that%
\begin{equation}
h(h_{\mu _{1},K}(\lambda _{0}))=h(h_{\mu _{1},K}(\lambda _{1})).  \label{MT}
\end{equation}

We shall prove that $\lambda _{1}<\lambda _{0}$. \vskip2pt plus 1pt minus 4pt

Suppose that $\lambda _{1}\geq \lambda _{0}$. It follows that $h_{\mu
_{1},K}(\lambda _{1})=\theta _{1}\geq h_{\mu _{1},K}(\lambda _{0})$. Since $%
\theta _{1}\in B_{3}^{{}}$, $\psi (\theta _{1},h(\theta _{1}))=0$, hence $%
\psi _{\mu _{1},K}(h_{\mu _{1},K}(\lambda _{0}),h(\theta _{1}))=0$, this
implies that $h(h_{\mu _{1},K}(\lambda _{0}))<h(\theta _{1})$ because by
Corollary~\ref{gn} and Remark~\ref{lu}, $h_{\mu _{1},K}(\lambda _{0})\in
B_{2}^{{}}$, it is impossible by (\ref{L}) and (\ref{MT}). \vskip2pt plus
1pt minus 4pt

Thus $\lambda_1<\lambda_0$. \vskip 2pt plus 1pt minus 4pt

Let us show that $\lambda_1\in \left[ \lambda_0\right] $. \vskip 2pt plus
1pt minus 4pt

Observe that \vskip2pt plus 1pt minus 4pt 
\begin{eqnarray*}
\mu _{2}(\left] \lambda _{0},\lambda _{1}\right] ) &=&\mu _{2}(\left\{ \eta
\in K_{\mu _{1}}(K);\lambda _{0}<\eta \leq \lambda _{1}\right\} ) \\
\noalign{\vskip3pt plus 2pt minus 2pt} &=&\mu _{1}(\left\{ \theta \in
K;\lambda _{0}<\Pi _{2}(\theta )\leq \lambda _{1}\right\} ) \\
\noalign{\vskip3pt plus 2pt minus 2pt} &\leq &\mu _{1}(\left\{ \theta \in
K;h_{\mu _{1},K}(\lambda _{0})\leq h_{\mu _{1},K}(\Pi _{2}(\theta ))\leq
h_{\mu _{1},K}(\lambda _{1})\right\} ) \\
\noalign{\vskip3pt plus 2pt minus 2pt} &=&\mu _{1}(\left\{ \theta \in
K;h_{\mu _{1},K}(\lambda _{0})\leq \theta \leq h_{\mu _{1},K}(\lambda
_{1})\right\} ) \\
\noalign{\vskip3pt plus 2pt minus 2pt} &=&\mu (\left\{ t\in E;h_{\mu
_{1},K}(\lambda _{0})\leq \theta _{t}\leq h_{\mu _{1},K}(\lambda
_{1})\right\} ) \\
\noalign{\vskip3pt plus 2pt minus 2pt} &\leq &\mu (\left\{ t\in E;h(h_{\mu
_{1},K}(\lambda _{0}))\leq h(\theta _{t})\leq h(h_{\mu _{1},K}(\lambda
_{1}))\right\} ) \\
\noalign{\vskip3pt plus 2pt minus 2pt} &=&\mu (\left\{ t\in E;h(h_{\mu
_{1},K}(\lambda _{0}))\leq t\leq h(h_{\mu _{1},K}(\lambda _{1}))\right\} )=0.
\end{eqnarray*}%
\vskip2pt plus 1pt minus 4pt \vskip2pt plus 1pt minus 4pt \noindent We
conclude that $\lambda _{0}\in \left[ \lambda _{1}\right] $. Let us show
that $\eta \in K_{\mu _{1}}(K)\setminus \{\theta _{K}^{\prime },\theta
_{K}^{\prime \prime },\lambda _{0},\lambda _{1}\}$. We show that $\eta
\notin \left[ \lambda _{0}\right] $. \vskip2pt plus 1pt minus 4pt \vskip2pt
plus 1pt minus 4pt

\emph{Case 1:} $\eta >\lambda_0$. \vskip 2pt plus 1pt minus 4pt \vskip 2pt
plus 1pt minus 4pt

\emph{Step 1:} We shall show that $h(h_{\mu _{1},K}(\eta ))>h(h_{\mu
_{1},K}(\lambda _{1}))=h(h_{\mu _{1},K}(\lambda _{0}))$.

Indeed, suppose that $h(h_{\mu _{1},K}(\eta ))=h(h_{\mu _{1},K}(\lambda
_{1}))$ $=h(h_{\mu _{1},K}(\lambda _{0}))$. By Corollary~\ref{gn} and Remark~%
\ref{lu}, $h_{\mu _{1},K}(\lambda _{0})\in B_{2}^{{}}$ and $h_{\mu
_{1},K}(\lambda _{1})\in B_{3}^{{}}$. We distinguish two cases: \vskip2pt
plus 1pt minus 4pt

\emph{Case 1-a,} $\eta \in B_{2}^{\mu _{1},K}$: By Corollary~\ref{gn} and
Remark~\ref{lu}, $h_{\mu _{1},K}(\eta )\in B_{2}^{{}}$. On the other hand,
the restriction of $h$ to $B_{2}^{{}}$ is injective, hence $h_{\mu
_{1},K}(\eta )=h_{\mu _{1},K}(\lambda _{0})$. This is impossible because the
restriction of $h_{\mu _{1},K}$ to $B_{2}^{\mu _{1},K}$ is injective by
Remark~\ref{cnn}. \vskip2pt plus 1pt minus 4pt

\emph{Case 1-b,} $\eta \in B_{3}^{\mu _{1},K}$: Since $h_{\mu _{1},K}(\eta
)\in B_{3}^{{}}$ and the restriction of $h$ to $B_{3}^{{}}$ is injective, $%
h_{\mu _{1},K}(\eta )=h_{\mu _{1},K}(\lambda _{1})$, this implies that $\eta
=\lambda _{1}$, this is impossible (we have shown previously that $\lambda
_{0}>\lambda _{1})$. We conclude that $h(h_{K}^{{}}(\eta ))>h(h_{\mu
_{1},K}(\lambda _{1}))=h(h_{\mu _{1},K}(\lambda _{0}))$. $\blacksquare $%
\vskip2pt plus 1pt minus 4pt

\emph{Step2:} Let us show that $\mu _{2}(\left] \lambda _{0},\eta \right[
)>0 $. \medskip

We have \vskip2pt plus 1pt minus 4pt \vskip2pt plus 1pt minus 4pt 
\begin{eqnarray*}
\mu _{2}(\left] \lambda _{0},\eta \right[ ) &=&\mu _{2}(\left\{ \omega \in
K_{\mu _{1}}(K));\lambda _{0}<\omega <\eta \right\} ) \\
\noalign{\vskip4pt minus2pt} &\geq &\mu _{2}(\left\{ \omega \in K_{\mu
_{1}}(K);h_{K}(\lambda _{0})<h_{\mu _{1},K}(\omega )<h_{\mu _{1},K}(\eta
)\right\} ) \\
\noalign{\vskip4pt minus2pt} &=&\mu _{1}(\left\{ \theta \in K;h_{\mu
_{1},K}(\lambda _{0})<h_{\mu _{1},K}(\Pi _{2}(\theta ))<h_{\mu _{1},K}(\eta
)\right\} ) \\
\noalign{\vskip4pt minus2pt} &=&\mu _{1}(\left\{ \theta \in K_{\mu
_{1}}(K);h_{\mu _{1},K}(\lambda _{0})<\theta <h_{\mu _{1},K}(\eta )\right\} )
\\
&\geq &\mu _{1}(\left\{ \theta \in K_{\mu _{1}}(K);h_{K}^{{}}(h_{\mu
_{1},K}(\lambda _{0}))<h(\theta )<h(h_{\mu _{1},K}(\eta ))\right\} ) \\
\noalign{\vskip4pt minus2pt} &=&\mu (\left\{ t\in E;h(h_{\mu _{1},K}(\lambda
_{0}))<h(\theta _{t})\leq h(h_{\mu _{1},K}(\eta ))\right\} ) \\
\noalign{\vskip4pt minus2pt} &=&\mu (\left\{ t\in E;h(h_{\mu _{1},K}(\lambda
_{0}))<t<h(h_{\mu _{1},K}(\eta ))\right\} )>0.\text{ }\blacksquare
\end{eqnarray*}%
\vskip2pt plus 1pt minus 4pt \noindent Thus $\eta \notin \left[ \lambda _{0}%
\right] $. \vskip2pt plus 1pt minus 4pt

\emph{Case 2, }$\eta <\lambda _{0}$: By a similar argument, one shows that $%
\eta \notin \lbrack \lambda _{0}]$. It follows that $[\lambda
_{0}]=\{\lambda _{0},\lambda _{1}\}$. $\blacksquare $

Let $A_{2}$ be the $C^{\ast }$-subalgebra which is generated by the family $%
\{f^{\eta };\eta \in K_{\mu _{1}}(K)\}$ (where $f^{\eta }$ the
characteristic function of $\{\xi \in K_{\mu _{1}}(K);\xi >\eta \}$) in $%
L^{\infty }(K_{\mu _{1}}(K),\mu _{2}),$ ($A_{2}=C(K_{\mu _{2}}(K_{\mu
_{1}}(K)))$).

Recall that $B_{1}^{\mu _{2},K_{\mu _{1}}(K)}$ is the set of $\xi \in K_{\mu
_{2}}(K_{\mu _{1}}(K))$ such that the map $\eta \in (K_{\mu
_{1}}(K))\setminus \{\theta _{K}^{\prime },\theta _{K}^{\prime \prime
}\}\rightarrow \psi _{\mu _{2},K_{\mu _{1}}(K)}(\xi ,\eta )=\xi (f^{\eta })$
is continuous and $B_{j}^{\mu _{2},K_{\mu _{1}}(K)}=\left\{ \xi \in K_{\mu
_{2}}(K_{\mu _{1}}(K))\setminus B_{1}^{\mu _{2},K_{\mu _{1}}(K)};\psi _{\mu
_{1},K}(\xi ,h_{K_{\mu _{1}}(K)}(\xi ))=3-j\right\} $, $j=2,3,$ where $%
h_{\mu _{2},K_{\mu _{1}}(K)}^{{}}(\xi )=\sup \left\{ \eta \in K_{\mu
_{1}}(K);\psi _{\mu _{2},K_{\mu _{1}}(K)}(\xi ,\eta )=1\right\} $.

\begin{corollary}
\label{wy} The map $h_{\mu _{2},K_{\mu _{1}}(K)}^{{}}:K_{\mu _{2}}(K_{\mu
_{1}}(K))\rightarrow K_{\mu _{1}}(K)$ is a homeomorphism.
\end{corollary}

\noindent \textit{Proof}. Consider the reference set $\widetilde{E}=K_{\mu
_{1}}(K)$, ($E)_{1}=B_{2}^{\mu _{2},K_{\mu _{1}}(K)}$, $E_{2}=B_{3}^{\mu
_{2},K_{\mu _{1}}(K)}$. By Theorem~\ref{vl} and Lemma~\ref{nb}, $h_{\mu
_{2},K_{\mu _{1}}(K)}^{{}}:K_{\mu _{2}}(K_{\mu _{1}}(K))\rightarrow K_{\mu
_{1}}(K)$ is a homeomorphism. $\blacksquare $

\section{General properties of ($C(L),\protect\tau _{p})$ when L is a
Rosenthal compact set}

In this part, we show that the funtion $\psi :(K\times (E,\tau
_{0}))\rightarrow \left\{ 0,1\right\} $ is not Borel under some assumptions.
Let $L$ be Hausdorff compact. In \cite{Tal(2)} one shows under the Martin's
Axiom that $(C(L),\tau _{p})$ is universally measurable. During this part,
we show that $(C(L),\tau _{p})$ is universally measurable without Martin's
Axioms, if $L$ is a Rosenthal compact set. Finally, we show that if $K$ is
separable, then $(C(K,\tau _{p})$ is not measure-compact.

\begin{definition}
\label{ae} Let $(X, \tau)$ be a Hausdorff space. We say that $(X, \tau)$ is
a Radon space \cite[p.~117]{Sch(1)}, if for every positive finite Borel
measure on $(E, \tau)$ is a Radon measure.
\end{definition}

\begin{definition}
\label{DR} Let $(X,\tau )$ be a Hausdorff space. A positive finite Borel
measure $\nu $ on $(X,\tau )$ is said to be normal if for every net $%
(O_{\alpha })$ of open sets is increasing to open set $O$, $\underset{\alpha 
}{\lim }\,\nu (O_{\alpha })=\nu (O)$.
\end{definition}

\begin{remark}
\label{hht} Let $(X, \tau)$ be a strongly Lindel\"{o}f space. Then every
positive Borel measure on $(X, \tau)$ is normal.
\end{remark}

\begin{remark}
\label{p} Let $(X,\tau )$ be a Hausdorff topological space and let $%
X_{1}\subset X$. Then $\mathop{\rm Bor}\nolimits(X_{1},\tau )=\left\{ C\cap
X_{1};C\in \mathop{\rm Bor}\nolimits(X,\tau )\right\} $.
\end{remark}

\begin{remark}
\label{mmm} Let $L$ be a compact strongly Lindel\"{o}f space. Then $K$ is a
Radon space.
\end{remark}

\noindent \textit{Proof}. As $L$ is a strongly Lindel\"{o}f space, by Remark~%
\ref{hht}, every positive finite Borel measure on $K$ is normal. On the
other hand, $K$ is universally measurable, by \cite[th.~3.2]{Sch(2)} every
normal positive finite Borel measure on $K$ is a Radon measure, hence $L$ is
a Radon-space. $\blacksquare $

\begin{proposition}
\label{ynm} Let $(E,\tau _{0})$ be a strongly Lindel\"{o}f space such that $%
-\infty $, $+\infty \notin E$; let $E^{\prime }$ be a subset of $E$ such
that $(E^{\prime },\tau _{1})$ be completely regular space. Assume that
there exists a measure $\nu $ strictly positive on $(E^{\prime },\tau _{0})$
such that $\nu $ is diffuse on $E$. Then $(E^{\prime },\tau _{1})$ is not
universally measurable space.
\end{proposition}

\noindent \textit{Proof}. By Lemma~\ref{cw}, $\mathop{\rm Bor}%
\nolimits(E,\tau )=\mathop{\rm Bor}\nolimits(E,\tau _{1})$. This implies by
Remark~\ref{p} that $\mathop{\rm Bor}\nolimits(E^{\prime },\tau _{0})=%
\mathop{\rm
Bor}\nolimits(E^{\prime },\tau _{1})$. Hence, there exists $V\in 
\mathop{\rm
Bor}\nolimits(E^{\prime },\tau _{1})$ such that $0<\nu (V)<+\infty $. One
defines the measure $\nu _{1}$ on $E^{\prime }$ by, $\nu _{1}(U)=\nu (U\cap
V)$, $U\in \mathop{\rm Bor}\nolimits(E^{\prime },\tau _{1})$. $\nu _{1}$ is
a finite Borel measure on $(E^{\prime },\tau _{1})$. On the other hand, by
Lemma~\ref{cw}, $(E,\tau _{1})$ is a strongly Lindel\"{o}f space, hence $%
(E^{\prime },\tau _{1})$ has the same property by Remark~\ref{nl}, we
conclude by Remark~\ref{hht} that $\nu _{1}$ is normal. Suppose now that $%
(E^{\prime },\tau _{1})$ is universally measurable. By \cite[th.~3.2]{Sch(2)}%
, $\nu _{1}$ is a Radon measure. It follows that there exists a compact set $%
L_{0}$ of $(E^{\prime },\tau _{1})$ such that $\nu _{1}(L_{0})>0$. Put $%
x_{0}=\max \left\{ t;t\in L_{0}\right\} $. By hypothesis on $\nu $, we have $%
\nu _{1}(L_{0}-\left\{ x_{0}\right\} )>0$, then there exists a compact set $%
L_{1}$ of $L_{0}-\left\{ x_{0}\right\} $ satisfying $\nu _{1}(L_{1})>0$.
Denote $x_{1}=\max \left\{ t;t\in L_{1}\right\} $. By induction one
constructs a strictly decreasing sequence $(x_{n})_{n\geq 0}$ of $L_{0}$.
Let $\mathcal{U}$ be a non trivial ultrafilter on $\mathbb{N}$. Since $L_{0}$
is compact, there exists $a\in L_{0}$ such that $x_{n}\underset{\mathcal{U}}{%
\rightarrow }a$ in $(E^{\prime },\tau _{1})$. Let $O=\left\{ t\in E^{\prime
};t\leq a\right\} $. $O$ is an open neighborhood of $a$, this implies that $%
\left\{ n\in \mathbb{N};x_{n}\in O\right\} =\left\{ n\in \mathbb{N}%
;x_{n}\leq a\right\} \in \mathcal{U}$. It is impossible because $x_{k}\geq a$%
, for every $k\in \mathbb{N}$ and the sequence $(x_{n})_{n\geq 0}$ is
strictly decreasing. $\blacksquare $

By Remark~\ref{cs} and Proposition~\ref{ynm}, one has the following
corollary:

\begin{corollary}
\label{ui} Let ($E, \tau_0)$ be a strongly Lindel\"{o}f space. Assume that $%
\nu $ strictly positive on $(E, \tau_0)$ and diffuse on $E$, $(E, \leq)$ is
complete and $-\infty, +\infty \notin E$. Then $(E, \tau_1)$ is not
universally measurable.
\end{corollary}

Recall that 
\begin{equation*}
E_2 = \left\{ a \in E; \left\{ x \in E; x < a \right\} \text{ is not a
closed subset of } E\right\}.
\end{equation*}

\begin{corollary}
\label{qj} Assume that $(E,\tau _{0})$ is a strongly Lindel\"{o}f space,
that $-\infty ,+\infty \notin E$, that $E_{2}=E$, the measure $\mu $ is
diffuse on $E$ and that $\mu (]t\wedge t^{\prime },t\vee t^{\prime }])>0$
for $t\neq t^{\prime }$. Then $B_{3}^{{}}$ is not universally measurable.
\end{corollary}

\noindent \textit{Proof}. By Remark~\ref{wj}, $(E,\tau _{1})$ is completely
regular. On the other hand Proposition~\ref{ynm} shows us that $(E,\tau
_{1}) $ is not universally measurable and by step~2 of Proposition~\ref{ee}%
-I), $B_{3}^{{}}$ is homeomorphic to $(h(B_{3}),\tau _{1})$. But $%
h(B_{3}^{{}})=E$ by Lemma~\ref{mq}, we conclude that $B_{3}^{{}}$ is not
universally measurable. $\blacksquare $ \vskip2pt plus 1pt minus 4pt

\begin{remark}
\label{mms} In Corollary~\ref{qj}, we can replace the hypothesis $E_{2}=E,$
by the hypothesis $h_{3}^{{}}:B_{3}^{{}}\rightarrow E$ is onto.
\end{remark}

\vskip 2pt plus 1pt minus 4pt

\begin{corollary}
\label{hh} Let $(E, \tau_0)$ be a strongly Lindel\"{o}f space, and suppose
that $-\infty, +\infty \notin E$. Assume that there exists a diffuse measure 
$\nu $ on $E$ and there exists a compact set $L_0$ of $(E, \tau_0)$
verifying $0 < \nu(L_0) < +\infty$. Then $(E, \tau_1)$ does not embed in any
analytic set.
\end{corollary}

\vskip 2pt plus 1pt minus 4pt

\noindent \textit{Proof}. Consider the identity map $i:(E,\tau
_{0})\rightarrow (E,\tau _{1})$. By Lemma~\ref{cw}, $i$ is Borel. Suppose
now that there exists an analytic set $Y$ such that $(E,\tau _{1})$ is a
topological subspace of $Y$. Consider the measure $\nu _{1}$ defined on $%
L_{0}$ by $\nu _{1}(B)=\nu (B)$, $B\in \mathop{\rm Bor}\nolimits(L_{0})$.
Since $(L_{0},\tau _{0})$ is a strongly Lindel\"{o}f space, $\nu _{1}$ is
normal. On the other hand $(L_{0},\tau _{0})$ is universally measurable, by 
\cite[th.3.2]{Sch(2)} $\nu _{1}$ is a Radon measure. By \cite[th.~14]{Sch(3)}%
, there is a compact subset $L_{1}$ of $L_{0}$ of measure strictly positive
such that the restriction of $i$ to $L_{1}$ is continuous (with values in $%
Y) $. It follows that $i(L_{1})=L_{1}$ is compact in $(E,\tau _{1})$ and $%
\nu _{1}(L_{1})>0$. using Proposition~\ref{ynm} for $E^{\prime }=L_{1}$, we
obtain that $L_{1}$ is not universally measurable, which is impossible. $%
\blacksquare $ \vskip2pt plus 1pt minus 4pt

\begin{corollary}
\label{rr} Let $(E,\tau _{0})$ be a strongly Lindel\"{o}f such that $-\infty
,+\infty \notin E$. Suppose that $h_{3}^{{}}$ is onto and the hypothesis on $%
\nu $ in corollary~\ref{hh} are satisfied for $\nu =\mu $. Then $B_{3}^{{}}$
does not embed in any analytic set.
\end{corollary}

\vskip 2pt plus 1pt minus 4pt

\noindent \textit{Proof}. By Proposition~\ref{ee}-I) $B_{3}^{{}}$ is
homeomorphic to $(E,\tau _{1})$. For getting the corollary it suffices to
use Corollary~\ref{hh}. \vskip2pt plus 1pt minus 4pt

\begin{proposition}
\label{av} Suppose that $(E,\tau _{0})$ is separable, satisfying the
condition~$(\ast )$, and the measure $\mu $ is diffuse on $E$, $%
h_{3}^{{}}:B_{3}^{{}}\rightarrow E$ is surjective, $-\infty ,+\infty \notin
E $ and $\mu (\left] t\wedge t^{\prime },t\vee t^{\prime }\right] )>0$ for $%
t\neq t^{\prime }$. Then the map $\psi :(K\times (E,\tau _{0}))\rightarrow
\{0,1\}$ is not Borel.
\end{proposition}

\vskip 2pt plus 1pt minus 4pt

\noindent\textit{Proof}.

\emph{Step 1:} Let us show that $\mathop{\rm Bor}\nolimits(K) \otimes %
\mathop{\rm Bor}\nolimits(E) = \mathop{\rm Bor}\nolimits(K\times E)$.

It is obvious that $\mathop{\rm Bor}\nolimits(K)\otimes \mathop{\rm Bor}%
\nolimits(E)\subset \mathop{\rm Bor}\nolimits(K\times E)$. We show the
reverse inclusion. Let $(O_{k})_{k\geq 0}$ be a countable basis of $(E,\tau
_{0})$ (a such basis exists by Remark~\ref{gd}) and let $V=\underset{i\in I}{%
\cup }U_{i}\times W_{i}$, where the $U_{i}$ are open subsets of $K$ and the $%
W_{i}$ are open subsets of $E$. For every $i\in I$, there exists a subset $%
M_{i}$ of $\mathbb{N}$ such that $W_{i}=\underset{k\in M_{i}}{\cup }O_{k}$,
hence $V=\underset{k\in \mathbb{N}}{\cup }\left[ (\underset{i\in Y_{k}}{\cup 
}U_{i})\times O_{k})\right] $, where $Y_{k}=\left\{ i;k\in M_{i}\right\} $.
Since $\underset{i\in Y_{k}}{\cup }U_{i}$ is an open subset of $K$, then 
\begin{equation*}
V\in \mathop{\rm Bor}\nolimits(K)\otimes \mathop{\rm Bor}\nolimits(E).
\end{equation*}%
It follows that $\mathop{\rm Bor}\nolimits(K\times E)\subset \mathop{\rm Bor}%
\nolimits(K)\otimes \mathop{\rm Bor}\nolimits(E)$. $\blacksquare $\smallskip %
\vskip2pt plus 1pt minus 4pt

\emph{Step 2:} Let us show that every Borel subset of $K$ is universally
measu\-rable.

By Corollary~\ref{mm}, $K$ is a strongly Lindel\"{o}f space, using Remark~%
\ref{mmm}, one obtains that $K$ is a Radon space. We conclude that every
Borel subset of $K$ is relatively universally measurable in $K$, by \cite[%
th.~3.2]{Sch(2)} every Borel subset of $K$ is universally measu\-rable. $%
\blacksquare $ \vskip2pt plus 1pt minus 4pt

Suppose now that $\psi $ is Borel. Consider the map $\sigma :K\setminus
B_{1}^{{}}\rightarrow K\times E$ defined by $\sigma (\theta )=(\theta
,h(\theta ))$, $\theta \in K\setminus B_{1}^{{}}$. This map $\sigma $ is
Borel because $\mathop{\rm Bor}\nolimits(K)\otimes \mathop{\rm Bor}%
\nolimits(E)=\mathop{\rm Bor}\nolimits(K\times E)$ and $h$ is continuous by
Proposition~\ref{cc}. It follows that $\psi \circ \sigma $ is Borel, hence $%
(\psi \circ \sigma )^{-1}\{0\}=B_{3}^{{}}$ is Borel subset of $K$. By the
step~2, $B_{3}^{{}}$ is universally measurable, this is impossible by Remark~%
\ref{mms}. $\blacksquare $ \vskip2pt plus 1pt minus 4pt

Let $(Y, \tau)$ be a Hausdorff topological space. We denote by $%
\mathop{\rm
Bair}\nolimits(Y)$ the $\sigma$-algebra generated by the continuous
functions on $Y$ and by $P_\sigma (Y)$ the set of probability measures on $%
(Y, \mathop{\rm Bair}\nolimits(Y))$. \bigskip

\begin{definition}
\label{nln} (\cite{Mor1}-\cite{Mor2}-\cite{Var}) Let $X$ be a completely
regular space:

a) $X$ is measure-compact, if for every measure $\nu \in P_\sigma (X)$ and
every net ($f_\alpha )_{\alpha \in I}$ of bounded continuous decreasing
functions on $(X, \tau)$ such that $f_\alpha \rightarrow 0$, $\underset{%
\alpha }{\lim} \int_{X} f_\alpha d\nu =0$.

b) $X$ is strongly measure-compact, if for every measure $\nu \in P_{\sigma
}(X)$ and every $\varepsilon >0$, there is a compact $L$ of $X$ such that $%
\nu^*(L)>1-\varepsilon $, where 
\begin{equation*}
\nu^*(L) = \inf \left\{ \nu (H); H \in \mathop{\rm Bair}\nolimits(X)\text{
and } L \subset H\right\} .
\end{equation*}
\end{definition}

\vskip 2pt plus 1pt minus 4pt

Note that if $X$ is strongly measure-compact, then $X$ is measure-compact.

\begin{remark}
\label{xx} Suppose that $(E, \leq)$ satisfies the condition~$(*)$. Then $(E,
\tau_0)$ and $(E, \tau_1)$ are regular spaces.
\end{remark}

\noindent \textit{Proof}. Let us show for example that $(E,\tau _{1})$ is
regular. Indeed, let $V=\left] a,b\right] $ an open neighborhood of $x$.
Since $(E,\leq )$ verifies the condition~$(\ast )$, there is $c\in E$ such
that $a<c<x\leq b$, hence $[c,b]=\overline{\left] c,b\right] }\subset V$ and 
$]c,b]$ is an open neighborhood of $x$. We can use the same argument if $%
V=]a,b[$, $a<b$. $\blacksquare $ \bigskip

\begin{proposition}
\label{nh} Let $(E, \tau_0)$ be a separable strongly Lindel\"{o}f space
satisfying the condition~$(*)$. Suppose there is a probability measure $\nu $
on $(E, \tau_0)$, diffuse on $E$ and $-\infty $, $+\infty \notin E$. Then $%
(E, \tau_1)$ is not strongly measure-compact.
\end{proposition}

\noindent\textit{Proof}. Since $E$ satisfies $(*)$, $E = E_2$. By Remark~\ref%
{wj}, $(E, \tau_1)$ is completely regular.

\emph{Step 1:} Show that every open subset of $(E, \tau_1)$ is in $%
\mathop{\rm Bair}\nolimits(E, \tau_1)$.

Pick $a\in E$. The set $\left\{ t\in E; t\leq a\right\} $ is open and closed
in $(E, \tau_1)$, hence the characteristic function of $\left\{ t\in E;
t\leq a\right\} $ is continuous on $(E, \tau_1)$. It follows that $\left\{
t\in E; t\leq a\right\} \in \mathop{\rm Bair}\nolimits(E, \tau_1)$.

Since $(E,\tau _{0})$ is strongly Lindel\"{o}f and regular (by Remark~\ref%
{xx}), by \cite[chap.~II, prop.~1]{Sch(1)}, $\{a\}$ is $G_{\delta }$ of $%
(E,\tau _{0})$. Thus there exists a sequence of open subsets $(O_{n})_{n\geq
0}$ of $(E,\tau _{0})$ such that $\left\{ a\right\} =\underset{n\geq 0}{\cap 
}O_{n}$. For every $n\geq 0$, there is an open interval $I_{n}$ of $O_{n}$
such that $a\in I_{n}\subset O_{n}$, this implies that there exists $%
a_{n}\in E$ such that $\left] a_{n},a\right] \subset I_{n}\subset O_{n}$.
Therefore $\left\{ a\right\} =\underset{n\geq 0}{\cap }\left] a_{n},a\right] 
$.

On the other hand, $\left] a_{n},a\right] =\left\{ t\in E;t\leq a\right\}
-\left\{ t\in E;t\leq a_{n}\right\} \in \mathop{\rm Bair}\nolimits(E,\tau
_{1})$. It follows that $\left\{ a\right\} \in \mathop{\rm Bair}%
\nolimits(E,\tau _{1})$. Thus for every $a\in E$, $\left\{ t\in
E;t<a\right\} =\left\{ t\in E;t\leq a\right\} -\left\{ a\right\} \in %
\mathop{\rm Bair}\nolimits(E,\tau _{1})$.

Since $(E,\tau _{1})$ is strongly Lindel\"{o}f, by Lemma~\ref{cw}, every
open subset of $(E,\tau _{1})$ is in $\mathop{\rm Bair}\nolimits(E,\tau
_{1}) $. $\blacksquare $

\emph{Step 2:} Let $L$ be a compact subset of $(E,\tau _{1})$. Let us show
that $\nu (L)=0$.

Sup\-pose $\nu (L)>0$. By Proposition~\ref{ynm}, $L$ is not universally
measurable, which is impossible. $\blacksquare $

\emph{Step 3}: Let us show that $(E,\tau _{1})$ is not strongly
measure-compact.

Note first by step~1 that $\mathop{\rm Bor}\nolimits(E,\tau _{1})=%
\mathop{\rm Bair}\nolimits(E,\tau _{1})$, hence $\nu \in P_{\sigma }(E,\tau
_{1})$.

Suppose now that $(E, \tau_1)$ is a strongly measure-compact. Then there
exists a compact $L$ subset of $(E, \tau_1)$ such that $(\nu)^* (L) > 0 $.
By step~1, $L \in \mathop{\rm Bair}\nolimits(E, \tau_1)$, this implies that $%
\nu (L) = \nu^*(L)>0$, this is impossible by step~2. $\blacksquare $

Note that if $X$ a Lindel\"{o}f space, then $X$ is measure-compact \cite{Var}%
. In \cite{Mor1}-\cite{Mor2} one shows that there exists a measure-compact
space which is not strongly measure-compact. The following remark gives
another type of example.

\begin{remark}
\label{nnn} There exists a Lindel\"{o}f space $Y$, separable paracompact
which is not strongly measure-compact.
\end{remark}

\noindent\textit{Proof}. It suffices to apply Proposition~\ref{nh} to $E=%
\left] 0,1\right[ $, $\nu $ is the Lebesgue measure, hence $Y=(\left] 0,1%
\right[, \tau_1)$ is not strongly measure-compact. $\blacksquare $

Let $L$ be a Hausdorff compact space. In \cite{Tal(2)} one shows under the
Martin's Axiom that $(C(L),\tau _{p})$ is universally measurable. In the
following Lemma, one shows that $(C(L),\tau _{p})$ is universally measurable
without Martin's Axioms if $L$ is a Rosenthal compac set. \bigskip

\begin{proposition}
\label{mar} Let $L$ be a Rosenthal compact set. Then $(C(L),\tau _{p})$ is
universally measurable.
\end{proposition}

\noindent\textit{Proof}. We shall show Lemma~\ref{mar} if $C(L)$ is the real
continuous functions space. The proof is similar when one replaces $\mathbb{R%
}$ by $\mathbb{C}$. By \cite[th.~3.2]{Sch(2)} it suffices to show that every
normal probability measure on $(C(L), \tau_p)$ is a Radon measure. Let $\nu $
a normal probability measure on $(C(L), \tau_p)$.

Denote, for every $n\in \mathbb{N}^{\ast }$, $B_{n}=\left\{ g\in C(L);\Vert
g\Vert _{C(L)}\leq n\right\} $ and $\nu _{n}$ the measure defined by $\nu
_{n}(C)=\nu (C\cap B_{n})$, $C\in \mathop{\rm Bor}\nolimits(C(L),\tau _{p})$%
. Since $\nu _{n}(C)\rightarrow _{n\rightarrow +\infty }\nu (C)$ and $B_{n}$
is $\tau _{p}$ closed, it is enough to show that $\nu _{n}$ is a Radon
measure for all $n\geq 1$.

\emph{Step 1:} Let us show that for every $n\geq 1$, $\nu _{n}$ is a normal
measure.

Let $n\in \mathbb{N}^{\ast }$ and let $(F_{\alpha })$ be a net of closed
subsets of $(C(L),\tau _{p})$ which is decreasing to a closed set $F$ of $%
(C(L),\tau _{p})$. Since $B_{n}$ is $\tau _{p}$-closed, the net $(F_{\alpha
}\cap B_{n})$ is decreasing to $B_{n}\cap F$ and $\nu $ is normal, it
follows that

$\underset{\alpha }{\lim }\nu _{n}(F_{\alpha })=\underset{\alpha }{\lim }\nu
(F_{\alpha }\cap B_{n})=\nu (F\cap B_{n}))=\nu ((\underset{i\in I}{\cap }%
F_{i})\cap B_{n})=\nu _{n}(F)$. $\blacksquare $

\emph{Step 2: }Pick $n\in \mathbb{N}$. Let us show that $\nu _{n}$ is a
Radon measure.

We can suppose that $\mathop{\rm supp}\nolimits(\nu _{n})$ is included in
the closed unit ball of $(C(L),\left\Vert .\right\Vert )$. Denote $V$ the
vector subspace spanned by $L$ in the dual of $(C(L),\Vert .\Vert )$. For
every $v\in V$ let $H_{v}$ the hyperplane defined by $\nu $ in $C(L)$, i.e. 
\begin{equation*}
H_{v}=\left\{ f\in C(L);(f,v)=v(f)=0\right\} .
\end{equation*}%
Consider the set $I^{\prime }$ formed of all elements $v$ of $V$ such that $%
\nu _{n}(H_{v})=1$ and $E_{\nu _{n}}$ is the intersection of $H_{v}$ when $%
v\in I^{\prime }$ (note that $0\in I^{\prime }$).

Let us show that $\nu _{n}(E_{\nu _{n}})=1$. Let $\widetilde{I}=\left\{ S;S%
\text{ is a finite subset of }I^{\prime }\right\} $ and $H_{S}=\underset{%
v\in S}{\cap }H_{v}$, $S\in \widetilde{I}$. It is clear that $(H_{S})_{S\in 
\widetilde{I}}$ is a net of decreasing closed subsets of $(C(L),\tau _{p})$;
since $\nu _{n}$ is normal, 
\begin{equation*}
\nu _{n}(E_{\nu _{n}})=\nu _{n}(\underset{v\in I^{\prime }}{\cap }H_{v})=\nu
_{n}(\underset{S\in \widetilde{I}}{\cap }H_{S})=\inf_{S\in \widetilde{I}}\nu
_{n}(H_{S})=1,
\end{equation*}%
because $\nu _{n}(H_{S})=1$, for every $S\in \widetilde{I}$. For $y\in L$
put 
\begin{equation*}
\left[ y\right] =\left\{ x\in L;(f,y)=(f,x),\forall f\in E_{\nu _{n}}\right\}
\end{equation*}%
and $L_{\nu _{n}}=\left\{ \left[ y\right] ;y\in L\right\} $ ($L_{\nu _{n}}$
as a topological subspace of $\left[ -1,+1\right] ^{E_{\nu _{n}}}$, where $%
\left[ y\right] (f)=(f,y)$, $f\in E_{\nu _{n}})$. We shall show that $L_{\nu
_{n}}$ is compact. It is enough to show that $L_{\nu _{n}}$ is closed in $%
\left[ -1,+1\right] ^{E_{\nu _{n}}}$. Let $(\left[ y_{\alpha }\right]
)_{\alpha \in I}$ be a net in $L_{\nu _{n}}$ such that $\left[ y_{\alpha }%
\right] \rightarrow z\in \left[ -1,+1\right] ^{E_{\nu _{n}}}$. Since $L$ is
compact, there exists $y\in L$ such that $y_{\alpha }\rightarrow y\in L$,
this implies that $\left[ y_{\alpha }\right] (f)=(f,y_{\alpha })\rightarrow
(f,y)=\left[ y\right] (f)$ for every $f\in E_{\nu _{n}}$. On the other hand, 
$\left[ y_{\alpha }\right] (f)=(f,y_{\alpha })\rightarrow z(f)$, it follows
that $z(f)=\left[ y\right] (f)$ for every $f\in E_{\nu _{n}}$. We conclude
that $z=\left[ y\right] $. \vskip2pt plus 1pt minus 1pt

\emph{Step 2-a):} Suppose that $\left[ y\right] =\left[ y^{\prime }\right] $%
, $\nu _{n}$-almost-everywhere, $(y,y^{\prime }\in L)$ and show that $\left[
y\right] =\left[ y^{\prime }\right] $ everywhere.

Indeed, for almost all $f\in C(L)$, $(y,f)=(y^{\prime },f)$. Thus ( $%
y-y^{\prime },f)=0$ for almost all $f\in C(L)$, i.e. the hyperplane $%
H_{y-y^{\prime }}$ is of measure one, hence $H_{y-y^{\prime }}$ contains $%
E_{\nu _{n}}$. We deduce that for every $f\in E_{\nu _{n}}$, $%
(f,y)=(f,y^{\prime })$. $\blacksquare $

We define the map $\Sigma :L\rightarrow L_{\nu _{n}}$, by $\Sigma (y)=\left[
y\right] $, $y\in L$ (it is clear that $\Sigma $ is continuous) and the map $%
U:L_{\nu _{n}}\rightarrow L^{1}(E_{\nu _{n}},\nu _{n})$ by $U(\left[ x\right]
)(f)=f(x)$, $x\in L$, $f\in E_{\nu _{n}}$. Note that $U$ is into by step 3.

\emph{Step 2-b): }Sow that $U$ is continuous. For that let $F$ be a closed
subset of $L^{1}(E_{\nu _{n}},\nu _{n})$ and let $\left[ x\right] $ be a
point in the closure of $U^{-1}(F)$. There exists a net $(\left[ x_{\alpha }%
\right] )_{i\in I}$ in $U^{-1}(F)$ such that $\left[ x_{\alpha }\right]
\rightarrow \left[ x\right] $. Since $L$ is compact, there exists $x^{\prime
}\in L$ such that $x_{\alpha }\rightarrow x^{\prime }$ in $L$, by continuity
of $\Sigma $, $\left[ x_{\alpha }\right] \rightarrow \left[ x^{\prime }%
\right] $ in $L_{\nu _{n}}$. It follows that $[x^{\prime }]=[x]$.

On the other hand, by \cite{Bour-Frem-Tal} $L$ is angelic, hence there is a
countable subsequence $(x_{i_{k}})_{k\geq 0}$ such that $x_{i_{k}}%
\rightarrow _{k\rightarrow +\infty }x^{\prime }$. By the continuity of $%
\Sigma $ we have $[x_{i_{k}}]\rightarrow _{k\rightarrow +\infty }[x^{\prime
}]=[x]$ in $L_{\nu _{n}}$. By dominated convergence theorem, $U(\left[
x_{i_{n}}\right] )\rightarrow _{n\rightarrow +\infty }U(\left[ x\right] )$
in $L^{1}(E_{\nu _{n}},\nu _{n})$. Since $U(\left[ x_{i_{k}}\right] )\in F$,
for every $k\in \mathbb{N}$, $U([x])\in F$. It follows that $U$ is
continuous. $\blacksquare $

Thus $L_{\nu _{n}}$ is metrisable, because $U$ is a homeomorphism on its
range. Therefore, $(C(L_{\nu _{n}}),\left\Vert .\right\Vert )$ is a Polish
space, and $(C(L_{\nu _{n}}),\tau _{p})$ is defined by a separating family
of seminorms, by \cite[TG.~9, 68-6, th.~14]{Bourb} $\mathop{\rm Bor}%
\nolimits(C(L_{\nu _{n}}),\left\Vert .\right\Vert )=\mathop{\rm Bor}%
\nolimits(C(L_{\nu _{n}}),\tau _{p})$. Since $(C(L_{\nu _{n}}),\left\Vert
.\right\Vert )$ is a Radon space \cite[chap.~II, th.~9]{Sch(1)} and $(E_{\nu
_{n}},\tau _{p})$ is a closed subspace of $(C(L_{\nu _{n}}),\tau _{p})$, $%
(E_{\nu _{n}},\tau _{p})$ is a Radon space. We deduce that $\nu _{n}$ is a
Radon measure. $\blacksquare $ \smallskip

For every locally topological vector space $X$, we denote by $%
\mathop{\rm
Cyl}\nolimits(X)$ the $\sigma$-algebra generated by continuous linear forms
on~$X$. \bigskip

Recall that $(X,weak)$ coincides with $\mathop{\rm Cyl}\nolimits(X)$ \cite[%
th.~2.3]{Ed}.

\begin{lemma}
\label{mpn} Let $L$ be a Rosenthal compact set. Then $(C(L),\tau _{p})$ is
measure-compact if and only if $(C(L),weak)$ is measure-compact.
\end{lemma}

\noindent \textit{Proof}. Show first that $\mathop{\rm Cyl}%
\nolimits((C(L),weak))=\mathop{\rm Cyl}\nolimits((C(L),\tau _{p}))$. (It
follows by \cite[th.~2.3]{Ed} that 
\begin{equation}
\mathop{\rm Bair}\nolimits(C(L),weak))=\mathop{\rm Bair}\nolimits(C(L),\tau
_{p})\text{)}.  \label{bo}
\end{equation}

Indeed, it is clear that $\mathop{\rm Cyl}\nolimits((C(L),\tau _{p}))\subset %
\mathop{\rm Cyl}\nolimits((C(L),weak))$. Let us show the converse inclusion.
Let $y^{\ast }\in M^{+}(L)$. By \cite{Gode}, $M^{+}(L)$ is a Rosenthal
compact set, hence by \cite{Bour-Frem-Tal} it is angelic , we deduce that
there exists a sequence $(y_{n}^{\ast })_{n\geq 0}$ of finite support on $L$
such that $y_{n}^{\ast }\underset{n\rightarrow +\infty }{\rightarrow }%
y^{\ast }$ for the $w^{\ast }$-topology. For each $n\in \mathbb{N}$, $%
y_{n}^{\ast }$ is measurable with respect to $\mathop{\rm Cyl}%
\nolimits((C(L),\tau _{p}))$, this implies that $y^{\ast }$ is measurable
with respect to $\mathop{\rm Cyl}\nolimits((C(L),\tau _{p}))$. If $y^{\ast
}\in B_{(C(L))^{\ast }}$, there exists $\mu _{1,}\mu _{2},\mu _{3},\mu
_{4}\in M^{+}(L)$ such that $y^{\ast }=\mu _{1}-\mu _{2}-i(\mu _{3}-\mu _{4})
$, hence $y^{\ast }$ is measurable with respect to $\mathop{\rm Cyl}%
\nolimits((C(L),\tau _{p}))$, thus $\mathop{\rm
Cyl}\nolimits((C(L),weak))\subset \mathop{\rm Cyl}\nolimits((C(L),\tau
_{p})) $.

Suppose now that $(C(L),\tau _{p})$ is measure-compact. Let $\nu $ be a
measure defined on $\mathop{\rm Cyl}\nolimits((C(L),weak))$. By \cite[th.~2]%
{Sun}, $\nu $ extends to a normal measure $\widetilde{\nu }$ on $(C(L),\tau
_{p})$. On the other hand, by Proposition~\ref{mar} $(C(L),\tau _{p})$ is
universally measurable, hence $\widetilde{\nu }$ is a Radon measure on $%
(C(L),\tau _{p})$. Let $\varepsilon >0$ and let $C\in \mathop{\rm Bor}%
\nolimits(C(L),\tau _{p})$. There exists a $\tau _{p}$-compact set $H\subset
C$ such that $\widetilde{\nu }(H)>\widetilde{\nu }(C)-\varepsilon $. For
every integer $n\in \mathbb{N}^{\ast }$, consider $B_{n}=\bigl\{g\in
C(L);\Vert g\Vert _{C(L)}\leq n\bigr\}$. Since $\widetilde{\nu }(H\cap
B_{n})\rightarrow _{n\rightarrow +\infty }\widetilde{\nu }(H)$, there is $%
n_{0}\in \mathbb{N}^{\ast }$ such that $\widetilde{\nu }(H\cap B_{n_{0}})>%
\widetilde{\nu }(C)-\varepsilon $.

Note that $H\cap B_{n_{0}}$ is an uniformly bounded $\tau _{p}$-compact
subset of~$H$, hence by \cite{Groth} $H\cap B_{n_{0}}$ is a weakly compact
subset. By \cite[chap.~I, th.~16]{Sch(1)}, one obtains that $\widetilde{\nu }
$ extends to a Radon-measure $\nu ^{\prime }$ on $(C(L),weak)$. In
particular $\nu $ extends to a normal measure on $(C(L),weak)$. Note that,
if $f$ is a bounded continuous function on $(C(L),weak)$, then we have 
\begin{equation*}
\int_{C(L)}f(x)\,d\widetilde{\nu }(x)=\int_{C(L)}f(x)\,d\nu ^{\prime
}(x)=\int_{C(L)}f(x)\,d\nu (x).
\end{equation*}%
Let now $(f_{\alpha })_{\alpha \in I}$ be a decreasing net of bounded
continuous functions on $(C(L),weak)$ such that $f_{\alpha }\rightarrow 0$
for the pointwise topology. By \cite[th.~2]{Sun} $\mathop{\displaystyle \int}%
\limits_{C(L)}f_{\alpha }(x)\,d\nu (x)=\mathop{\displaystyle \int}%
\limits_{C(L)}f_{\alpha }(x)\,d\widetilde{\nu }(x)\rightarrow 0$. It follows
that $(C(L),weak)$ is measure-compact.

Conversely, suppose that $(C(L),weak)$ is measure-compact. Since, by above $%
\mathop{\rm Bair}\nolimits(C(L),weak)=\mathop{\rm Bair}\nolimits(C(L),\tau
_{p})$, it follows that $(C(L),\tau _{p})$ is measure compact. $\blacksquare 
$

\begin{proposition}
\label{knn} Suppose that $K$ is separable. Then ($C(K),\tau _{p})$ is not
measure-compact.
\end{proposition}

\noindent \textit{Proof}. It is clear that for every $\theta \in
K=B_{1}^{{}}\cup B_{2}^{{}}\cup B_{3}^{{}}$, the map $t\rightarrow \psi
(\theta ,t)$ is measurable. \vskip2pt plus 1pt minus 1pt

By an argument similar to that of Proposition~\ref{ai}, we show that for
every $\xi \in C(K)^{\ast }$, the map $t\in E\rightarrow \xi (f^{t})$ is
measurable. \vskip2pt plus 1pt minus 1pt

Suppose now that ($C(K),\tau _{p})$ is measure-compact. Lemma~\ref{mpn}
shows us that $(C(K),weak)$ is measure-compact. Since the map $t\rightarrow
f_{t}$ is scalarly measurable and $(C(K),weak)$ is measure-compact, by
Theorem of \cite{Ed}, there exists a strongly measurable function $g$ with
values in $C(K)$ such that for every $\xi \in C(K)^{\ast }$ $\xi (f^{t})=\xi
(g(t))$, for almost all $t\in E$. \vskip2pt plus 1pt minus 1pt

Let $(\theta _{n})_{n\geq 0}$ be a dense sequence in $K$. For every $n\in 
\mathbb{N}$ there exists a measurable subset $D_{n}$ of $E$ such that $\mu
(D_{n})=0$ and $\theta _{n}(f^{t})=\theta _{n}(g(t))$, for all $t\in
(D_{n})^{c}$. Denote $D=\underset{n\geq 0}{\cup }D_{n}$. Since $\theta
_{n}(f^{t})=\theta _{n}(g(t))$, for all $n\in \mathbb{N}$ and all $t\in
D^{c} $, $\theta (f^{t})=\theta (g(t))$, for all $\theta \in K$ and all $%
t\in D^{c} $. It follows that $f^{t}=g(t)$, for almost all $t\in E$, this is
impossible, because $\left\Vert f_{t}-f_{t^{\prime }}\right\Vert _{L^{\infty
}(E,\mu )}\geq 1$, if $t\neq t^{\prime }$ and $g$ has almost values in a
separable subspace of $C(K)$. Thus $(C(K),\tau _{p})$ is not
measure-compact. $\blacksquare $

\section{The case of a separable connected abelian group verifying the
condition~(**)}

In this part, $E$ is a separable connected abelian group verifying the
condition~(**) and $E$ is locally compact. \bigskip

By Remark~\ref{gd}, $E$ has a countable basis. As $E$ is locally compact, by 
\cite[chap.~II, th.~6]{Sch(1)}, $E$ is a Polish space.

Note that $-\infty, +\infty \notin E$. For $g\in C(K)\subset L^{\infty }(E,
\mu)$, $\theta \in K$ and $t\in E$, one defines $g_t\in C(K)$ by $%
g_t(u)=g(u-t)$, for almost all $u\in E$ and $\theta_t\in K$ by $%
\theta_t(r)=\theta (r_t)$, $r\in C(K)$.

Let $f$ \ be the characteristic function of $\left\{ t\in E;t>0\right\} $.

\begin{lemma}
\label{ab} For every $\theta \in B_{2}^{{}}$, $h(\theta _{u})=h(\theta )-u$,
for all $u\in E$.
\end{lemma}

\noindent \textit{Proof}. Since $E$ satisfies the condition~$(\ast \ast )$, $%
f_{t}=f^{t}$, for every $t\in E$. Let $\theta \in B_{2}^{{}}$ and let $%
t,u\in E$ such that $\psi (\theta _{u},t)=1$. Hence we have $\psi (\theta
_{u},t)=\theta _{u}(f_{t})=\theta \left[ (f_{t})_{u}\right] =\theta
(f_{t+u})=\psi (\theta ,t+u)=1$. It follows that $t+u\leq h(\theta )$. The
condition~$(\ast \ast )$ implies that $t\leq h(\theta )-u$. By Lemma~\ref{kt}%
, one has $h(\theta _{u})=\sup \left\{ t;\psi (\theta _{u},t)=1\right\} \leq
h(\theta )-u$.

Conversely, let $t\in E$ such that $\psi (\theta ,t)=1$. For every $u\in E$,
we have $\psi (\theta _{u},t-u)=\theta _{u}(f_{t-u})=\theta
((f_{t-u})_{u})=\theta (f_{t})=\psi (\theta ,t)=1$, hence $h(\theta
_{u})\geq t-u$, since $E$ satisfies the condition~$(\ast \ast $), $h(\theta
_{u})+u\geq t$. By Lemma~\ref{kt}, we obtain that $h(\theta _{u})+u\geq
h(\theta )$, i.e. $h(\theta _{u})\geq h(\theta )-u$, we deduce that $%
h(\theta _{u})=h(\theta )-u$. $\blacksquare $

\begin{remark}
\label{mp} In the previous lemma, one can replace $B_{2}^{{}}$ by $%
B_{3}^{{}} $ (using Lemma~\ref{oo}).
\end{remark}

\begin{remark}
\label{dq} $B_{2}^{{}}$ is notempty.
\end{remark}

\noindent \textit{Proof}. Indeed, suppose that $B_{2}^{{}}=\emptyset $.
Proposition~\ref{bv} tells us that $B_{3}^{{}}=\emptyset $, hence $%
K=B_{1}^{{}}$ which is metrisable (by Proposition~\ref{gf}-ii), It is
impossible, because $(E,\mu )$ satisfies the condition~(\ref{x}). $%
\blacksquare $

\begin{lemma}
\label{hmo} For every $\theta \in K\setminus B_{1}^{{}}$: 
\begin{equation*}
\sup \left\{ -t\in E;\psi (\theta ,t)=0\right\} =-\inf \left\{ t\in E;\psi
(\theta ,t)=0\right\} .
\end{equation*}
\end{lemma}

\noindent \textit{Proof}. Consider $\theta \in K\setminus B_{1}^{{}}$. Put 
\begin{equation*}
\beta =-\inf \left\{ t\in E;\psi (\theta ,t)=0\right\} .
\end{equation*}

Let $t\in E$ such that $\psi (\theta ,t)=0$. By definition of $-\beta $ we
have $-\beta \leq t$. The condition~(**) implies that $-\beta +(\beta
-t)\leq t+(\beta -t)$, hence $-t\leq \beta $, to complete the proof, it
suffices to show that for every $z\in E$ such that $-t\leq z$, and if $\psi
(\theta ,t)=0$, then $\beta \leq z$. Let a such $z\in E$. Since $E$
satisfies the condition~(**), $t\geq -z$, for every $t\in E$ such that $\psi
(\theta ,t)=0$. We deduce that inf$\left\{ t\in E;\psi (\theta ,t)=0\right\}
\geq -z$, i.e. $\beta \leq z$. $\blacksquare $

For every $\theta \in K$, one defines $\overset{\vee }{\theta }\in K$, by $%
\overset{\vee }{\theta }(g)=\theta (\overset{\vee }{g})$, where $\overset{%
\vee }{g}(x)=g(-x)$, for almost $x\in E$, $g\in C(K)\subset L^{\infty }(E,m)$%
.

\begin{proposition}
\label{bv} Let $\theta \in K\setminus B_{1}^{{}}$. Then:

i) $h(\theta )=-h(\overset{\vee }{\theta })$.

ii) $\varphi \in B_{2}^{{}}$, if and only if $\overset{\vee }{\varphi }\in
B_{3}^{{}}$.

iii) $\overset{\vee}{(\theta^{\prime})} = \theta^{\prime\prime}$.
\end{proposition}

\noindent\textit{Proof}.

\textit{i)}. The measure $\mu =m$ is invariant by translation, then $\mu
(\{t\})=0$, for every $t\in E$. Thus we have 
\begin{align*}
& \mskip25mu\mu \Bigl(\bigl\{y\in E;\ \overset{\vee }{(f_{-t})}(y)=0\bigr\}%
\Bigr) \\
& =\mu (\left\{ y\in E;f_{-t}(-y)=0\right\} ) \\
& =\mu (\left\{ y\in E;f(t-y)=0\right\} )=\mu (\left\{ y\in E;\text{ }y\geq
t\right\} ) \\
& =\mu (\left\{ y\in E;y>t\right\} )=\mu (\left\{ y\in E;f_{t}(y)=1\right\} )
\\
& =\mu (\left\{ y\in E;e(y)-f_{t}(y)=0\right\} ).
\end{align*}%
By a similar argument we show that 
\begin{equation*}
\mu \bigl(\bigl\{y\in E;\overset{\vee }{(f_{-t}})(y)=1\bigr\}\bigr)=\mu
(\{y\in E;e(y)-f_{t}(y)=1\}).
\end{equation*}%
Thus for every $t\in E$ 
\begin{equation}
\overset{\vee }{(f_{-t})}=e-f_{t}.  \label{xp}
\end{equation}%
Let $\theta \in K\setminus B_{1}$ and $t\in E$. By (\ref{xp}), 
\begin{equation*}
\psi (\overset{\vee }{\theta },-t)=\overset{\vee }{\theta }(f_{-t})=\theta {%
\overset{\vee }{(f_{-t})}=\theta (e)-\psi (\theta ,t)}=1-\psi (\theta ,t),
\end{equation*}%
\textit{i.~e.} 
\begin{equation}
\psi (\overset{\vee }{\theta },-t)=0\text{ if and only if }\psi (\theta
,t)=1.  \label{n}
\end{equation}%
Using (\ref{n}) and Lemma~\ref{kt}, one obtains 
\begin{align*}
h(\theta )& =\sup \{t\in E;\psi (\theta ,t)=1\}=\sup \bigl\{t\in E;\psi (%
\overset{\vee }{\theta },-t)=0\bigr\} \\
& =\sup \bigl\{-t\in E;\psi (\overset{\vee }{\theta },t)=0\bigr\}.
\end{align*}%
\vskip2pt plus 1pt minus 4pt \vskip2pt plus 1pt minus 4pt

On the other hand, by Lemma~\ref{hmo} and Lemma~\ref{oo} 
\begin{eqnarray*}
&&\sup \bigl\{-t\in E;\psi (\overset{\vee }{\theta },t)=0\bigr\} \\
&=&-\inf \bigl\{t\in E;\psi (\overset{\vee }{\theta },t)=0\bigr\}=-h(\overset%
{\vee }{\theta }).
\end{eqnarray*}%
Therefore $h(\theta )=-h(\overset{\vee }{\theta })$. $\blacksquare $ \vskip%
2pt plus 1pt minus 4pt

\textit{ii)}. By (\ref{n}$)$ and (i), $\psi (\varphi ,h(\varphi ))=1$ if and
only if $\psi (\overset{\vee }{\varphi },h(\overset{\vee }{\varphi }))=\psi (%
\overset{\vee }{\varphi },-h(\varphi ))=0$. Thus $\varphi \in B_{2}^{{}}$ if
and only if $\overset{\vee }{\varphi }\in B_{3}^{{}}$. $\blacksquare $ \vskip%
2pt plus 1pt minus 4pt

\textit{iii)}. Obviously by (\ref{xp}). $\blacksquare $ \vskip 2pt plus 1pt
minus 4pt

\begin{lemma}
\label{ww} Let $\theta \in B_{j}^{{}}$ and $t\in E$. Then $\theta _{t}\in
B_{j}^{{}}$, $j\in \{2,3\}$.
\end{lemma}

\vskip 2pt plus 1pt minus 4pt

\noindent \textit{Proof}. Suppose that $\theta \in B_{2}^{{}}$. By Lemma~\ref%
{ab} we have $\psi (\theta _{t},h(\theta _{t}))=\psi (\theta _{t},h(\theta
)-t)=\psi (\theta ,h(\theta ))=1$, hence $\theta _{t}\in B_{2}^{{}}$. By a
similar argument, one shows that $\theta _{t}\in B_{3}^{{}}$, if $\theta \in
B_{3}^{{}}$. $\blacksquare $ \vskip2pt plus 1pt minus 4pt

\begin{proposition}
\label{lv} For $\theta \in B_{2}^{{}}$, the map $U_{\theta }^{\prime
}:E\rightarrow K\times K$, $t\rightarrow (\theta _{t},\overset{\vee }{%
(\theta )_{t}})$ is not Borel.
\end{proposition}

\vskip 2pt plus 1pt minus 4pt

\noindent \textit{Proof}. Suppose that there exists $\theta \in B_{2}^{{}}$
such that $U_{\theta }^{\prime }$ is a Borel function. Choose a subset $H$
of $E$ which is not in the Borel $\sigma $-algebra of $E$. We define the
open subset $V$ of $K\times K$ by 
\begin{align*}
V=\bigcup_{u\in H}\Bigl[& \left\{ \varphi \in K;\psi (\varphi ,-u+h(\theta
))=1\right\} \\
& \times \bigl\{\varphi \in K;\psi (\varphi ,-u+h(\overset{\vee }{\theta }%
))=0\bigr\}\Bigr].
\end{align*}%
\vskip2pt plus 1pt minus 4pt \noindent By Lemma~\ref{ww}, for every $t\in E$%
, $\theta _{t}\in B_{2}^{{}}$ and $(\overset{\vee }{\theta })_{t}\in
B_{3}^{{}}$ (because $\overset{\vee }{\theta }\in B_{3}^{{}}$, by
Proposition~\ref{bv}). Thus we have 
\begin{align*}
& \mskip25mu(U_{\theta }^{\prime })^{-1}(V) \\
& =\bigcup_{u\in H}\bigl\{t\in E;\psi (\theta _{t},-u+h(\theta ))=1\text{
and }\psi ((\overset{\vee }{\theta })_{t},-u+h(\overset{\vee }{\theta }))=0%
\bigr\} \\
& =\bigcup_{u\in H}\bigl\{t\in E;h(\theta _{t})\geq -u+h(\theta )\text{ and }%
h((\overset{\vee }{\theta })_{t})\leq -u+h(\overset{\vee }{\theta })\bigr\}
\\
& =\bigcup_{u\in H}\bigl\{t\in E;h(\theta )-t\geq -u+h(\theta )\text{ and }h(%
\overset{\vee }{\theta })-t\leq -u+h(\overset{\vee }{\theta })\bigr\} \\
& =\bigcup_{u\in H}\{t\in E;-u\leq -t\leq -u\}=\bigcup_{u\in H}\{u\}=H,
\end{align*}%
because $h(\theta _{t})=h(\theta )-t$ and $h((\overset{\vee }{\theta )}%
_{t})=h(\overset{\vee }{\theta })-t$, for every $t\in E$, by Lemma~\ref{ab}.
This implies that $H$ is a Borel subset of $E$, which is impossible. $%
\blacksquare $

\begin{corollary}
\label{zu} The Borel $\sigma $-algebra of the product space $K\times K$
strictly contains $\mathop{\rm Bor}\nolimits(K)\otimes \mathop{\rm Bor}%
\nolimits(K)$.
\end{corollary}

\noindent \textit{Proof}. Consider the maps $\delta _{1}:t\in E\rightarrow
\theta _{t}\in B_{2}^{{}}\subset K$, $\delta _{2}:t\in E\rightarrow (\overset%
{\vee }{\theta })_{t}\in B_{3}^{{}}\subset K$.

Let us show that $\delta _{1},\delta _{2}$ are Borel functions. Let $u,t\in
E $. It is clear that the subsets of the form 
\begin{eqnarray}
&&(\delta _{1})^{-1}\left[ \left\{ \varphi \in K;\psi (\varphi ,u)=0\right\} %
\right] ,  \label{bn} \\
&&(\delta _{1})^{-1}\left[ \left\{ \varphi \in K;\psi (\varphi ,t)=1\right\} %
\right] ,  \notag \\
&&(\delta _{2}^{{}})^{-1}\left[ \left\{ \varphi \in K;\psi (\varphi
,u)=0\right\} \right] ,  \notag \\
&&(\delta _{2})^{-1}\left[ \left\{ \varphi \in K;\psi (\varphi ,t)=1\right\} %
\right] ,  \notag
\end{eqnarray}%
are Borel subsets of $E$ (note for example that 
\begin{eqnarray*}
(\delta _{1})^{-1}\left[ \left\{ \varphi \in K;\psi (\varphi ,u)=0\right\} %
\right] &=&\left\{ t\in E;\psi (\theta _{t},u)=0\right\} \\
=\left\{ t\in E;h(\theta _{t})=h(\theta )-t<u\right\} &=&\left\{ t\in
E;t>h(\theta )-u\right\} ).
\end{eqnarray*}%
On the other hand by Proposition~\ref{wb}, $K$ is a strongly Lindel\"{o}f
space, hence $\delta _{1}$ and $\delta _{2}$ are Borel functions, we deduce
that for every $S\in \mathop{\rm Bor}\nolimits(K)\otimes \mathop{\rm Bor}%
\nolimits(K)$, $(U_{\theta }^{\prime })^{-1}(S)\in \mathop{\rm Bor}%
\nolimits(E)$.

Suppose now that $\mathop{\rm Bor}\nolimits(K\times K)=\mathop{\rm Bor}%
\nolimits(K)\otimes \mathop{\rm Bor}\nolimits(K)$. By above, $U_{\theta
}^{\prime }$ is a Borel function, this is impossible by Proposition~\ref{lv}%
. $\blacksquare $

L. Lindenstrauss and C. Stegall \cite{Lin-Steg} showed that there exists a
Banach space $X$ which does not contain $\ell^\infty $ isomorphically and a
map $\sigma :\left\{ 0,1\right\}^{\mathbb{N}}\rightarrow X$ scalarly
measurable, but $\sigma $ is not weakly equivalent to a measurable function.
In the following proposition, one gives an other type of examples. In the
fact, the considered function (in Proposition~\ref{zp}) is locally
integrable in the sense of Riemann. \vskip 2pt plus 1pt minus 4pt

\begin{proposition}
\label{zp} There exists a Banach space $X$ and a function $g:E\rightarrow X$
such that

I) $X$ does not contain $\ell ^{\infty }$ isomorphically.

II) $g$ is scalarly measurable.

III) $g$ cannot be weakly equivalent to a strongly measurable function.

IV) $g$ is locally integrable in the sense of Riemann.
\end{proposition}

\vskip 2pt plus 1pt minus 4pt \vskip 2pt plus 1pt minus 4pt

\noindent\textit{Proof}.

I). Consider $X=C(K)$. Since for all $\xi \in B_{C(K)^{\ast }}$ there exist $%
\xi _{1},\xi _{2},\xi _{3},\xi _{4}\in M^{+}(K)$ such that $\xi =\xi
_{1}-\xi _{2}+i(\xi _{3}-\xi _{4})$, by \cite{Bour-Frem-Tal}-\cite{Gode} $%
B_{C(K)^{\ast }}$ is angelic, hence $X$ does not contain $\ell ^{\infty }$
isomorphically, because by \cite{Od-Ros} the closed unit ball of $(\ell
^{\infty })^{\ast }$ is not angelic. $\blacksquare $ \vskip2pt plus 1pt
minus 4pt

II). Let $g : E\rightarrow (C(K) \subset L^\infty(E)$ the function defined
by $g(t)=f_t$ for almost all $t\in E$. By Proposition~\ref{ai}, $g$ is
scalarly measurable. $\blacksquare $ \vskip 2pt plus 1pt minus 4pt

III). In the proof of Proposition~\ref{knn} we saw that $g$ cannot be weakly
equivalent to a strongly measurable function. $\blacksquare $ \vskip 2pt
plus 1pt minus 4pt

IV). Let $a,b\in E$ such that $a<b$. Note that \vskip2pt plus 1pt minus 4pt 
\begin{equation*}
\int_{a}^{b}\chi _{\lbrack t,+\infty \lbrack }(\omega )\,dt=\left\{ 
\begin{array}{l}
(a\vee \omega )-a,\text{ si }\omega \leq b \\ 
b-a,\text{ si }\omega >b.%
\end{array}%
\right.
\end{equation*}%
\vskip2pt plus 1pt minus 4pt \noindent Consider $t_{0}=a$, $t_{1},\ldots
,t_{n}=b\in \lbrack a,b]$ such that $t_{0}<t_{1}<\ldots <t_{n}$. For $%
t_{k(\omega )}\leq \omega <t_{k(\omega )+1}$, we have $\sum\limits_{k=1}^{n}%
\chi _{\left[ t_{k},+\infty \right[ }(\omega
)(t_{k}-t_{k-1})=\sum\limits_{k=1}^{k(\omega )}\chi _{\lbrack t_{k},+\infty
\lbrack }(\omega )(t_{k}-t_{k-1})=\sum\limits_{k=1}^{k(\omega
)}(t_{k}-t_{k-1})=t_{k(\omega )}-a$. One deduces that \vskip2pt plus 1pt
minus 4pt 
\begin{equation*}
\sum\limits_{k=1}^{n}\chi _{\lbrack t_{k},+\infty \lbrack }(\omega
)(t_{k}-t_{k-1})=\left\{ 
\begin{array}{l}
t_{k(\omega )}-a,\text{ si }t_{k(\omega )}\leq \omega <t_{k(\omega )+1} \\ 
0,\text{ si }\omega <a \\ 
\sum\limits_{k=1}^{n}(t_{k}-t_{k-1})=b-a,\text{ si }\omega >b.%
\end{array}%
\right.
\end{equation*}%
\vskip2pt plus 1pt minus 4pt \noindent Thus 
\begin{align*}
& \mskip40mu\Bigl\|\sum\limits_{k=1}^{n}\chi _{\lbrack t_{k},+\infty \lbrack
}(.)(t_{k}-t_{k-1})-\int\limits_{a}^{b}\chi _{\lbrack t,+\infty \lbrack
}(.)\,dt\Bigr\|_{L^{\infty }(E)} \\
& \leq \underset{\omega \in \lbrack a,b]}{\sup }|\omega -t_{k(\omega )}|\leq 
\underset{1\leq k\leq n}{\sup }\left\vert t_{k}-t_{k-1}\right\vert .
\end{align*}%
It follows that $g$ is integrable in the sense of Riemann on $[a,b]$. $%
\blacksquare $ \vskip2pt plus 1pt minus 4pt

Let $M$ be the subset of $C((K)$ formed by even functions with values in $%
\left\{ 0,1,2\right\} $ ($g$ is even function on $K$ if $g=\overset{\vee }{g}%
).$

\begin{theorem}
\label{vh} There exists a Hausdorff topology $\tau ^{\prime }$ on $C(K)$, a $%
C^{\ast }$-subalgebra $Y_{1}$ of $C(K)$ and a Banach subspace $Y_{2}$ of $%
C(K)$ such that

\textbf{1}$)$ $\tau_p$ is finer than $\tau^{\prime}$.

\textbf{2) } $(C(K),\tau ^{\prime })$ is not universally measurable.

\textbf{3) }$Y_{1}$ and $Y_{2}$ are $\tau _{p}$ closed and $(C(K),\tau
_{p})=(Y_{1},\tau _{p})\oplus (Y_{2},\tau _{p})$.

\textbf{4) } $(Y_{1},\tau _{p})$ is isomorphic to $(Y_{2},\tau _{p})$.

\textbf{5) } $(Y_{j},\tau _{p})=(Y_{j},\tau ^{\prime })$, $j\in \{1,2\}$.

\textbf{6) } $(Y_{1},\tau _{p})$ contains a closed discrete uncountable
subset which is contain in $M$.

\textbf{7) } $\mathop{\rm Bor}\nolimits(C(K),\tau ^{\prime })\otimes %
\mathop{\rm Bor}\nolimits(C(K),\tau ^{\prime })\neq \mathop{\rm Bor}%
\nolimits(C(K)\times C(K),\tau ^{\prime }\times \tau ^{\prime })$.

\textbf{8) } The projection $P:(C(K)),\tau ^{\prime })\rightarrow
(Y_{1},\tau ^{\prime })$ is not a Borel function.

\textbf{9) } The set $\{H:E\rightarrow (C(K)),\tau ^{\prime });H\text{ is a
Borel function}\}$ is not a vector subspace of $C(K)^{E}$.
\end{theorem}

\noindent \textit{Proof}.

1).

\emph{Step 1: }Let $\theta \in B_{2}^{{}}$. Let us show that $\left\{ \theta
_{t};t\in E\right\} =B_{2}^{{}}$.

First observe that $B_{2}^{{}}\neq \emptyset $ by Remark~\ref{dq} . Let $%
\varphi \in B_{2}^{{}}$. There exists $t\in E$ such that $h(\varphi
)=h(\theta _{t})$ (because $h(\theta _{t})=h(\theta )-t$ by Lemma~\ref{ab}).
By Lemma~\ref{ww}, $\theta _{t}\in B_{2}^{{}}$. Since $h_{2}^{{}}:B_{2}^{0}%
\rightarrow E$ is into, by Remark~\ref{cnn}, $\varphi =\theta _{t}$. $%
\blacksquare $\vskip2pt plus 1pt minus 4pt

\emph{Step 2 : }Show that there exists a sequence in $B_{2}^{{}}$ which is
dense in~$K$. \vskip2pt plus 1pt minus 4pt

Since $h_{2}^{{}}$ is onto, there is a sequence $(\theta _{n})_{n\geq 0}$ in 
$B_{2}^{{}}$ such that $(h(\theta _{n}))_{n\geq 0}$ is dense in $E$. \vskip%
2pt plus 1pt minus 4pt

We will show that $(\theta _{n})_{n\geq 0}$ is dense in $K$. For that, let $%
t\in E$. Denote $V_{t}=\left\{ \theta \in K;\psi (\theta ,t)=1\right\} $.
Choose $\theta \in K\setminus B_{1}^{{}}$ such that $t<h(\theta )$, there
exists $\theta _{m}\in B_{2}^{{}}$ satisfying $t<h(\theta _{m})$, because
the sequence $(h(\theta _{n}))_{n\geq 0}$ is dense in $E$. Since $\psi
(\theta _{m},h(\theta _{m}))=1$ and $t<h(\theta _{m})$, $\theta _{m}\in
V_{t} $. Let $u\in E$ and let $U_{u}=\left\{ \theta \in K;\psi (\theta
,u)=0\right\} $. Choose $\theta \in K\setminus B_{1}^{{}}$ such that $%
u>h(\theta )$, there is $\theta _{m^{\prime }}\in B_{2}^{{}}$ satisfying $%
u>h(\theta _{m^{\prime }})$. It follows that $\psi (\theta _{m}^{\prime
},u)=0$, because $\psi (\theta _{m}^{\prime },h(\theta _{m}^{\prime }))=1$.
Thus $\theta _{m^{\prime }}\in U_{u}$. \vskip2pt plus 1pt minus 4pt

Let $t,u\in E$ and let 
\begin{equation*}
\theta _{0}\in W_{t,u}=\{\theta \in K;\psi (\theta ,t)=1\text{ and }\psi
(\theta ,u)=0\}.
\end{equation*}%
We note that $t<u$. Since $E$ satisfies the condition~$(\ast )$ and $h$ is
onto, there exists $\theta _{1}\in K\setminus B_{1}^{{}}$ such that $%
t<h(\theta _{1})<u$. We deduce that there exists $\theta _{n}\in B_{2}^{{}}$
verifying $t<h(\theta _{n})<u$. Thus $\theta _{n}\in W_{t,u}$. \smallskip $%
\blacksquare $

Let $\tau ^{\prime }$ be the topology defined on $C(K)$ by the condition
that $u_{\alpha }\rightarrow u$ in $(C(K),\tau ^{\prime })$, if and only if $%
u_{\alpha }(\theta )\rightarrow u(\theta )$ for all $\theta \in
B_{2}^{{}}\cup B_{1}^{{}}$. Therefore $\tau _{p}$ is finer than $\tau
^{\prime }$, and by step~2 $\tau ^{\prime }$ is Hausdorff.

2). Put $D=\left\{ f_t; t\in E\right\} \subset (C(K), \tau^{\prime})$.
\smallskip

\emph{Step 1}: Show that $(D,\tau ^{\prime })$ is a strongly Lindel\"{o}f
space.

Pick $\theta \in B_{2}^{{}}$ and define the map $\pi :(D,\tau ^{\prime
})\rightarrow B_{2}^{{}}$, by $\pi (f_{t})=\theta _{t}$, $t\in E$. It
suffices to show that $\pi $ is a homeomorphism, because $B_{2}^{{}}$ is
strongly Lindel\"{o}f. \smallskip

Let us first show that $\pi $ is injective. For that, let $u,v\in E$ such
that $\theta _{u}=\theta _{v}$. By Lemma~\ref{ab}, $h(\theta )-u=h(\theta
)-v $, hence $u=v$. Let $u\in E$. Put $V_{u}=\left\{ \theta \in
B_{2}^{{}};\psi (\theta ,u)=0\right\} $ and $W_{u}=\left\{ \theta \in
B_{2}^{{}};\ \psi (\theta ,u)=1\right\} $. Thus we have $\pi
^{-1}(V_{u})=\left\{ f_{t};\psi (\theta _{t},u)=0\right\} =\left\{
f_{t};\theta _{t}(f_{u})=0\right\} =\left\{ f_{t};\theta
_{u}(f_{t})=0\right\} $ is an open subset of $(D,\tau ^{\prime })$. In the
same way, we have $\pi ^{-1}(W_{u})=\left\{ f_{t};\theta
_{u}(f_{t})=1\right\} $ is an open subset of ($D,\tau ^{\prime })$. We
deduce that $\pi $ is continuous.

We shall show that $\pi ^{-1}$ is continuous. Indeed, let $t_{1},t_{2}\in E$%
. Denote $W_{1}=\left\{ f_{t};f_{t}(\theta _{t_{1}})=1\right\} $ and $%
W_{2}=\left\{ f_{t};f_{t}(\theta _{t_{2}})=0\right\} $. It is enough to show
that $\pi (W_{j})$ is an open subset of $B_{2}^{{}}$, $j\in \left\{
1,2\right\} $, because the open sets of the previous form generated the
topology $\tau ^{\prime }$. \vskip2pt plus 1pt minus 1pt

Note that $\pi (W_{1})=\left\{ \theta _{t};f_{t}(\theta _{t_{1}})=1\right\}
=\left\{ \theta _{t};f_{t_{1}}(\theta _{t})=1\right\} $ and $\pi
(W_{2})=\left\{ \theta _{t};f_{t}(\theta _{t_{2}})=0\right\} =\left\{ \theta
_{t};f_{t_{2}}(\theta _{t})=0\right\} $. We deduce that $\pi (W_{1})$, $\pi
(W_{2})$ are open subsets of $B_{2}^{{}}$. \vskip2pt plus 1pt minus 1pt Thus 
$\pi $ is a homeomorphism. $\blacksquare $

Let $F_{1}$ be the closure of $D$ in $(C(K),\tau ^{\prime })$ and let $%
\sigma _{1}:(E,\tau _{0})\rightarrow F_{1}\subset (C(K),\tau ^{\prime })$
the map defined by $\sigma _{1}(t)=f_{t}$, $t\in E$. We will show that $%
\sigma _{1}$ is Borel. \vskip2pt plus 1pt minus 1pt

Indeed, let $O$ be an open subset of $(D,\tau ^{\prime })$. There exists a
sequence $(O_{i})_{i\in I}$ of $D$ such that $O=\underset{i\in I}{\cup }%
O_{i} $ and each $O_{i}$ is a finite intersection of open subsets, under the
form $W_{1},W_{2}$, where $W_{1}=\left\{ f_{t};f_{t}(\theta
_{t_{1}})=1\right\} $, $W_{2}=\left\{ f_{t};f_{t}(\theta _{t_{2}})=0\right\} 
$, $t_{1},t_{2}\in E$. Note that for every $\theta \in B_{2}^{{}}\cup
B_{1}^{{}}$ the map $t\in E\rightarrow f_{t}(\theta )=\psi (\theta ,t)$ is
Borel, hence $(\sigma _{1}^{-1})(O_{i})$ is Borel, for all $i\in I$. By
step~1, $(D,\tau ^{\prime })$ is strongly Lindel\"{o}f, one concludes that
there is a countable subset $I_{1}$ of $I$ such that $O=\underset{i\in I_{1}}%
{\cup }O_{i}$. It follows that $(\sigma _{1})_{{}}^{-1}(O)$ is a Borel
subset of $(E,\tau _{0})$. Finally if $O^{\prime }$ is an open subset of $%
F_{1}$, then $(\sigma _{1})^{-1}(O^{\prime })=(\sigma _{1})^{-1}(O^{\prime
}\cap D)$ and $O=O^{\prime }\cap D$ is an open subset of $(D,\tau ^{\prime
}) $. \vskip2pt plus 1pt minus 4pt

\emph{Step 2:} Show that every compact of $F_{1}$ is metrisable.

Note that if $g\in F_1$, $g(\theta) \in \{0, 1\}$.

Let $L$ be a compact subset of $F_{1}$ and let $(\theta _{n})_{n\geq 0}$ be
a sequence in $B_{2}^{{}}$ which is dense in $K$ (a such sequence exists by
step~2 of~1). Define the map $\xi :L\rightarrow \left\{ 0,1\right\} ^{%
\mathbb{N}}$ by $\xi (g)=(g(\theta _{n}))_{n\geq 0}$, $g\in L$. It is clear
that $\xi $ is continuous and injective, this implies that $L$ is
metrisable. $\blacksquare $\vskip2pt plus 1pt minus 4pt

Observe now that the topology $\tau ^{\prime }$ is defined by a family of
semi-norms which separates the points of $C(K)$, hence $(C(K),\tau ^{\prime
})$ is completely regular space. Since every subspace of a completely
regular space is completely regular, $(F_{1},\tau ^{\prime })$ is completely
regular. \vskip2pt plus 1pt minus 4pt

Suppose that $(C(K), \tau^{\prime})$ is universally measurable. \smallskip

\emph{Step 3:} Show that $F_{1}$ is universally measurable.

By \cite[th.~3.2]{Sch(2)}, it is enough to show that every normal
probability measure on $F_{1}$ is a Radon measure. Let $\nu $ be a normal
probability measure on $F_{1}$. Consider $\nu ^{\prime }$ the measure
defined on $(C(K),\tau ^{\prime })$ by $\nu ^{\prime }(B)=\nu (B\cap F_{1})$%
, where $B\in \mathop{\rm Bor}\nolimits(C(K),\tau ^{\prime })$. Let us show
that $\nu ^{\prime }$ is a normal measure on $(C(K),\tau ^{\prime })$. \vskip%
2pt plus 1pt minus 4pt

Let $i:F_{1}\rightarrow (C(K),\tau ^{\prime })$ be the canonical injection.
Since $i$ is continuous and $\nu $ is normal, by \cite{Sch(2)}, $i(\nu )=\nu
^{\prime }$ is a normal measure. On the other hand $(C(K),\tau ^{\prime })$
is universally measurable and $\nu ^{\prime }$ is normal, hence by \cite[%
th.~3.2]{Sch(2)}, $\nu ^{\prime }$ is a Radon measure. Now we let $B^{\prime
}$ be a Borel subset of $F_{1}$. Note that $B^{\prime }$ is a Borel subset
of $(C(K),\tau ^{\prime })$. \vskip2pt plus 1pt minus 4pt

Thus we have 
\begin{eqnarray*}
\nu (B^{\prime }) &=&\nu ^{\prime }(B^{\prime })=\sup \{\nu ^{\prime
}(L^{\prime });L^{\prime }\text{ compact in }B^{\prime }\} \\
&=&\sup \{\nu (L^{\prime });L^{\prime }\text{ compact in }B^{\prime }\}.
\end{eqnarray*}%
Therefore $\nu $ is a Radon measure. It follows that $F_{1}$ is universally
measurable. $\blacksquare $\vskip2pt plus 1pt minus 1pt

Let $L$ be a compact subset of $(E,\tau _{0})$ of strictly positive measure.
We can suppose that $\mu (L)=1$. Consider the map $\eta _{1}:L\rightarrow
F_{1}$, the restriction of $\sigma _{1}$ to $L$ and $\nu _{1}=\eta _{1}(\mu
_{L})$, where $\mu _{L}$ the measure defined on $L$ by 
\begin{equation*}
\mu _{L}(C)=\mu (C)=m(C),\quad C\in \mathop{\rm Bor}\nolimits(L)\subset %
\mathop{\rm Bor}\nolimits(E).
\end{equation*}%
\vskip2pt plus 1pt minus 1pt

\emph{Step 4:} Show that $\nu _{1}$ is a normal measure on~$F_{1}$.

Let $(U_{\alpha })$ be a net of open subsets in $F_{1}$ which is increasing
to an open subset of~$F_{1}$. Observe that for every Borel subset $B$ of~$%
F_{1}$, $\nu _{1}(B)=\mu (\{t\in L;f_{t}\in B\})=\mu (\{t\in L;f_{t}\in
B\cap D\})$.

By step~1 of 2), $(D,\tau ^{\prime })$ is strongly Lindel\"{o}f, hence there
exists a countable subset $I_{1}$ of $I$ such that $\underset{i\in I}{\cup }%
U_{i}\cap D=\underset{i\in I_{1}}{\cup }U_{i}\cap D$. Thus 
\begin{align*}
\nu _{1}(\underset{i\in I}{\cup }U_{i})& =\mu (\bigl\{t\in L;f_{t}\in (%
\underset{i\in I}{\cup }U_{i})\cap D\bigr\}) \\
& =\mu (\bigl\{t\in L;f_{t}\in \underset{i\in I}{\cup }U_{i}\cap D\bigr\}%
)=\mu (\bigl\{t\in L;f_{t}\in \underset{i\in I_{1}}{\cup }U_{i}\cap D\bigr\})
\\
& =\mu (\bigl\{t\in L;f_{t}\in \underset{i\in I_{1}}{(\cup }U_{i})\cap D%
\bigr\}) \\
& =\nu _{1}(\underset{i\in I_{1}}{\cup }U_{i}))=\sup_{i\in I_{1}}\nu
_{1}(U_{i})\leq \sup_{i\in I}\nu _{1}(U_{i}).
\end{align*}%
This implies that $\nu _{1}$ is a normal measure on $F_{1}$. $\blacksquare $%
\vskip2pt plus 1pt minus 1pt

By step 3 of 2), $F_{1}$ is universally measurable, since $\nu _{1}$ is a
normal probability measure on $F_{1}$, by \cite[th.~3.2]{Sch(2)}, $\nu _{1}$
is a Radon measure, we deduce that there exists a compact subset $G_{1}$ of $%
F_{1}$ such that $\nu _{1}(G_{1})>\frac{1}{2}$. We denote by $\eta _{2}$ the
restriction of $\sigma _{1}$ to ($\eta _{1}^{{}})^{-1}(G_{1})\subset L$. the
map $\eta _{2}$ is Borel with values in $G_{1}$ and $m((\eta
_{1}^{{}})^{-1}(G_{1}))=\nu _{1}(G_{1})>\frac{1}{2}$. The step~2 of~2) shows
us that $G_{1}$ is metrisable. By \cite{Frem}, there exists a compact subset 
$L_{1}$ of ($(\eta _{1})^{-1}(G_{1}),\tau _{0})$ such that $\mu (L_{1})>%
\frac{1}{3}$ and the restriction of $\sigma _{1}$ to $L_{1}$ is continuous. %
\vskip2pt plus 1pt minus 1pt

Since $\mu $ is diffuse on $E$, by an argument similar to that of
Proposition~\ref{ynm}, one constructs a strictly decreasing sequence $%
(t_{n})_{n\geq 0}$ in $L_{1}$. There exists $t_{0}\in L_{1}$ such that $t_{n}%
\underset{n\rightarrow +\infty }{\rightarrow }t_{0}$ (modulo subsequence).
On the other hand, the restriction of $\sigma _{1}$ to $L_{1}$ is
continuous, it follows that $f_{t_{n}}\underset{n\rightarrow +\infty }{%
\rightarrow }f_{t_{0}}$ in $(C(K),\tau ^{\prime })$. For $s=h(\theta )-t_{0}$%
, we have $f_{t_{n}}(\theta _{s})\underset{n\rightarrow +\infty }{%
\rightarrow }f_{t_{0}}(\theta _{s})=\psi (\theta _{s},t_{0})=\psi (\theta
,h(\theta ))=1$. Thus for $n$ large enough, $\psi (\theta _{s},t_{n})=1$,
hence $t_{n}\leq h(\theta _{s})=h(\theta )-s=t_{0}$. For seeing a
contradiction, it suffices to note that the sequence $(t_{n})_{n\geq 0}$ is
strictly decreasing. $\blacksquare $ \vskip2pt plus 1pt minus 1pt

3). Let $Y_{1}$ be the vector subspace formed by even functions in $%
C(K)\subset L^{\infty }(E)$ and $Y_{2}$ be the vector space formed by odd
functions in $C(K)\subset L^{\infty }(E)$. Observe that if $g\in C(K)$, $%
\overset{\vee }{g}\in C(K)$, because $\overset{\vee }{f}_{t}=e-f_{-t}$, $%
\forall t\in E$. \vskip2pt plus 1pt minus 1pt

Let us show that $Y_1$ is $\tau_p$ closed in $C(K)$. For that let $%
(r_\alpha)_{\alpha \in I}$ be a net in $Y_1$ such that $r_{\alpha
}\rightarrow r$ in $(C(K), \tau_p)$. We must prove that $r\in Y_1$. \vskip%
2pt plus 1pt minus 1pt

For $\alpha \in I$ and $\theta \in K$, we have $r_\alpha (\overset{\vee }{%
\theta })=\overset{\vee }{r}_\alpha (\theta)=r_\alpha (\theta)$, by passing
to the limit we obtain that $\overset{\vee }{r}(\theta)=r(\theta)$, hence $%
\overset{\vee }{r}=r\in Y_1$. By a similar argument, we show that $Y_2$ is $%
\tau_p$ closed in $C(K)$. \vskip2pt plus 1pt minus 1pt

On the other hand, for $g\in C(K)$, $g=(g+\overset{\vee }{g})/2\oplus (g-%
\overset{\vee }{g})/2$, this implies that $(C(K),\tau _{p})=(Y_{1},\tau
_{p})\oplus (Y_{2},\tau _{p})$ (note that the map $g\in C(K)\rightarrow 
\overset{\vee }{g}\in C(K)$ is continuous). $\blacksquare $ \vskip2pt plus
1pt minus 1pt

4). Define the operator $U:(Y_{1},\tau _{p})\rightarrow (Y_{2},\tau _{p})$,
by 
$U(g)=2gf-g$, $g\in Y_{1}$. Let us show that $U$ has values in $Y_{2}$. For
that let $t\in G.$ Observe that%
\begin{equation*}
(2gf-g)(-t)=2g(t)f(-t)-g(t),
\end{equation*}
but $f(-t)=e-f(t),$ hence%
\begin{eqnarray*}
(2gf-g)(-t) &=&2g(t)(e-f(t)-g(t) \\
&=&-\left[ 2g(t)f(t)-g(t)\right] .
\end{eqnarray*}

Since $2(2gf-g)f-(2gf-g)=g$, $U^{-1}(g)=2gf-g$, $g\in Y_{2}$. We observe
that $U,U^{-1}$ are continuous. Thus $U$ is an isomorphism. $\blacksquare $
\smallskip

5). Show that $(Y_{1},\tau ^{\prime })=(Y_{1},\tau _{p})$. For $\theta \in
B_{2}^{{}}\cup B_{1}^{{}}$ and $g\in Y_{1}$ we have $g(\overset{\vee }{%
\theta })=(\overset{\vee }{g})(\theta )=g(\theta )$. Since $\Bigl\{\overset{%
\vee }{\theta };\theta \in B_{2}^{{}}\cup B_{1}^{{}}\Bigr\}=B_{3}^{{}}\cup
B_{1}^{{}}$ by Proposition~\ref{bv}, it follows that $(Y_{1},\tau ^{\prime
})=(Y_{1},\tau _{p})$.

Le us show that $(Y_{2},\tau ^{\prime })=(Y_{2},\tau _{p})$. Let $(g_{\alpha
})_{\alpha \in I}$ be a net in $Y_{2}$ and let $g\in Y_{2}$ such that $%
g_{\alpha }(\theta )\rightarrow g(\theta )$, for every $\theta \in
B_{2}^{{}}\cup B_{1}^{{}}$. We must prove that $g_{\alpha }\rightarrow g$ in 
$(Y_{2},\tau _{p})$. Note that $U^{-1}(g_{\alpha })(\theta )=(2g_{\alpha
}\times f-g_{\alpha })(\theta )$ and $U^{-1}(g)(\theta )=(2g\times
f-g)(\theta )$, hence $U^{-1}(g_{\alpha })(\theta )\rightarrow
U^{-1}(g)(\theta )$, for every $\theta \in B_{2}^{{}}\cup B_{1}^{{}}$. Since 
$(Y_{1},\tau ^{\prime })=(Y_{1},\tau _{p})$, $U^{-1}(g_{\alpha
})(.)\rightarrow U^{-1}(g)(.)$ in $(Y_{1},\tau _{p})$, we deduce that $%
g_{\alpha }\rightarrow g$ in $(Y_{2},\tau _{p})$, because $U$ is a
homeomorphism. Thus $(Y_{2},\tau ^{\prime })=(Y_{2},\tau _{p})$. $%
\blacksquare $

6). For every $t\in E$ denote $w_{t}=f_{t}+\overset{\vee }{(f_{t})}%
=f_{t}+1-f_{-t}$ and $V=\left\{ \omega _{t};t\in E\right\} $. We shall show
that $(V,\tau _{p})$ is discrete and closed in $(Y_{1},\tau _{p})$. Indeed,
let $t\in E;$ choose $\theta _{1}\in B_{2}^{{}}$ and $\theta _{2}\in
B_{3}^{{}}$ such that $h(\theta _{1})=h(\theta _{2})=t$. Denote 
\begin{equation*}
W=\Bigl\{\omega _{s};\bigl|\omega _{s}(\theta _{j})-\omega _{t}(\theta _{j})%
\bigr|<\frac{1}{2},j\in \{1,2\}\Bigr\}.
\end{equation*}%
We observe that $W$ is an open subset of $(V,\tau _{p})$. \medskip

Let us show that $W=\left\{ \omega_t\right\} $. \vskip2pt plus 1pt minus 1pt

\emph{Case 1,} $t<0$: Let $\omega _{s}\in W$. Note that $\omega _{t}(\theta
_{1})=1+1-0=2$ and $\omega _{t}(\theta _{2})=0+1-0=1$, hence $\omega
_{s}(\theta _{1})=2$ and $\omega _{s}(\theta _{2})=1$, this implies that $%
\psi (\theta _{1},s)=1,$ $\psi (\theta _{1},-s)=0$ and $\psi (\theta
_{2},s)=\psi (\theta _{2},-s)$.

The case $\psi (\theta _{2},s)=\psi (\theta _{2},-s)=1$ is excluded, because
it will lead that $t>0$, which is impossible. We have then $\psi (\theta
_{1},s)=1$ and $\psi (\theta _{2},s)=\psi (\theta _{2},-s)=0$. It follows
that $s\leq h(\theta _{1})=t\leq s$. It means that $t=s$.

\emph{Case 2,} $t>0$: Let $\omega _{s}\in W$. Note that $\omega _{t}(\theta
_{1})=1+1-1=1$ and $\omega _{t}(\theta _{2})=0+1-1=0$, hence $\omega
_{s}(\theta _{1})=\psi (\theta _{1},s)+1-\psi (\theta _{1},-s)=1$ and $%
\omega _{s}(\theta _{2})=\psi (\theta _{2},s)+1-\psi (\theta _{2},-s)=0$.
Therefore $\psi (\theta _{1},s)=\psi (\theta _{1},-s)$ and $\psi (\theta
_{2},s)=0$ (note that $\psi (\theta _{2},-s)=1)$. The case $\psi (\theta
_{1},s)=\psi (\theta _{1},-s)=0$ leads that $t<0$, which is impossible,
hence $\psi (\theta _{1},s)=1$ and $\psi (\theta _{2},s)=0$. It follows that 
$s=t$.

\emph{Case 3,} $t=0$: Let $\omega _{s}\in W$. It is clear that $\omega
_{0}(\theta _{1})=\omega _{0}(\theta _{2})=1$, which implies that $\psi
(\theta _{1},s)=\psi (\theta _{1},-s)$ and $\psi (\theta _{2},s)=\psi
(\theta _{2},-s)$. Note that the case $\psi (\theta _{1},s)=\psi (\theta
_{1},-s)=0$ and the case $\psi (\theta _{2},s)=\psi (\theta _{2},-s)=1$ are
excluded. Thus $\psi (\theta _{1},s)=1$ and $\psi (\theta _{2},s)=0$. This
implies that $s=0$. Therefore $(V,\tau _{p})$ is a discrete space in $%
(Y_{1},\tau _{p})$. \vskip2pt plus 1pt minus 1pt

It remains to show that $V$ is a closed subset in $(Y_{1},\tau _{p})$. \vskip%
2pt plus 1pt minus 1pt

Let $(\omega _{t_{\alpha }})_{\alpha \in I}$ be a net in $V$ such that $%
\omega _{t_{\alpha }}\rightarrow g\in Y_{1}$. Choose a fixed point $\theta $
in $B_{2}^{{}}$. The sequence $(\theta _{t_{\alpha }})_{\alpha \in I}$ is a
net in $K$, hence, there exists $\theta _{0}\in K$ such that $\theta
_{t_{\alpha }}\rightarrow \theta _{0}\in K$ (with respected to an
ultrafilter $\mathcal{U}$ on $I$). \medskip

\emph{Case 1,} $t_{\alpha }\rightarrow +\infty $: Let $\varphi \in
K\setminus \left\{ \theta ^{\prime },\theta ^{\prime \prime }\right\} $.
There exist $v\in E^{+}$ and $u\in E^{-}$ such that $-u>v$, $\psi (\varphi
,u)=1$ and $\psi (\varphi ,v)=0$. Put $W=\left\{ t\in E;t>v\text{ and }%
-t<u\right\} $. We observe that $S=\left\{ \alpha \in I;t_{\alpha }\in
W\right\} $ belongs to $\mathcal{U}$, hence for all $\alpha \in S$, $%
f_{t_{\alpha }}(\varphi )=\psi (\varphi ,t_{\alpha })=0$ and $f_{-t_{\alpha
}}(\varphi )=\psi (\varphi ,-t_{\alpha })=1$, it follows that $\omega
_{t_{\alpha }}(\varphi )=0+1-1=0$, for every $\alpha \in S$. We conclude
that $g(\varphi )=0$. It is obvious that $\omega _{t_{\alpha }}(\theta
^{\prime })=$ $\omega _{t_{\alpha }}(\theta ^{\prime \prime })=1$, which
implies that $g(\theta ^{\prime })=g(\theta ^{\prime \prime })=1$.
Therefore, we have $\left\{ \theta ^{\prime },\theta ^{\prime \prime
}\right\} =g^{-1}(\left\{ 1\right\} )$, hence this subset is closed and
open, because $g$ has values in $\left\{ 0,1,2\right\} $ (note that $%
K=\left\{ \theta ^{\prime },\theta ^{\prime \prime }\right\} \cup
g^{-1}(\left\{ 0,2\right\} ))$. Since $\left\{ \theta ^{\prime }\right\}
=\left\{ \theta ^{\prime },\theta ^{\prime \prime }\right\} \setminus
\{\theta ^{\prime \prime }\}$, $\left\{ \theta ^{\prime }\right\} $ is an
open subset, this means that there is $t_{0}\in E$ such that $\left\{ \theta
^{\prime }\right\} =\left\{ \theta \in K;\psi (\theta ,t_{0})=1\right\} $,
this is impossible, because the set $\left\{ \theta \in K;\psi (\theta
,t_{0})=1\right\} $ contains the set $\left\{ \theta \in K;h(\theta
)<t_{0}\right\} $. \vskip2pt plus 1pt minus 4pt

\emph{Case 1,} $t_{\alpha }\rightarrow -\infty $: By an argument similar to
the previous case, we show that this case can be excluded, because $%
g(\varphi )\in \left\{ 0,1,2\right\} $, for every $\varphi \in K$. \vskip2pt
plus 1pt minus 4pt

Since the cases~1 and~2 are excluded and $E$ is a locally compact, there
exists $t_{0}\in E$ such that $t_{\alpha }\rightarrow t_{0}\in E$. By Lemma~%
\ref{ab}, $h(\theta _{t_{\alpha }})=h(\theta )-t_{\alpha }$, which implies
that $t_{\alpha }\rightarrow h(\theta )-h(\theta _{0})=t_{0}$ (because $h$
is continuous by Proposition \ref{xh}). \vskip2pt plus 1pt minus 4pt

Show that $g=\omega _{t_{0}}$. Pick $\varphi \in K\setminus B_{1}^{{}}$. %
\vskip2pt plus 1pt minus 4pt \goodbreak

\emph{Case 1,} $t_0>0$. \vskip 2pt plus 1pt minus 4pt

\emph{Case 1-a), }$-t_{0}<h(\varphi )<t_{0}$: Let $W=\left\{ t\in
E;-t<h(\varphi )<t\right\} $. $W$ is an open neighborhood of $t_{0}$, hence 
\begin{equation*}
S=\left\{ \alpha \in I;-t_{\alpha }<h(\varphi )<t_{\alpha }\right\} \in 
\mathcal{U}.
\end{equation*}%
For every $\alpha \in S$ we have $f_{t_{\alpha }}(\varphi )=\psi (\varphi
,t_{\alpha })=\psi (\varphi ,t_{0})=0$ and $\psi (\varphi ,-t_{\alpha
})=\psi (\varphi ,-t_{0})=1$, this implies that $\omega _{t_{\alpha
}}(\varphi )=\omega _{t_{0}}(\varphi )=0+1-1=0$ for all $\alpha \in S$. Thus 
$\omega _{t_{\alpha }}(\varphi )\rightarrow \omega _{t_{0}}(\varphi
)=g(\varphi )$.

\emph{Case 1-b), }$h(\varphi )<-t_{0}$. Put 
\begin{equation*}
W=\left\{ t\in E;h(\varphi )<-t<0\ \right\} .
\end{equation*}%
$W$ is an open neighborhood of $t_{0}$, hence $S=\left\{ \alpha \in
I;h(\varphi )<-t_{\alpha }<0\right\} $ belongs to $\mathcal{U}$. It is clear
that $\omega _{t_{\alpha }}(\varphi )=\omega _{t_{0}}(\varphi )=0+1-0=1$
(because $t_{\alpha }>-t_{\alpha }$ and $\psi (\varphi ,-t_{\alpha })=0,$
hence $\psi (\varphi ,t_{\alpha })=0)$. We deduce that $\omega _{t}(\varphi
)\rightarrow \omega _{t_{0}}(\varphi )=1=g(\varphi )$.

\emph{Case 1-c),} $h(\varphi )>t_{0}$: By an argument similar to case 1-b),
one shows that $\omega _{t_{\alpha }}(\varphi )\rightarrow \omega
_{t_{0}}(\varphi )=g(\varphi )$.

\emph{Case 2,} $t_{0}<0$: The case $t_{0}<h(\varphi )<-t_{0}$, $h(\varphi
)<-t_{0}$ and the case $h(\varphi )>t_{0}$ will be are treated by arguments
similar to the previous cases.

\emph{Case 3)}, $h(\varphi )=t_{0}$, $t_{0}\in E$: Without losing the
generality, we can suppose that $\varphi \in B_{2}^{{}}$. In the proof of $%
1) $, we showed that $B_{2}^{{}}$ is dense in $K$. Hence there exists a
sequence $(\varphi _{n})_{n\geq 0}$ in $B_{2}^{{}}$ such that $\varphi
_{n}\rightarrow _{n\rightarrow +\infty }\varphi $. We can choose the
sequence $(\varphi _{n})_{n\geq 0}$ in $B_{2}^{{}}\setminus \{\varphi \}$,
because $K$ does not contain isolated points. Let us show that for $n$ large
enough $t_{0}<h(\varphi _{n})$ (note that the restriction of $h$ to $%
B_{2}^{{}}$ is injective, hence $h(\varphi _{n})\neq t_{0})$.

Indeed, suppose that for every $n\in \mathbb{N}$, there exists $m_{n}\in 
\mathbb{N}$ such that $t_{0}>h(\varphi _{m_{n}})$. Since $\psi (\varphi
_{m_{n}},h(\varphi _{m_{n}}))=1$, $\psi (\varphi _{m_{n}},t_{0})=0$. By
passing to the limit, we obtain that $\psi (\varphi ,t_{0})=\psi (\varphi
,h(\varphi ))=0$, which is impossible, hence there exists $n_{0}\in \mathbb{N%
}$ such that $t_{0}<h(\varphi _{n})$ for all $n\geq n_{0}$.

Applying the case 1-c, to $\varphi =\varphi_n$, we obtain $\omega_{t_\alpha
}(\varphi_n) \rightarrow \omega_{t_0}(\varphi_n) = g(\varphi_n)$ for every $%
n\geq n_0$. It follows that $\omega _{t_0}(\varphi)=g(\varphi)$, because $g$
is continuous.

\emph{Case 4,} $h(\varphi )>t_{0}=0$: Put 
\begin{equation*}
W=\left\{ t\in E;h(\varphi )>t\text{ and }-h(\varphi )<t\right\} .
\end{equation*}%
$W$ is a open neighborhood of $t_{0}=0$, hence 
\begin{equation*}
S=\left\{ \alpha \in I;h(\varphi )>t_{\alpha }\text{ and }-h(\varphi
)<t_{\alpha }\right\} \in \mathcal{U}.
\end{equation*}%
If $\alpha \in S$, $\omega _{t_{\alpha }}(\varphi )=\omega _{0}(\varphi
)=1+1-1=1$. We conclude that $\omega _{t_{\alpha }}(\varphi )=1\rightarrow
\omega _{0}(\varphi )=1=g(\varphi )$.

The case $h(\varphi )<t_{0}=0$ is treated by a similar argument. If $\varphi
\in B_{1}^{{}}$, it is clear that $\omega _{t_{\alpha }}(\varphi
)\rightarrow \omega _{t}(\varphi )$. Thus $V$ is a closed subset in $%
(Y_{1},\tau _{p})$. $\blacksquare $ \vskip2pt plus 1pt minus 1pt

7). We will show that the map $\sigma _{2}:t\in E\rightarrow \overset{\vee }{%
(f_{t})}\in (C(K),\tau ^{\prime })$ is Borel. \vskip2pt plus 1pt minus 1pt

Note that for every $t\in E$, $\overset{\vee }{(f_{t})}=e-f_{-t}$. Consider
the map $\Delta :E\rightarrow E$ defined by $\Delta (t)=-t$, $t\in E$, $%
\Delta $ is a Borel function. Thus $\sigma _{2}=e-\sigma _{1}\circ \Delta $, 
$\sigma _{2}$ is a Borel function, because in the proof of 2) we showed that 
$\sigma _{1}:t\in E\rightarrow f_{t}$ is Borel. \vskip2pt plus 1pt minus 1pt

Consider the map $\phi_1:t\in E\rightarrow (f_t, \overset{\vee }{(f_t)})\in
C(K)\times C(K)$. It is obvious that if $W \in \mathop{\rm Bor}%
\nolimits(C(K), \tau^{\prime}) \otimes \mathop{\rm Bor}\nolimits(C(K),
\tau^{\prime})$, then $\phi_1^{-1}(W)$ is a Borel subset of $E$.

Suppose now that $\mathop{\rm Bor}\nolimits(C(K), \tau^{\prime}) \otimes %
\mathop{\rm Bor}\nolimits(C(K), \tau^{\prime}) = \mathop{\rm Bor}%
\nolimits(C(K) \times C(K), \tau^{\prime}\times \tau^{\prime})$. Therefore
the map $\phi_1 : E \rightarrow (C(K) \times C(K), \tau^{\prime}\times
\tau^{\prime})$, $t \rightarrow (f_t, \overset{\vee }{(f_t)})$ is Borel.
\medskip

Define the map $\phi _{2}:(C(K)\times C(K),\tau ^{\prime }\times \tau
^{\prime })\rightarrow (C(K),\tau ^{\prime })$, by $\phi _{2}(g,u)=g+u$, $%
(g,u)\in C(K)\times C(K)$. Observe that $\phi _{2}$ is continuous, hence it
is Borel. It follows that $\phi _{2}\circ \phi _{1}$ is Borel. We deduce
that the map $t\in E\rightarrow w_{t}=f_{t}+\overset{\vee }{(f_{t})}\in
(Y_{1},\tau ^{\prime })=(Y_{1},\tau _{p})$ is Borel. But $\mathop{\rm
Bor}\nolimits(C(K),\tau _{p})=\mathop{\rm Bor}\nolimits(C(K),\left\Vert
.\right\Vert )$ by Corollary~\ref{mr}, hence $\mathop{\rm Bor}%
\nolimits(Y_{1},\tau _{p})=\mathop{\rm Bor}\nolimits(Y_{1},\Vert .\Vert )$,
which implies that the map $t\in (E,\tau _{0})\rightarrow \omega _{t}\in
Y_{1}$ is strongly measurable, this is impossible, because in the proof of
6), we saw that $(V,\tau _{p})$ is uncountable discrete space, this means
that the map $t\in E\rightarrow \omega _{t}$ cannot be almost everywhere
valued in a separable subspace of~$C(K)$. $\blacksquare $

8). Suppose that the projection $P:(C(K),\tau ^{\prime })\rightarrow
(Y_{1},\tau ^{\prime })=(Y,\tau _{p}$) is Borel. Since the map $\sigma
_{1}:t\in E\rightarrow f_{t}$ $\in (C(K),\tau ^{\prime })$ is Borel, the map 
$t\in E\rightarrow P\circ \sigma _{1}(t)=Pf_{t}=w_{t}/2$ $\in (C(K),\tau
^{\prime })$ is Borel, but in the proof of 7) we saw that this is
impossible. Thus $P$ is not a Borel function. $\blacksquare $

9). Put $H_{1}=\sigma _{1}:t\in E\rightarrow f_{t}\in (C(K),\tau ^{\prime })$
and $H_{2}=\sigma _{2}:t\in E\rightarrow \overset{\vee }{(f_{t})}\in
(C(K),\tau ^{\prime })$. We proved previously that $H_{1}$, $H_{2}$ are
Borel and in the end of the proof of $8)$ we proved that $H_{1}+H_{2}$ is
not Borel. $\blacksquare $

\begin{definition}
\label{bh} Let $(X,\tau )$ be a topological Hausdorff space. $(X,\tau )$ is
called a $K$-analytic set \cite{Bress-Si} if $(X,\tau )$ is a continuous
image of a set belonging to the family $K_{\sigma \delta }$.
\end{definition}

By \cite[th.~2.7.1]{Ro-Ja}, every $K$-analytic set is a Lindel\"{o}f space.
On the other hand if $L$ is separable and compact such that $(C(L),\tau
_{p}) $ is a Lindel\"{o}f space, then $L$ is metrisable \cite{Ar} (under $%
MA+\urcorner CH$).

\begin{corollary}
\label{jj} $(C(K), \tau^{\prime})$ is not $K$-analytic.
\end{corollary}

\noindent\textit{Proof}. Suppose that $(C(K), \tau^{\prime})$ is $K$%
-analytic. By \cite{Bress-Si} $(C(K), \tau^{\prime})$ is universally
measurable, which is impossible by Theorem~\ref{vh} 2). $\blacksquare $

\begin{corollary}
\label{ma} There exists two subspaces $Z_{1}$, $Z_{2}$ of $(C(K),\tau _{p})$%
, $\tau _{p}$-closed such that $Z_{1}\cap Z_{2}=\left\{ 0\right\} $, $%
(Z_{j},\tau ^{\prime })$ is universally measurable, $j\in \left\{
1,2\right\} $ and $(Z_{1},\tau ^{\prime })$ is isomorphic to $(Z_{2},\tau
^{\prime })$, but $(Z_{1}\oplus Z_{2},\tau ^{\prime })$ is not universally
measurable.
\end{corollary}

\noindent \textit{Proof}. Let $Z_{j}=Y_{j}$, $j\in \left\{ 1,2\right\} $. By
Theorem~\ref{vh}, one has $(Z_{j},\tau ^{\prime })=(Z_{j},\tau _{p})$, $%
Z_{j} $ is $\tau _{p}$ closed in $(C(K),\tau _{p})$, $j\in \left\{
1,2\right\} $ and $(Z_{1},\tau ^{\prime })$ is isomorphic to $(Z_{2},\tau
^{\prime })$. On the other hand, Lemma~\ref{mar} shows us that $(C(K),\tau
_{p})$ is universally measurable, hence $(Z_{j},\tau ^{\prime })=(Z_{j},\tau
_{p})$ is universally measurable (because a closed subset of an universally
measurable set is universally measurable). Since $(Z_{1}\oplus Z_{2},\tau
^{\prime })=(C(K),\tau ^{\prime })$, by Theorem~\ref{vh} 2), $(Z_{1}\oplus
Z_{2},\tau ^{\prime })$ is not universally measurable. $\blacksquare $
\medskip

\begin{remark}
\label{ka} By Corollary~\ref{mr}, $\mathop{\rm Bor}\nolimits(C(K), \tau_p)
\otimes \mathop{\rm Bor}\nolimits(C(K), \tau_p) $ is equal to $%
\mathop{\rm
Bor}\nolimits(C(K), \left\| .\right\| ) \otimes \mathop{\rm Bor}%
\nolimits(C(K), \left\| .\right\| )$.
\end{remark}

On the other hand by \cite[th.~3]{Tal(1)} (with the continuous hypthothesis)
one has 
\begin{equation*}
\mathop{\rm Bor}\nolimits(C(K),\Vert .\Vert )\otimes \mathop{\rm Bor}%
\nolimits(C(K),\Vert .\Vert )=\mathop{\rm Bor}\nolimits(C(K)\times
C(K),\Vert .\Vert ).
\end{equation*}%
Thus $\mathop{\rm Bor}\nolimits(C(K)\times C(K),\tau _{p}\times \tau
_{p})\subset \mathop{\rm Bor}\nolimits(C(K)\times C(K),\Vert .\Vert )=$

\noindent $=\mathop{\rm Bor}\nolimits(C(K),\Vert .\Vert )\otimes 
\mathop{\rm
Bor}\nolimits(C(K),\Vert .\Vert )=\mathop{\rm Bor}\nolimits(C(K),\tau
_{p})\otimes \mathop{\rm Bor}\nolimits(C(K),\tau _{p})$. But $%
\mathop{\rm
Bor}\nolimits(C(K),\tau _{p})\otimes \mathop{\rm Bor}\nolimits(C(K),\tau
_{p})\subset \mathop{\rm Bor}\nolimits(C(K)\times C(K),\tau _{p}\times \tau
_{p})$, one deduces that 
\begin{equation*}
\mathop{\rm Bor}\nolimits(C(K),\tau _{p})\otimes \mathop{\rm Bor}%
\nolimits((C(K),\tau _{p})=\mathop{\rm Bor}\nolimits(C(K)\times C(K),\tau
_{p}\times \tau _{p}).
\end{equation*}

\textbf{Acknowledgement. }I would like to thank Professors Bernard Maurey
and Gilles Godefroy for the time they devoted to me during the preparation
of this work.

\end{document}